\pdfoutput=1
\documentclass[11pt]{article}
\usepackage{amsfonts,amssymb,amsmath,amsthm}
\usepackage{thm-restate}

\usepackage{bm}
\usepackage{palatino}
\usepackage{mathpazo}
\usepackage{inconsolata}
\usepackage{caption}
\usepackage[hidelinks]{hyperref}
\hypersetup{
	pdftitle={Enumerating Pattern-Avoiding Involutions using Combinatorial Exploration},
	pdfauthor={Christian Bean, Anthony J. Guttmann, Jay Pantone},
}
\usepackage{xcolor}
\usepackage{tikz}
\usetikzlibrary{calc}
\usetikzlibrary{patterns}
\usetikzlibrary{shapes.geometric}
\usepackage{tabularray}
\usepackage{colortbl}
\usepackage{enumitem}
\usepackage{tilings}

\begingroup
    \makeatletter
    \@for\theoremstyle:=definition,remark,plain\do{%
        \expandafter\g@addto@macro\csname th@\theoremstyle\endcsname{%
            \addtolength\thm@preskip\parskip
            }%
        }
\endgroup

\usepackage[margin=1in]{geometry}

\theoremstyle{definition}
\newtheorem{theorem}{Theorem}
\newtheorem{proposition}[theorem]{Proposition}

\newtheorem{conjecture}[theorem]{Conjecture}

\newcommand*\patchAmsMathEnvironmentForLineno[1]{%
  \expandafter\let\csname old#1\expandafter\endcsname\csname #1\endcsname
  \expandafter\let\csname oldend#1\expandafter\endcsname\csname end#1\endcsname
  \renewenvironment{#1}%
     {\linenomath\csname old#1\endcsname}%
     {\csname oldend#1\endcsname\endlinenomath}}%
\newcommand*\patchBothAmsMathEnvironmentsForLineno[1]{%
  \patchAmsMathEnvironmentForLineno{#1}%
  \patchAmsMathEnvironmentForLineno{#1*}}%
\AtBeginDocument{%
\patchBothAmsMathEnvironmentsForLineno{equation}%
\patchBothAmsMathEnvironmentsForLineno{align}%
\patchBothAmsMathEnvironmentsForLineno{flalign}%
\patchBothAmsMathEnvironmentsForLineno{alignat}%
\patchBothAmsMathEnvironmentsForLineno{gather}%
\patchBothAmsMathEnvironmentsForLineno{multline}%
}

\usepackage{parskip}
\usepackage{lineno}
\usepackage[small]{titlesec}
\usepackage{authblk}

\usepackage[square,comma,numbers,sort&compress]{natbib}

\usepackage{color}
\definecolor{todocolor}{RGB}{205,235,139}
\definecolor{todo-idea}{RGB}{120,180,255}
\definecolor{todo-error}{RGB}{208,31,60}
\definecolor{todo-question}{RGB}{255,255,136}

\usepackage[colorinlistoftodos, color=todocolor, textsize=small]{todonotes}

\newcommand{\CC}{\mathcal{C}}
\newcommand{\TT}{\mathcal{T}}
\newcommand{\OO}{\mathcal{O}}
\newcommand{\RR}{\mathcal{R}}
\newcommand{\GG}{\mathcal{G}}
\newcommand{\ds}{\displaystyle}
\DeclareMathOperator{\Av}{Av}
\DeclareMathOperator{\I}{\Av^I}
\DeclareMathOperator{\Grid}{Grid}
\newcommand{\GridI}{\Grid^I}
\DeclareMathOperator{\gr}{gr}
\DeclareMathOperator{\ugr}{\overline{\gr}}
\DeclareMathOperator{\lgr}{\underline{\gr}}
\newcommand{\N}{\mathbb{N}}
\newcommand{\mc}{\mathcal}

\tikzset{treeedge/.style={out=-90, in=90}}

\renewenvironment{abstract}{
	\begin{list}{}%
	{\setlength{\rightmargin}{1in}%
	\setlength{\leftmargin}{1in}}%
	\item[]\ignorespaces\begin{small}}%
	{\end{small}\unskip\end{list}%
}

\title{Enumerating Pattern-Avoiding Involutions using Combinatorial Exploration}
\author[1]{Christian Bean}
\author[2]{Anthony J.\ Guttmann}
\author[3]{Jay Pantone}
\affil[1]{School of Computer Science and Mathematics, Keele University. \texttt{c.n.bean@keele.ac.uk}}
\affil[2]{School of Mathematics and Statistics, The University of Melbourne. \texttt{guttmann@unimelb.edu.au}}
\affil[3]{Department of Mathematical and Statistical Sciences, Marquette University. \texttt{jay.pantone@marquette.edu}}

\date{}

\begin{document}
\maketitle

\begin{abstract}
	The enumeration of pattern-avoiding permutations has been a popular area of study over the past several decades, but comparatively little attention has been given to the topic of pattern-avoiding involutions. In this paper, we derive the algebraic generating functions of two Wilf-equivalence classes of involutions avoiding a single pattern of length $4$, $\I(2431)$ and $\I(3421)$. We then adapt the Mosaic method, a fast counting algorithm for permutations, to count involutions and apply it to substantially extend the known initial terms of the counting sequences for the remaining two Wilf-equivalence classes avoiding a pattern of length $4$, $\I(1324)$ and $\I(4231)$. Based on these extended sequences, we empirically analyze the asymptotic behavior of the counting sequences of these two classes.
\end{abstract}


\section{Introduction}

A \emph{permutation} of length $n$ is a bijection from $\{1, 2, \ldots, n\}$ to itself, written in one-line notation as $\pi = \pi(1)\pi(2)\cdots\pi(n)$. The \emph{inverse} of $\pi$ is the permutation $\pi^{-1}$ satisfying $\pi^{-1}(i) = j$ if and only if $\pi(j) = i$. An \emph{involution} is a permutation satisfying $\pi = \pi^{-1}$.

We say that a permutation $\pi$ of length $n$ \emph{contains} the permutation $\sigma$ of length $k$ as a \emph{pattern}, and write $\sigma \leq \pi$, if there exist indices $1 \leq i_1 < i_2 < \cdots < i_k \leq n$ such that the subsequence $\pi(i_1)\pi(i_2)\cdots\pi(i_k)$ is order-isomorphic to $\sigma$; that is, $\pi(i_a) < \pi(i_b)$ if and only if $\sigma(a) < \sigma(b)$. If $\pi$ does not contain $\sigma$, we say that $\pi$ \emph{avoids} $\sigma$. For example, the permutation $24513$ contains $132$ because the subsequence $\pi(1)\pi(3)\pi(5) = 253$ is order-isomorphic to $132$.

\begin{figure}[ht]
	\centering
	\begin{tikzpicture}[scale=0.45]
		\draw[gray!30] (0.5,0.5) grid (9.5,9.5);
		\foreach \i in {1,...,9} {
			\node[below] at (\i, 0.5) {\small $\i$};
			\node[left]  at (0.5, \i) {\small $\i$};
		}
		\foreach \x/\y in {1/3, 2/7, 3/4, 4/1, 5/9, 6/6, 7/8, 8/2, 9/5} {
			\fill (\x, \y) circle (4pt);
		}
	\end{tikzpicture}
	\hspace{1.5cm}
	\begin{tikzpicture}[scale=0.45]
		\draw[gray!30] (0.5,0.5) grid (9.5,9.5);
		\draw[gray!60, dashed] (0.5,0.5) -- (9.5,9.5);
		\foreach \i in {1,...,9} {
			\node[below] at (\i, 0.5) {\small $\i$};
			\node[left]  at (0.5, \i) {\small $\i$};
		}
		\foreach \x/\y in {1/3, 2/7, 3/1, 4/9, 5/6, 6/5, 7/2, 8/8, 9/4} {
			\fill (\x, \y) circle (4pt);
		}
	\end{tikzpicture}
	\caption{On the left, the plot of the permutation $374196825$. On the right, the plot of the involution $371965284$.}
	\label{figure:perm-plots}
\end{figure}
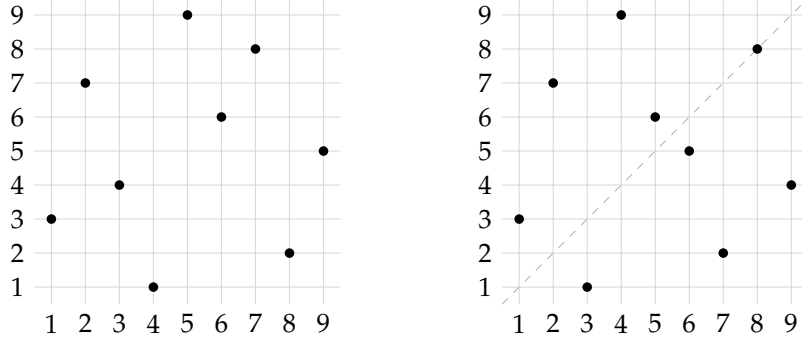

We identify a permutation $\pi$ of length $n$ with its \emph{plot}, the set of points $\{(i, \pi(i)) : 1 \leq i \leq n\}$ in the Cartesian plane. With this identification in mind, we freely use directional language when discussing entries of a permutation: we say that $\pi(j)$ is \emph{to the right} of $\pi(i)$ when $i < j$, that $\pi(j)$ is \emph{above} $\pi(i)$ when $\pi(i) < \pi(j)$, and so on. Note that the plot of an involution $\pi$ is symmetric about the diagonal line $y = x$, since $(i, j)$ belongs to the plot if and only if $(j, i)$ does. See Figure~\ref{figure:perm-plots} for the plot of the permutation $374196825$ and the involution $371965284$.

A \emph{permutation class} $\CC$ is a set of permutations that is closed downward under the containment relation, i.e., if $\pi \in \CC$ and $\sigma \leq \pi$, then we must have $\sigma \in \CC$. The \emph{basis} of a permutation class $\CC$ is the unique set of minimal permutations not in $\CC$. We write $\Av(B)$ for the class of permutations avoiding every element of $B$, and from here on we drop the set braces, writing $\Av(3124, 4312)$ instead of $\Av(\{3124, 4312\})$. For a permutation class $\CC$, we write $\CC_n$ for the set of elements of $\CC$ having length $n$. Two permutation classes with the same counting sequence $(|\CC_n|)_{n \geq 0}$ are said to be \emph{Wilf-equivalent}. For deeper background on permutation classes, we refer the interested reader to B\'ona's textbook~\cite{bona:perms-book} and Vatter's survey~\cite{vatter:perm-survey}.

In this paper we address the enumeration of pattern-avoiding involutions. We define\[
	\I(B) = \{\text{involutions } \pi : \text{$\pi$ avoids every $\beta \in B$}\}
\]
to be the \emph{involution class} containing all involutions avoiding the patterns in $B$, which may themselves not be involutions. Although an involution class is no longer a downward-closed set of permutations, it remains trivially true that if $\pi \in \I(B)$ and $\sigma$ is an involution with $\sigma \leq \pi$, then $\sigma \in \I(B)$. We let $\I_n(B)$ denote the subset of involutions in $\I(B)$ that have length $n$, and we again say that two involution classes are \emph{Wilf-equivalent} if their counting sequences are equal. Other notions of pattern avoidance for involutions may also be considered; see, for example, Fang, Hamaker, and Troyka~\cite{fang:matchings-involutions}.

Unlike the basis of a permutation class, the set $B$ is not uniquely determined by the involution class $\I(B)$. For example, $\I(231) = \I(312) = \I(231,312)$ because an involution avoids the pattern $231$ if and only if it avoids the pattern $231^{-1}$. Throughout this work, we take the convention of defining an involution class with only the lexicographically smaller permutation of each $\{\beta, \beta^{-1}\}$ pair.

For unrestricted permutations, all six patterns of length $3$ are Wilf-equivalent, each enumerated by the Catalan numbers. For patterns of length $4$, Stankova~\cite{stankova:forbidden-subseqs} and Backelin, West, and Xin~\cite{backelin:wilf-eq} showed that the twenty-four patterns of length $4$ condense into exactly three Wilf-equivalence classes, represented by $\Av(1342)$, $\Av(1234)$, and $\Av(1324)$. B\'ona~\cite{bona:1342} enumerated the first, which has an algebraic generating function, and Gessel~\cite{gessel:symmetric-functions} found a D-finite generating function for the second. The exact enumeration of $\Av(1324)$ remains open.

The study of pattern-avoiding involutions was initiated by Simion and Schmidt~\cite{simion:restricted-permutations}, who resolved all cases for patterns of length $3$, proving that $|\I_n(\beta)| = 2^{n-1}$ for $n \geq 1$ and $\beta \in \{231, 312\}$, and that $|\I_n(\beta)| = \binom{n}{\lfloor n/2 \rfloor}$ for $\beta \in \{123, 132, 213, 321\}$. Already here the Wilf classification is finer than in the unrestricted case, with two equivalence classes rather than one.

For involution classes defined by avoiding a single pattern of length $4$, the classification becomes finer still, because symmetries that identify classes in the unrestricted setting may no longer apply. For example, $\Av(1324)$ and $\Av(4231)$ are Wilf-equivalent for general permutations, yet $\I(1324)$ and $\I(4231)$ have different counting sequences. The work of Guibert~\cite{guibert:combinatoire-involutions}, Guibert, Pergola, and Pinzani~\cite{guibert:vexillary-involutions}, and Jaggard~\cite{jaggard:prefix-exchanging} together established that there are exactly eight Wilf-equivalence classes in this setting.

Progress on enumerating these eight classes has been gradual. In 1990, Gessel~\cite{gessel:symmetric-functions} showed that $\I(1234)$ is counted by the Motzkin numbers. Later, in 2008, Brignall, Huczynska, and Vatter~\cite{brignall:simple-perms} enumerated the separable involutions, $\I(2413)$. In 2016, B\'ona, Homberger, Pantone, and Vatter~\cite{bona:pattern-avoiding-involutions} then used the substitution decomposition to obtain exact generating functions for $\I(1342)$ and $\I(2341)$, bringing the total to four solved classes out of eight. 

Before the present work, four involution classes remained without exact enumerations: $\I(2431)$, $\I(3421)$, $\I(1324)$, and $\I(4231)$. The concluding remarks of~\cite{bona:pattern-avoiding-involutions} suggested that substitution decomposition might be applicable to $\I(2431)$ and $\I(3421)$, but deemed the method ``much less promising'' for $\I(1324)$ and $\I(4231)$, in which the number of simple permutations grows dramatically.

In this paper, we take a different approach. Combinatorial Exploration, developed by Albert, Bean, Claesson, Nadeau, Pantone, and Ulfarsson~\cite{combinatorial-exploration}, is an algorithmic framework for automatically discovering structural decompositions of combinatorial objects and deriving their generating functions. In the domain of permutation patterns, the framework represents sets of gridded permutations using finite data structures called \emph{tilings} and systematically applies decomposition strategies to build \emph{combinatorial specifications} from which generating functions can be extracted. It has produced exact enumerations for thousands of permutation classes, including all classes defined by avoiding any subset of length $4$ patterns except for $\Av(1324)$ and $\Av(4231)$. The website \url{https://permpal.com}~\cite{permpal-bibtex} catalogs these results.

Here, we adapt Combinatorial Exploration to the involution setting by introducing \emph{involution tilings}: tilings that are symmetric about the main diagonal and that represent only gridded involutions. We develop involution-preserving versions of the key decomposition strategies: factoring, row placement, row and column separation, and fusion. Then, we apply the Combinatorial Exploration search framework to find combinatorial specifications for $\I(3421)$ (in Section~\ref{subsec:I-3421}) and $\I(2431)$ (in Section~\ref{subsec:I-2431}), from which we compute their algebraic generating functions, proving the following.

\begin{restatable}{theorem}{smallclasstheorem}
\label{thm:3421-gf}
The generating function for $\I(3421)$ is
\[
\frac{2-5x+x^{2}-2x^{3}+x(1-x)\sqrt{1-4x^{2}}}{2(1-x)(1-2x-x^{2}-2x^{3})}.
\]
\end{restatable}

\begin{restatable}{theorem}{bigclasstheorem}
\label{thm:2431-gf}
The generating function for $\I(2431)$ is
\[
\frac{1-2x-5x^2+8x^3+3x^4+2x^5+x^4\sqrt{1-4x^2}}{(1-x)(1-2x-6x^2+8x^3+8x^4+4x^5)}.
\]
\end{restatable}

For the remaining two classes, $\I(1324)$ and $\I(4231)$, exact enumeration remains unsurprisingly out of reach. Instead, in Section~\ref{section:mosaic} we adapt the \emph{Mosaic method}, a fast counting algorithm for classical permutation avoidance~\cite{bean:mosaic-method}, to involutions by once again using involution tilings. We then use significant computational resources to compute the number of involutions in $\I(1324)$ of length up to 34 and the number of involutions in $\I(4231)$ of length up to 40. In Section~\ref{subsection:asympt-est} we use asymptotic analysis techniques to estimate their growth rates.\footnote{It is known from prior work that $\I(4231)$ has a proper growth rate, while for $\I(1324)$ this remains open; see Section~\ref{subsection:old-results}.} This analysis leads to the following two conjectures.

\begin{conjecture}
\label{conjecture:1324}
	Let $a_n = |\I_n(1324)|$. There are constants $B$, $\mu$, $g$, and $\mu_1$ with $0 < \mu_1 < 1$ such that
	\[
		a_n \sim B\,\mu^n\cdot \mu_1^{\sqrt{n}}\cdot n^g.
	\]
	Moreover, we conjecture that the quantity $\mu$, which is called the \emph{exponential growth rate} of $\I(1324)$, is exactly the square root of the corresponding quantity for $\Av(1324)$, which is estimated by Conway, Guttmann, and Zinn-Justin~\cite{conway:4231-50-terms} to be $11.600 \pm 0.003$. We numerically estimate $\mu \approx 3.40$, $\mu_1 \approx 0.200$, and $g \approx -0.1$.
\end{conjecture}

\begin{conjecture}
\label{conjecture:4231}
	Let $a_n = |\I_n(4231)|$. There are constants $B$, $\mu$, $g$, and $\mu_1$ with $0 < \mu_1 < 1$ such that
	\[
		a_n \sim B\,\mu^n\cdot \mu_1^{\sqrt{n}}\cdot n^g.
	\]
	We numerically estimate $\mu \approx 3.477$.
\end{conjecture}

The presence of the \emph{stretched exponential} term $\mu_1^{\sqrt{n}}$ in both conjectures distinguishes the predicted asymptotic behavior of these two involution classes from the known asymptotic behavior of the other six, which all take the simple power-law form $A\,\mu^n n^g$. 

We summarize our results by replicating Table 3 from~\cite{bona:pattern-avoiding-involutions}, here as Table~\ref{table:results}, updating the final rows with our results. The columns are ordered by the number of involutions of length $20$ in each class. Their work was motivated by the realization that the columns in a similar table in Jaggard~\cite{jaggard:prefix-exchanging}, which took each enumeration up to length $11$, were out of order in the sense that the ordering would likely be much different if the enumerations were taken further. They further remark that the ordering of the columns in their own table is ``likely still incorrect.'' We discuss this further in Section~\ref{section:concluding-remarks}.

\begin{table}
  \centering
  \scriptsize
  \begin{tblr}{
    colspec   = {l *{8}{r}},
    cells     = {mode=math},
    column{1} = {halign=l, rightsep=10pt},
    row{1-13}    = {halign=c},
    rowsep    = 2pt,
    colsep    = 6pt,
  }
    \hline
    {} & \bm{2431} & \bm{2341} & \bm{1342} & \bm{1234} & \bm{1324} & \bm{3421} & \bm{4231} & \bm{2413} \\[1pt]\hline

    |\Av^I_{12}(\beta)| & 
    	16238 &
    	18477 &
    	18322 &
    	15511 &
    	15272 &
    	22878 &
    	16716 &
    	27246 \\[1pt]\hline
    |\Av^I_{13}(\beta)| & 
    	40914 &
		46825 &
		47560 &
		41835 &
		40758 &
		60794 &
		46246 &
		77132
     \\[1pt]\hline
    |\Av^I_{14}(\beta)| & 
    	103954 &
		118917 &
		124358 &
		113634 &
		112280 &
		161668 &
		128414 &
		221336
     \\[1pt]\hline
    |\Av^I_{15}(\beta)| & 
    	262298 &
		301734 &
		323708 &
		310572 &
		304471 &
		429752 &
		361493 &
		635078
     \\[1pt]\hline
    |\Av^I_{16}(\beta)| & 
    	665478 &
		766525 &
		846766 &
		853467 &
		852164 &
		1142758 &
		1020506 &
		1839000 \\[1pt]\hline
    |\Av^I_{17}(\beta)| & 
   		1680726 &
		1946293 &
		2208032 &
		2356779 &
		2341980 &
		3038173 &
		2913060 &
		5331274 \\[1pt]\hline
    |\Av^I_{18}(\beta)| &
    	4260262 &
		4944614 &
		5777330 &
		6536382 &
		6640755 &
		8078606 &
		8335405 &
		15555586 \\[1pt]\hline
    |\Av^I_{19}(\beta)| & 
	  	10766470 &
		12557685 &
		15082372 &
		18199284 &
		18460066 &
		21479469 &
		24067930 &
		45465412
     \\[1pt]\hline
    |\Av^I_{20}(\beta)|  &
		27274444 &
		31900554 &
		39469786 &
		50852019 &
		52915999 &
		57113888 &
		69646035 &
		133517130 \\[1pt]\hline

    \text{growth rate}
      & \approx 2.53041
      & \approx 2.53999
      & \approx 2.61803
      & 3
      & \text{conj: }\approx 3.40559
      & \approx 2.65897
      & \text{conj: } \approx 3.477
      & \approx 3.14626 \\
    \text{reference}
      & \text{Section~\ref{subsec:I-2431}}
      & \text{\cite{bona:pattern-avoiding-involutions}}
      & \text{\cite{bona:pattern-avoiding-involutions}}
      & \text{\cite{regev:asymp-young}}
      & \text{Section~\ref{subsection:asympt-est}}
      & \text{Section~\ref{subsec:I-3421}}
      & \text{Section~\ref{subsection:asympt-est}}
      & \text{\cite{brignall:simple-perms}} \\[1pt]
    \text{OEIS}
      & {A230551}
      & {A230552}
      & {A230553}
      & {A001006}
      & {A230554}
      & {A230555}
      & {A230556}
      & {A121704} \\[1pt]\hline
  \end{tblr}
  \caption{A replication of Table 3 from~\cite{bona:pattern-avoiding-involutions} with updated results. For $\I(1324)$, it has not been proved that a proper growth rate exists, but the upper growth rate is known to lie in the interval $(3.20, 4.84)$ combining the results of~\cite{bona:pattern-avoiding-involutions} and~\cite{bevan:1324-staircase}. The conjectured value $3.40559$ is the square root of $9 + 3\sqrt{3}/2$, an algebraic candidate for $\gr(\Av(1324))$. See Section~\ref{subsection:old-results} for further discussion.}
  \label{table:results}
\end{table}

\section{Enumerating 3421-avoiding involutions and 2431-avoiding involutions}
\label{section:two-exact-enumerations}

\subsection{Combinatorial Exploration}
\label{subsection:CE-overview}

Combinatorial Exploration, introduced by Albert, Bean, Claesson, Nadeau, Pantone, and Ulfarsson~\cite{combinatorial-exploration}, is a domain-agnostic algorithmic framework for discovering structural decompositions of combinatorial objects and deriving their counting sequences and generating functions. The framework generalizes the classical symbolic method of Flajolet and Sedgewick~\cite{flajolet:ac}. Where the symbolic method works with a fixed collection of constructive operations (disjoint union, Cartesian product, sequence, and so on), Combinatorial Exploration accommodates any decomposition rule satisfying certain formal properties, and searches for applicable rules automatically, systematically decomposing a combinatorial set into simpler pieces until every piece is fully understood.

The central concepts in the framework are \emph{strategies}, \emph{rules}, and \emph{combinatorial specifications}. A \emph{strategy} is a function that takes a combinatorial set and, when it applies, produces a decomposition of that set into one or more simpler combinatorial sets, together with a formula for counting the objects in the original set from the counts of the simpler sets. For instance, a disjoint union strategy decomposes a set $\CC$ into two sets $\CC_1$ and $\CC_2$ such that $\CC = \CC_1 \sqcup \CC_2$, while a Cartesian product strategy writes $\CC$ as the set of pairs from two sets $\CC_1$ and $\CC_2$.

A \emph{rule} is a single application of a strategy to a particular set, and a \emph{combinatorial specification} is a collection of rules in which every set appearing on the right-hand side of some rule either appears on the left-hand side of exactly one rule or is a set whose enumeration is already known. When such a specification exists, it provides a complete structural description of the initial set, from which one can extract a polynomial-time counting algorithm and, typically, a system of equations whose solution is the generating function.

As a brief example, consider the permutation class $\CC = \Av(132)$. Every nonempty permutation in $\CC$ has a topmost entry, and since $\CC$ avoids $132$, every entry to the left of this topmost entry must lie above every entry to its right. One can then decompose any nonempty permutation in $\CC$ by locating its topmost entry and observing that the entries to its left form a permutation $\alpha \in \CC$ and the entries to its right also form a permutation $\beta \in \CC$. Conversely, starting with any $\alpha, \beta \in \CC$ one can assemble a unique permutation in $\CC$ by placing all entries of $\alpha$ above and to the left of all entries in $\beta$ and putting a new topmost point between them. Therefore, the nonempty permutations in $\CC$ are in bijection with $\CC \times \{1\} \times \CC$. Together with the base case (the empty permutation), and using $\CC_{\geq 1}$ to denote the nonempty permutations in $\CC$, this gives the specification
\begin{align*}
	\CC &= \{\varepsilon\} \;\sqcup\; \CC_{\geq 1}, \\
	\CC_{\geq 1} &= \CC \;\times\; \{1\} \;\times\; \CC,
\end{align*}
Writing $C(x) = \sum_{n \geq 0} |\CC_n| x^n$, this specification translates into the equation $C(x) = 1 + x\,C(x)^2$, whose solution is the generating function for the Catalan numbers.

The algorithmic search proceeds as follows. The framework maintains a queue of unexplored combinatorial sets, initially containing only the initial set, which we also call the \emph{root}. It applies strategies from a configurable \emph{strategy pack} to each set in the queue, adding any newly created sets to the queue, and periodically checks whether the accumulated rules form a combinatorial specification. When they do, the search terminates. Different combinatorial domains require different strategies, but the search infrastructure and the specification-checking algorithm are domain-agnostic. We refer the reader to Chapters 3--5 of~\cite{combinatorial-exploration} for the formal framework and to Chapter 6 for its application to permutation patterns. In the present paper, we work exclusively in the permutation patterns domain, adapting the existing strategies to the involution setting.

\subsection{Tilings}
\label{subsection:tilings}

To apply Combinatorial Exploration in the domain of permutation patterns, one needs a finite, computationally efficient way to represent the (typically infinite) sets of objects that arise during the search. Albert et al.~\cite{combinatorial-exploration} introduced \emph{tilings} for this purpose. A tiling implicitly defines a structured set of \emph{gridded permutations}, which are permutations equipped with additional geometric structure specifying where each entry sits on a grid. In this subsection, we define gridded permutations and tilings, and then describe how we adapt them to the involution setting.

\subsubsection*{Gridded permutations}

A \emph{gridded permutation} of size $n$ is a pair $(\pi, P)$, where $\pi$ is a permutation of length $n$ and $P = (c_1, \ldots, c_n)$ is a tuple of \emph{positions}, where each $c_i \in \N \times \N$ records where an entry of $\pi$ is placed on a grid: the point corresponding to the entry $\pi(i)$ sits in the unit cell whose upper-right corner has coordinates $c_i$. For this placement to be consistent, we require that the left-to-right order of the entries agrees with the left-to-right order of their cells, and the bottom-to-top order of the values agrees with the bottom-to-top order of their cells. More precisely, writing $c_i = (x_i, y_i)$, we require $x_1 \leq x_2 \leq \cdots \leq x_n$ and $y_{\pi^{-1}(1)} \leq y_{\pi^{-1}(2)} \leq \cdots \leq y_{\pi^{-1}(n)}$.

\begin{figure}
	\centering
	\begin{tikzpicture}[scale=1]
		\draw (0,0) grid (3,2);
		\foreach \i in {0,...,2} {
			\pgfmathtruncatemacro{\lab}{\i+1}
			\node[below] at (\i+0.5, -0.05) {$\lab$};
		}
		\foreach \j in {0,...,1} {
			\pgfmathtruncatemacro{\lab}{\j+1}
			\node[left] at (-0.05, \j+0.5) {$\lab$};
		}
		\fill (0.3, 0.4) circle (2.5pt);
		\fill (0.5, 1.7) circle (2.5pt);
		\fill (0.7, 1.3) circle (2.5pt);
		\fill (1.3, 0.2) circle (2.5pt);
		\fill (1.7, 0.8) circle (2.5pt);
		\fill (2.3, 0.6) circle (2.5pt);
		\fill (2.7, 1.5) circle (2.5pt);
	\end{tikzpicture}
	\caption{The gridded permutation $(2751436,\; ((1,1),\, (1,2),\, (1,2),\, (2,1),\, (2,1),\, (3,1),\, (3,2)))$ drawn on a $3 \times 2$ grid.}
	\label{fig:gridded-perm-example}
\end{figure}
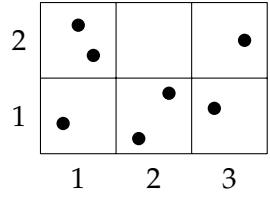

For example, Figure~\ref{fig:gridded-perm-example} shows the gridded permutation $(2751436,\; ((1,1),\,\allowbreak (1,2),\,\allowbreak (1,2),\,\allowbreak (2,1),\,\allowbreak (2,1),\,\allowbreak (3,1),\,\allowbreak (3,2)))$, which is consistent. Not every assignment of cells is consistent, e.g.,  $(21, ((1,1), (2,2)))$ is not a valid gridded permutation because the larger value is in a cell that is both below and to the left of the smaller value.

When all entries of a gridded permutation lie in the same cell $c$, we abbreviate the notation by writing $(\pi, c)$ instead of $(\pi, (c, c, \ldots, c))$, and say that the gridded permutation is \emph{localized} in cell $c$.

The notion of pattern containment extends naturally to gridded permutations: a gridded permutation $(\pi, P)$ \emph{contains} $(\sigma, Q)$ if there is a subsequence of $\pi$ that is order-isomorphic to $\sigma$ and whose entries occupy exactly the cells specified by $Q$. That is, the cells of the pattern must match exactly; there is no standardization of cells. If $(\pi, P)$ does not contain $(\sigma, Q)$, we say it \emph{avoids} $(\sigma, Q)$. For instance, the gridded permutation in Figure~\ref{fig:gridded-perm-example} contains $(312,\, ((1,2),\, (1,2),\, (3,2)))$, as witnessed by the entries $7$, $5$, $6$ at positions $2$, $3$, $7$, while it avoids $(312,\, ((1,1),\, (1,1),\, (2,1)))$.

\subsubsection*{Permutation tilings}

A \emph{tiling} is a triple $\TT = ((t,u),\, \OO,\, \RR)$, where $t$ and $u$ are positive integers specifying the dimensions of the grid, $\OO$ is a set of gridded permutations called \emph{obstructions}, and $\RR = \{\RR_1, \RR_2, \ldots, \RR_k\}$ is a set of sets of gridded permutations called \emph{requirements}. The tiling $\TT$ represents the set $\Grid(\TT)$ of all gridded permutations on a $t \times u$ grid that avoid every obstruction in $\OO$ and, for each $i$, contain at least one element of the requirement list $\RR_i$. Formally, using $\GG^{(t,u)}$ to denote the set of all gridded permutations whose positions lie in $\{1, \ldots, t\} \times \{1, \ldots, u\}$, we define
\[
	\Grid(\TT) = \big\{ g \in \GG^{(t,u)} : \text{$g$ avoids every $h \in \OO$, and for all $\RR_i \in \RR$, $g$ contains some $r \in \RR_i$} \big\}.
\]

\begin{figure}
	\centering
	\begin{tikzpicture}[baseline=(current bounding box.center)]
		\node at (0,0) {\tiling{1.0}{3}{2}{1/1}%
		{%
			{3/{(0.15, 0.25), (0.35, 0.75), (0.55, 0.50)}},%
			{3/{(2.45, 0.25), (2.65, 0.75), (2.85, 0.50)}},%
			{3/{(0.75, 0.1), (2.15, 0.6), (2.35, 0.675)}}%
		}{%
			{2/{(0.5, 0.85), (0.85, 0.5)}}%
		}};
	\end{tikzpicture}
	\caption{A tiling on a $3 \times 2$ grid, adapted from~\cite{combinatorial-exploration}. The cells $(1,2)$, $(2,1)$, and $(3,2)$ are empty (no entries allowed). The cell $(2,2)$ contains exactly one entry (indicated by the large dot). The cell $(1,1)$ has a blue requirement demanding at least two entries forming a $21$ pattern, and both cells $(1,1)$ and $(3,1)$ carry a red $132$ obstruction. There is also a \emph{crossing} obstruction, drawn in red across cells $(1,1)$ and $(3,1)$, forbidding a $123$ pattern whose first entry is in $(1,1)$ and whose last two entries are in $(3,1)$.}
	\label{fig:tiling-example}
\end{figure}
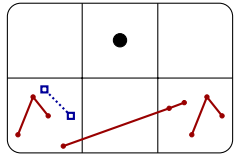

To make pictures of tilings readable, we adopt the visual conventions of~\cite{combinatorial-exploration}. Obstructions are drawn in red with solid lines and round points, while requirements are drawn in blue with dotted lines and hollow square points. The lines connecting the points of an obstruction or requirement are a visual aid indicating which points belong to the same gridded permutation. A cell containing a length $1$ obstruction is \emph{empty} (no entry of any gridded permutation in $\Grid(\TT)$ can appear there), and we draw such cells without the obstruction to reduce clutter. A cell that is required to contain exactly one entry (i.e., it has a length $1$ requirement and length $2$ obstructions $12$ and $21$) is drawn with a single large solid circle. These shortcuts are illustrated in Figure~\ref{fig:tiling-example}.

As an important special case, any permutation class $\Av(B)$ corresponds to the $1 \times 1$ tiling $\TT_{\Av(B)} = ((1,1),\, \{(\beta, (1,1)) : \beta \in B\},\, \emptyset)$ with no requirements. The map $\pi \mapsto (\pi, (1,1))$ is then a size-preserving bijection between $\Av(B)$ and $\Grid(\TT_{\Av(B)})$, and so enumerating the tiling is equivalent to enumerating the permutation class. This $1 \times 1$ tiling serves as the starting point, the root of the search tree, when Combinatorial Exploration is applied to a permutation class.

\subsubsection*{Involution tilings}

The tilings described above represent sets of gridded permutations. To adapt Combinatorial Exploration to the involution setting, we introduce \emph{involution tilings}.

First, we define the \emph{inverse} of a gridded permutation. Given $g = (\pi, (c_1, \ldots, c_n))$ with $c_i = (x_i, y_i)$, its inverse is $g^{-1} = (\pi^{-1}, (c'_1, \ldots, c'_n))$, where $c'_i = (y_j, x_j)$ for the unique $j$ satisfying $\pi(j) = i$. Informally, inverting a gridded permutation inverts the underlying permutation and reflects each cell position across the line $y = x$, swapping the cell $(x,y)$ with the cell $(y,x)$.

A \emph{gridded involution} is a gridded permutation that is equal to its own inverse. This requires both that the underlying permutation $\pi$ is an involution and that the gridding is symmetric: if the entry $\pi(i) = j$ lies in cell $(x_i, y_i)$, then the entry $\pi(j) = i$ must lie in the reflected cell $(y_i, x_i)$. If $i \neq j$, we call these two entries \emph{mirror images} of each other and together call them a \emph{mirror pair}. For this reflection to map the grid to itself, the grid must be square. The \emph{fixed points} of a gridded involution $(\pi, P)$ are the entries $i$ with $\pi(i) = i$, and the symmetry condition forces each fixed point to lie in a cell on the diagonal, i.e., a cell of the form $(x, x)$. Every non-fixed-point entry and its mirror image occupy a pair of cells that are reflections of each other across the diagonal.

An \emph{involution tiling} is a tiling $\TT = ((t,t),\, \OO,\, \RR)$ on a $t \times t$ grid that is \emph{symmetric}: the set of obstructions $\OO$ is closed under the inverse operation, and the collection of requirement lists $\RR = \{\RR_1, \ldots, \RR_k\}$ is closed under the map $\RR_i \mapsto \{g^{-1} : g \in \RR_i\}$. (That is, for each requirement list $\RR_i$, the set of inverses of its elements is also a requirement list in $\RR$, though it need not be $\RR_i$ itself.) Rather than all of $\Grid(\TT)$, an involution tiling represents only the gridded involutions it contains:
\[
	\GridI(\TT) = \big\{(\pi, P) \in \Grid(\TT) : (\pi, P)^{-1} = (\pi, P)\big\}.
\]
The distinction between gridded involutions and gridded permutations whose underlying permutation happens to be an involution is important. An involution $\pi$ may admit multiple griddings on a given tiling, and only those griddings that are themselves symmetric about the diagonal are contained in $\GridI(\TT)$. As we describe in Section~\ref{subsection:strategies}, our strategies always place points symmetrically, so all tilings constructed during our exploration are involution tilings. From here on, when the context is clear, we refer to involution tilings simply as \emph{tilings}.

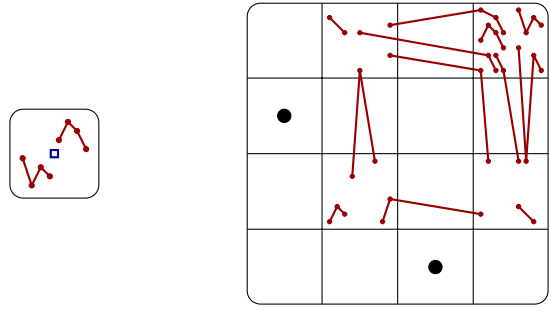
\begin{figure}
	\centering
	\begin{tikzpicture}[baseline=(current bounding box.center)]
		\node at (0,0) {\tiling{1.2}{1}{1}{}%
		{%
	        {4/{(0.550, 0.650), (0.650, 0.850), (0.750, 0.750), (0.850, 0.550)}},%
	        {4/{(0.150, 0.450), (0.250, 0.150), (0.350, 0.350), (0.450, 0.250)}}%
		}{%
			{1/{(0.5, 0.5)}}%
		}};
	\end{tikzpicture}
	\hspace{1.5cm}
	\begin{tikzpicture}[baseline=(current bounding box.center)]
		\node at (0, 0) {\tiling{1.0}{4}{4}{0/2, 2/0}%
    {%
        {2/{(1.100, 3.800), (1.300, 3.600)}},%
        {2/{(3.600, 1.300), (3.800, 1.100)}},%
        {3/{(1.100, 1.100), (1.200, 1.300), (1.300, 1.200)}},%
        {3/{(1.800, 1.100), (1.900, 1.400), (3.100, 1.200)}},%
        {3/{(1.400, 1.700), (1.500, 3.100), (1.700, 1.900)}},%
        {3/{(1.900, 3.300), (3.100, 3.100), (3.200, 1.900)}},%
        {3/{(1.500, 3.600), (3.200, 3.300), (3.300, 3.100)}},%
        {3/{(3.300, 3.300), (3.400, 3.100), (3.600, 1.900)}},%
        {4/{(1.900, 3.700), (3.100, 3.900), (3.300, 3.800), (3.400, 3.600)}},%
        {4/{(3.100, 3.500), (3.200, 3.700), (3.300, 3.600), (3.400, 3.400)}},%
        {4/{(3.600, 3.400), (3.700, 1.900), (3.800, 3.300), (3.900, 3.100)}},%
        {4/{(3.600, 3.900), (3.700, 3.600), (3.800, 3.800), (3.900, 3.700)}}%
    }%
    {%
    }};
	\end{tikzpicture}
	\caption{Two involution tilings. On the left, the tiling representing the nonempty involutions in $\I(2431)$, containing the obstructions $2431$ and $4132 = 2431^{-1}$ and a point requirement. On the right, one of the involution tilings obtained by placing the bottommost entry of the left tiling, which simultaneously places the leftmost entry by involution symmetry.}
	\label{fig:involution-tilings}
\end{figure}

Figure~\ref{fig:involution-tilings} shows two involution tilings. The left tiling represents the nonempty involutions in $\I(2431)$: a single cell containing obstructions $2431$ and its inverse $4132$, with a point requirement. Since $2431^{-1} = 4132$, both must appear as obstructions to satisfy the symmetry condition. The right tiling is a $4 \times 4$ involution tiling, one of those obtained by placing the bottommost entry of the left tiling (which, by involution symmetry, simultaneously places the leftmost entry).

The diagonal symmetry of the tiling is visible throughout. For example, cell $(2,2)$ contains a local $132$ obstruction, which is self-symmetric since $132^{-1} = 132$ (even though the visual aid of the lines connecting the points makes the picture of the obstruction look non-symmetric). The crossing obstructions respect symmetry as well. For instance, the obstruction $(2431, ((2,4), (4,4), (4,4), (4,4)))$ is paired with its inverse $(4132, ((4,4), (4,2), (4,4), (4,4)))$.

Letting $\TT$ denote the tiling on the right of Figure~\ref{fig:involution-tilings}, the gridded involution $(2143, ((1,3),\allowbreak (3,1),\allowbreak (4,4),\allowbreak (4,4)))$ is an element of $\GridI(\TT)$. On the other hand, $(312,\, ((1,3),\, (3,1),\, (4,2)))$ is not in $\GridI(\TT)$ because $312$ is not an involution. The gridded involution $(526413, ((1,3),\allowbreak (2,2), \allowbreak (2,4),\allowbreak (2,2),\allowbreak (3,1),\allowbreak (4,2)))$ is not in $\GridI(\TT)$ either: it contains the obstruction $(132, ((2,2), (2,4), (2,2)))$ and also its mirror image $(132, ((2,2), (2,2), (4,2)))$.

\subsection{Adapting permutation strategies to involutions}
\label{subsection:strategies}

The four strategies used to decompose involution tilings during the Combinatorial Exploration search are \emph{factoring}, \emph{row placement}, \emph{row and column separation}, and \emph{fusion}. Each is an adaptation of a strategy for general permutation classes from~\cite{combinatorial-exploration}, and the central requirement in each case is that the strategy preserves the diagonal symmetry of involution tilings, so that every child tiling produced is again an involution tiling. We introduce each strategy briefly here, and then concrete illustrations will follow in Sections~\ref{subsec:I-3421} and~\ref{subsec:I-2431}.

\textbf{Factor.} A tiling \emph{factors} when its nonempty cells partition into subsets that do not interact: no two subsets share a row or column, no obstruction spans cells from different subsets, and no requirement list involves cells from more than one subset. The gridded permutations on the tiling then decompose as a Cartesian product over the subsets, and the generating function is the corresponding product. For involution tilings, we additionally require that each subset be closed under cell reflection (if cell $(x, y)$ belongs to a subset, then cell $(y, x)$ must as well), ensuring that each factor is itself an involution tiling.

\textbf{Row placement.} Given a row $r$, \emph{row placement} produces a disjoint union of child tilings based on which cell of row $r$ contains the bottommost entry of that row, with one additional child for the case that row $r$ is empty. In each nonempty case, the bottommost entry is isolated into its own row, expanding the grid. (Topmost row placement is defined analogously.) For involution tilings, a bottommost row placement in row $r$ is paired with a simultaneous leftmost column placement in column $r$: each case in the disjoint union specifies which cell contains the bottommost entry of row $r$, and the symmetric operation is performed simultaneously in column $r$.

One subtlety arises at the diagonal cell $(r,r)$. When the bottommost entry of the row lies in this cell, we must distinguish whether the entry is a fixed point of the involution or not. In the non-fixed-point case, the simultaneous column placement proceeds as usual. In the fixed-point case, the entry is isolated in its own row and column, splitting the diagonal cell into four quadrants. The lower two are empty because the entry is bottommost in the row. The left two are also empty because for a fixed point, bottommost in the row implies, by the diagonal symmetry, leftmost in the column. This fixed-point sub-case appears as an additional child in the disjoint union.

\textbf{Row and column separation.} \emph{Row and column separation} applies when obstructions force all entries in certain cells of a row to lie below all entries in the remaining cells of that row. The row splits in two, with one row inheriting the cells containing lower-valued entries and the other inheriting the cells with higher-valued entries. This is a size-preserving bijection. Column separation is defined analogously. For involution tilings, separating row $r$ requires simultaneously separating column $r$ in the corresponding fashion, adding one new row and one new column to maintain a square grid and preserve the diagonal symmetry.

\textbf{Fusion.} \emph{Fusion} applies when two adjacent rows of a tiling are indistinguishable from the perspective of all other cells: every obstruction or requirement involving cells in one of the two rows has a duplicate involving the corresponding cells in the other row, and any obstruction or requirement involving entries from both rows appears together with every distribution of those entries across the two rows. Fusing merges the two rows into one. The reverse operation, unfusion, splits a row back into two rows. Since the boundary can be placed in any of the $k+1$ positions relative to the $k$ entries of the fused row, each gridded permutation with $k$ entries in the fused row corresponds to $k + 1$ gridded permutations on the unfused tiling. Column fusion is defined analogously. For involution tilings, fusing rows $i$ and $i + 1$ requires simultaneously fusing columns $i$ and $i + 1$.

With these four strategies in hand, we now turn to the Combinatorial Exploration search for $\I(3421)$ and $\I(2431)$.

\subsection{Involutions avoiding the pattern 3421}
\label{subsec:I-3421}

We present in this section a combinatorial specification for $\I(3421)$ that we found using Combinatorial Exploration. As in~\cite{combinatorial-exploration}, we present the specification as a \emph{proof tree}, a tree-shaped diagram in which each parent-to-children relationship represents that the set of gridded involutions on the parent tiling can be constructed by performing some operation on the gridded involutions on the child tilings.

\subsubsection{The proof tree}

The full specification is available as an interactive HTML rendering,\footnote{\url{https://jaypantone.github.io/pattern-avoiding-involutions/av3421_4312.html}} as an ancillary file accompanying the arXiv version of this article, and in the dataset archived on Zenodo~\cite{bean:invs-data-set}. We encourage the reader to open the HTML rendering and follow along. Clicking a tiling reveals its complete list of obstructions and requirements, including the crossing obstructions that we suppress in the figures below.

The root of the tree is the $1 \times 1$ tiling $\TT_0$ below. It represents $\GridI(\TT_0)$, the set of all gridded involutions avoiding $3421$ and $4312$ in a $1 \times 1$ grid. This set is clearly in bijection with $\I(3421)$, as all involutions avoiding $3421$ must also avoid $3421^{-1} = 4312$.

	\begin{center}
		$\TT_0$ \; \def\TilingScale{1.0}
\tiling{\TilingScale}{1}{1}{}%
{%
	{4/{(0.10, 0.5), (0.3, 0.7), (0.5, 0.3), (0.7, 0.10)}},%
	{4/{(0.3, 0.9), (0.5, 0.7), (0.7, 0.3), (0.90, 0.5)}}%
}{}%
	\end{center}
	
The first rule is a simple version of row placement, which is the strategy that decomposes a tiling into the disjoint union of several tilings based on whether a certain row is empty, or, if it is nonempty, which cell in the row contains the bottommost point. In this application, the entire tiling is a single cell, so the decomposition says that either the cell is empty or contains a bottommost point. The only case when the cell is empty is the \emph{empty gridded involution} of size $0$, and this is represented by the tiling $\TT_1$. If the tiling is nonempty, the bottommost point is either a fixed point, or not. The case where the bottommost point is a fixed point is represented by the tiling $\TT_2$. If the bottommost point is not a fixed point, then to ensure that all tilings involved are symmetric (recall that we always mean ``involution tiling'' when we say ``tiling''), we also simultaneously place the leftmost point in the cell. This produces the tiling\footnote{We preserve the labeling of tilings given by the output of the Combinatorial Exploration software, which numbers tilings according to a pre-order traversal of the tree.} $\TT_4$, which represents gridded versions of all involutions in this involution class whose bottommost point is not a fixed point. The rule below captures the obvious fact that $\I(3421)$ is the disjoint union of three sets: (1) the set containing the empty involution, (2) the nonempty $3421$-avoiding involutions in which the first entry is a fixed point, and (3) the nonempty $3421$-avoiding involutions in which the first entry is not a fixed point. However, it is not literally the case that $\GridI(\TT_0)$ is equal to the disjoint union $\GridI(\TT_1) \sqcup \GridI(\TT_2) \sqcup \GridI(\TT_4)$ because the gridded involutions in $\GridI(\TT_0)$ are on a $1 \times 1$ grid, while the gridded involutions in the children are not.
	
	\begin{center}
		\begin{tikzpicture}
			
			\node at (-0.25,0) {$\TT_0\;$ \def\TilingScale{1.0}
\tiling{\TilingScale}{1}{1}{}%
{%
	{4/{(0.10, 0.5), (0.3, 0.7), (0.5, 0.3), (0.7, 0.10)}},%
	{4/{(0.3, 0.9), (0.5, 0.7), (0.7, 0.3), (0.90, 0.5)}}%
}{}%
};
			
			\begin{scope}[shift={(-1in, -0.8in)}]
				\node at (0,0) {$\TT_1\;$\def\TilingScale{1.0}
\tiling{\TilingScale}{1}{1}{}{}{}%
};
			\end{scope}
			
			\begin{scope}[shift={(0in, -1in)}]
				\node at (-1.3,0.6) {$\TT_2$};
				\node at (0,0) {\def\TilingScale{1.0}
\tiling{\TilingScale}{2}{2}{0/0}%
{%
	{4/{(1.10, 1.5), (1.3, 1.7), (1.5, 1.3), (1.7, 1.10)}},%
	{4/{(1.3, 1.9), (1.5, 1.7), (1.7, 1.3), (1.90, 1.5)}}%
}{}%
};
			\end{scope}

			\begin{scope}[shift={(2in, -1.385in)}]
				\node at (-2.3,1.6) {$\TT_4$};
				\node at (0,0) {\def\TilingScale{1.0}
\tiling{\TilingScale}{4}{4}{0/2, 2/0}%
{%
	{2/{(1.9, 1.3), (3.1, 1.1)}},%
	{2/{(1.1, 3.1), (1.3, 1.9)}},%
	{3/{(1.15, 1.4), (1.35, 1.6), (1.55, 1.2)}},%
	{3/{(1.45, 1.8), (1.65, 1.4), (1.85, 1.6)}},%
	{3/{(1.7, 3.2), (1.8, 3.3), (1.9, 3.1)}},%
	{3/{(3.1, 1.9), (3.2, 1.7), (3.3, 1.8)}},%
	{3/{(1.4, 3.3), (1.55, 3.15), (3.15, 1.55)}},%
	{3/{(1.8, 3.9), (1.9, 3.8), (3.1, 3.6)}},%
	{3/{(1.45, 3.1), (3.05, 1.5), (3.2, 1.35)}},%
	{3/{(3.6, 3.1), (3.8, 1.9), (3.9, 1.8)}},%
	{4/{(3.3, 3.6), (3.45, 3.75), (3.6, 3.45), (3.75, 3.30)}},%
	{4/{(3.45, 3.9), (3.6, 3.75), (3.75, 3.45), (3.90, 3.6)}},%
	{4/{(3.40, 1.45), (3.55, 1.6), (3.7, 1.3), (3.85, 1.15)}},%
	{4/{(1.15, 3.85), (1.3, 3.7), (1.45, 3.4), (1.60, 3.55)}}%
}{}%
};
			\end{scope}
			
			\draw (0,-0.5) to[treeedge] (-2.5,-1.54);
			\draw (0,-0.5) to[treeedge] (0,-1.54);
			\draw (0,-0.5) to[treeedge] (5.08,-1.54);
		\end{tikzpicture}
	\end{center}
	
	As $\TT_4$ is the most complicated tiling we have seen, and as there will be many more to come, we should pause to point out a few important features of it. The first is that we have actually not drawn all of the crossing obstructions.\footnote{A \emph{crossing obstruction} is an obstruction that involves more than one cell. Recall that an obstruction that involves only one cell is called \emph{local}.} There are $10$ crossing obstructions of length $4$, each involving some combination of the cells $(2,4)$, $(4,2)$, and $(4,4)$. Five of them have the underlying permutation $3421$ while the other five have the underlying permutation $4312$. We omit them from the picture because their presence would obscure the other shorter, more interesting obstructions, but they can be seen in the interactive HTML rendering by clicking on the tiling $\TT_4$.
	
	The local obstructions $(231, (2,2))$ and $(312, (2,2))$ capture the fact that all entries in cell $(2,2)$ of any gridded involution on this tiling must avoid the patterns $231$ and $312$. Similarly, the entries in cell $(4,2)$ must avoid $312$ and $3421$ and the entries in cell $(2,4)$ must avoid $231$ and $4312$. The crossing obstruction $(21, ((2,2),(4,2)))$ implies that every entry in cell $(2,2)$ must lie below every entry in cell $(4,2)$, while its mirror image $(21, ((2,4),(2,2)))$ implies a symmetric fact about the entries in column $2$. Lastly, for now, the obstruction $(321, ((2,4), (2,4), (4,2)))$ implies that if there are two entries in cell $(2,4)$ that form a $21$ pattern, then cell $(4,2)$ must be empty.
	
	The root tiling has only two obstructions, both of length $4$. The shorter obstructions on $\TT_4$ are all the result of \emph{obstruction simplification}. Because cells $(1,3)$ and $(3,1)$ are point cells, any obstruction involving an entry in one of these cells reduces to a shorter pattern on the remaining entries. For instance, the crossing obstruction $(21, ((2,4), (2,2)))$ arises from a $3421$ pattern spanning cells $(1,3)$, $(2,4)$, $(2,2)$, and $(3,1)$: since the entries in the two point cells are always present, forbidding this $3421$ is equivalent to forbidding just the $21$ portion.

	Carrying on, there is no need to expand the tiling $\TT_1$, as the set of gridded involutions it represents (just the empty one of size $0$) is already completely understood.
		
	Next, we expand $\TT_2$. The two active cells\footnote{A cell is \emph{active} if some gridded permutation on the tiling has an entry in that cell.} of this tiling are the point cell $(1,1)$ and the cell $(2,2)$, which carries the same $3421$ and $4312$ obstructions as the root. These cells do not share a row or column, and there are no crossing obstructions on $\TT_2$ at all, so the tiling factors. Both cells lie on the diagonal, so each factor is closed under cell reflection and is therefore itself an involution tiling. The factor containing cell $(1,1)$ is the tiling $\TT_3$, a single point cell, and the factor containing cell $(2,2)$ is simply a copy of $\TT_0$, the root tiling.

	\begin{center}
		\begin{tikzpicture}
			\node at (-1.3, 0.6) {$\TT_2$};
			\node at (0,0) {\def\TilingScale{1.0}
\tiling{\TilingScale}{2}{2}{0/0}%
{%
	{4/{(1.10, 1.5), (1.3, 1.7), (1.5, 1.3), (1.7, 1.10)}},%
	{4/{(1.3, 1.9), (1.5, 1.7), (1.7, 1.3), (1.90, 1.5)}}%
}{}%
};

			\node at (-1.5,-2.4) {$\TT_3\;$\def\TilingScale{1.0}
\tiling{\TilingScale}{1}{1}{0/0}{}{}%
};

			\node at (1.,-2.4) {$\TT_0\;$\def\TilingScale{1.0}
\tiling{\TilingScale}{1}{1}{}%
{%
	{4/{(0.10, 0.5), (0.3, 0.7), (0.5, 0.3), (0.7, 0.10)}},%
	{4/{(0.3, 0.9), (0.5, 0.7), (0.7, 0.3), (0.90, 0.5)}}%
}{}%
};

			\draw (0,-1) to[treeedge] (-1.25,-1.9);
			\draw (0,-1) to[treeedge] (1.3,-1.9);
		\end{tikzpicture}
	\end{center}

	Unlike the previous rule, which represented a disjoint union, this rule expresses a Cartesian product: every gridded involution on $\TT_2$ is formed by choosing one gridded involution from $\TT_3$ (of which there is only one option, in this case) and one from $\TT_0$ and combining them in the positions prescribed by the parent tiling. This rule therefore expresses the fact that every involution in $\I(3421)$ whose bottommost entry is a fixed point consists of that fixed point together with an arbitrary involution in the class on the remaining entries above and to the right. As with $\TT_1$, the tiling $\TT_3$ needs no further expansion. 

	We now turn to expanding $\TT_4$. The crossing $21$ obstructions on $\TT_4$ force orderings between the pairs of active cells in row $2$ and column $2$, enabling a simultaneous row and column separation. The result is the $5 \times 5$ tiling $\TT_5$:

	\begin{center}
		\begin{tikzpicture}
			\node at (-1, 1.3) {$\TT_4$};
			\node at (1,0) {\def\TilingScale{0.8}
\tiling{\TilingScale}{4}{4}{0/2, 2/0}%
{%
	{2/{(1.9, 1.3), (3.1, 1.1)}},%
	{2/{(1.1, 3.1), (1.3, 1.9)}},%
	{3/{(1.15, 1.4), (1.35, 1.6), (1.55, 1.2)}},%
	{3/{(1.45, 1.8), (1.65, 1.4), (1.85, 1.6)}},%
	{3/{(1.7, 3.2), (1.8, 3.3), (1.9, 3.1)}},%
	{3/{(3.1, 1.9), (3.2, 1.7), (3.3, 1.8)}},%
	{3/{(1.4, 3.3), (1.55, 3.15), (3.15, 1.55)}},%
	{3/{(1.8, 3.9), (1.9, 3.8), (3.1, 3.6)}},%
	{3/{(1.45, 3.1), (3.05, 1.5), (3.2, 1.35)}},%
	{3/{(3.6, 3.1), (3.8, 1.9), (3.9, 1.8)}},%
	{4/{(3.3, 3.6), (3.45, 3.75), (3.6, 3.45), (3.75, 3.30)}},%
	{4/{(3.45, 3.9), (3.6, 3.75), (3.75, 3.45), (3.90, 3.6)}},%
	{4/{(3.40, 1.45), (3.55, 1.6), (3.7, 1.3), (3.85, 1.15)}},%
	{4/{(1.15, 3.85), (1.3, 3.7), (1.45, 3.4), (1.60, 3.55)}}%
}{}%
};

			\node at (3.3,0) {$\cong$};

			\node at (3.6, 1.7) {$\TT_5$};
			\node at (6,0) {\def\TilingScale{0.8}
\tiling{\TilingScale}{5}{5}{0/3, 3/0}%
{%
	{3/{(1.15, 1.4), (1.35, 1.6), (1.55, 1.2)}},%
	{3/{(1.45, 1.8), (1.65, 1.4), (1.85, 1.6)}},%
	{3/{(2.7, 4.2), (2.8, 4.3), (2.9, 4.1)}},%
	{4/{(2.15, 4.85), (2.3, 4.7), (2.45, 4.4), (2.60, 4.55)}},%
	{3/{(4.1, 2.9), (4.2, 2.7), (4.3, 2.8)}},%
	{4/{(4.40, 2.45), (4.55, 2.6), (4.7, 2.3), (4.85, 2.15)}},%
	{4/{(4.3, 4.6), (4.45, 4.75), (4.6, 4.45), (4.75, 4.30)}},%
	{4/{(4.45, 4.9), (4.6, 4.75), (4.75, 4.45), (4.90, 4.6)}},%
	{3/{(2.4, 4.3), (2.55, 4.15), (4.15, 2.55)}},%
	{3/{(2.8, 4.9), (2.9, 4.8), (4.1, 4.6)}},%
	{3/{(2.45, 4.1), (4.05, 2.5), (4.2, 2.35)}},%
	{3/{(4.6, 4.1), (4.8, 2.9), (4.9, 2.8)}}%
}{}%
};
		\end{tikzpicture}
	\end{center}

	The sets $\GridI(\TT_4)$ and $\GridI(\TT_5)$ are not literally the same, since the gridded involutions they involve lie on grids of different dimensions, but each gridded involution on $\TT_4$ maps to a unique gridded involution on $\TT_5$ by reassigning all entries in cell $(2,4)$ to cell $(3,5)$, all entries in cell $(4,2)$ to cell $(5,3)$, all entries in cell $(4,4)$ to cell $(5,5)$, and the points in cells $(1,3)$ and $(3,1)$ to cells $(1,4)$ and $(4,1)$, respectively. This map is a size-preserving bijection.
	
	The six active cells of $\TT_5$ can be partitioned into three involution-symmetric parts with no crossing obstructions between them, and the tiling factors as $\TT_5 = \TT_6 \times \TT_7 \times \TT_{12}$, where $\TT_6$ is the $2 \times 2$ tiling containing the points originally in cells $(1,4)$ and $(4,1)$, $\TT_7$ is the $1 \times 1$ tiling containing the cell that was originally $(2,2)$, and $\TT_{12}$ is the $2 \times 2$ tiling containing the cells that were originally $(3,5)$, $(5,3)$, and $(5,5)$.

	\begin{center}
		\begin{tikzpicture}
			\node at (-2.3, 1.7) {$\TT_5$};
			\node at (0,0) {\def\TilingScale{0.8}
\tiling{\TilingScale}{5}{5}{0/3, 3/0}%
{%
	{3/{(1.15, 1.4), (1.35, 1.6), (1.55, 1.2)}},%
	{3/{(1.45, 1.8), (1.65, 1.4), (1.85, 1.6)}},%
	{3/{(2.7, 4.2), (2.8, 4.3), (2.9, 4.1)}},%
	{4/{(2.15, 4.85), (2.3, 4.7), (2.45, 4.4), (2.60, 4.55)}},%
	{3/{(4.1, 2.9), (4.2, 2.7), (4.3, 2.8)}},%
	{4/{(4.40, 2.45), (4.55, 2.6), (4.7, 2.3), (4.85, 2.15)}},%
	{4/{(4.3, 4.6), (4.45, 4.75), (4.6, 4.45), (4.75, 4.30)}},%
	{4/{(4.45, 4.9), (4.6, 4.75), (4.75, 4.45), (4.90, 4.6)}},%
	{3/{(2.4, 4.3), (2.55, 4.15), (4.15, 2.55)}},%
	{3/{(2.8, 4.9), (2.9, 4.8), (4.1, 4.6)}},%
	{3/{(2.45, 4.1), (4.05, 2.5), (4.2, 2.35)}},%
	{3/{(4.6, 4.1), (4.8, 2.9), (4.9, 2.8)}}%
}{}%
};

			\node at (-4.6, -3.5) {$\TT_6$};
			\node at (-3.5, -4.0) {\def\TilingScale{0.8}
\tiling{\TilingScale}{2}{2}{0/1, 1/0}{}{}%
};

			\node at (-0.25, -3.6) {$\TT_7\;$\def\TilingScale{0.8}
\tiling{\TilingScale}{1}{1}{}%
{%
	{3/{(0.15, 0.4), (0.35, 0.6), (0.55, 0.2)}},%
	{3/{(0.45, 0.8), (0.65, 0.4), (0.85, 0.6)}}%
}{}%
};
			
			\node at (2.3, -3.5) {$\TT_{12}$};
			\node at (3.5, -4.0) {\def\TilingScale{0.8}
\tiling{\TilingScale}{2}{2}{}%
{%
	{3/{(0.7, 1.2), (0.8, 1.3), (0.9, 1.1)}},%
	{4/{(0.15, 1.85), (0.3, 1.7), (0.45, 1.4), (0.60, 1.55)}},%
	{3/{(1.1, 0.9), (1.2, 0.7), (1.3, 0.8)}},%
	{4/{(1.40, 0.45), (1.55, 0.6), (1.7, 0.3), (1.85, 0.15)}},%
	{4/{(1.3, 1.6), (1.45, 1.75), (1.6, 1.45), (1.75, 1.30)}},%
	{4/{(1.45, 1.9), (1.6, 1.75), (1.75, 1.45), (1.90, 1.6)}},%
	{3/{(0.4, 1.3), (0.55, 1.15), (1.15, 0.55)}},%
	{3/{(0.8, 1.9), (0.9, 1.8), (1.1, 1.6)}},%
	{3/{(0.45, 1.1), (1.05, 0.5), (1.2, 0.35)}},%
	{3/{(1.6, 1.1), (1.8, 0.9), (1.9, 0.8)}}%
}{}%
};

			\draw (0,-2) to[treeedge] (-3.5,-3.2);
			\draw (0,-2) to[treeedge] (0,-3.2);
			\draw (0,-2) to[treeedge] (3.5,-3.2);
		\end{tikzpicture}
	\end{center}

	The tiling $\TT_6$ represents a single gridded involution of size $2$, and so it does not need to be expanded further. The set $\GridI(\TT_7)$ is in size-preserving bijection with the set $\I(231, 312)$, whose enumeration was derived by Simion and Schmidt~\cite{simion:restricted-permutations}. So, although the algorithmic search rediscovers its enumeration by expanding $\TT_7$ into tilings $\TT_8$ through $\TT_{11}$, we omit this part of the proof tree in our exposition here. 
	
	The third child, $\TT_{12}$, is expanded by row placement into its bottommost row, which contains a single active cell, $(2,1)$. This cell is not on the diagonal, so any entry placed in it cannot be a fixed point. If the bottommost row is empty, then by involution symmetry the leftmost column is also empty, and the sole remaining active cell is $(2,2)$, which is a copy of the root tiling $\TT_0$. Otherwise, placing the bottommost entry in cell $(2,1)$ and its mirror image in cell $(1,2)$ produces the $5 \times 5$ tiling $\TT_{13}$. As with $\TT_4$, we have omitted the crossing length $4$ obstructions from the figure of $\TT_{13}$, all of which are instances of the basis patterns $3421$ and $4312$ spanning multiple cells.

	\begin{center}
		\begin{tikzpicture}
			\node at (-1.05, 0.4) {$\TT_{12}$};
			\node at (0,0) {\def\TilingScale{0.7}
\tiling{\TilingScale}{2}{2}{}%
{%
	{3/{(0.7, 1.2), (0.8, 1.3), (0.9, 1.1)}},%
	{4/{(0.15, 1.85), (0.3, 1.7), (0.45, 1.4), (0.60, 1.55)}},%
	{3/{(1.1, 0.9), (1.2, 0.7), (1.3, 0.8)}},%
	{4/{(1.40, 0.45), (1.55, 0.6), (1.7, 0.3), (1.85, 0.15)}},%
	{4/{(1.3, 1.6), (1.45, 1.75), (1.6, 1.45), (1.75, 1.30)}},%
	{4/{(1.45, 1.9), (1.6, 1.75), (1.75, 1.45), (1.90, 1.6)}},%
	{3/{(0.4, 1.3), (0.55, 1.15), (1.15, 0.55)}},%
	{3/{(0.8, 1.9), (0.9, 1.8), (1.1, 1.6)}},%
	{3/{(0.45, 1.1), (1.05, 0.5), (1.2, 0.35)}},%
	{3/{(1.6, 1.1), (1.8, 0.9), (1.9, 0.8)}}%
}{}%
};

			\node at (-2.5, -2.0) {$\TT_0\;$\def\TilingScale{0.7}
\tiling{\TilingScale}{1}{1}{}%
{%
	{4/{(0.10, 0.5), (0.3, 0.7), (0.5, 0.3), (0.7, 0.10)}},%
	{4/{(0.3, 0.9), (0.5, 0.7), (0.7, 0.3), (0.90, 0.5)}}%
}{}%
};

			\node at (1.4, -1.9) {$\TT_{13}$};
			\node at (3.5, -3.4) {\def\TilingScale{0.7}
\tiling{\TilingScale}{5}{5}{0/3, 3/0}%
{%
	{2/{(2.900, 2.100), (4.350, 1.250)}},%
	{2/{(2.900, 2.400), (4.050, 2.125)}},%
	{2/{(4.100, 1.250), (4.300, 1.100)}},%
	{3/{(1.725, 4.80), (1.875, 4.90), (2.100, 4.750)}},%
	{3/{(1.850, 4.600), (2.050, 4.70), (2.300, 4.450)}},%
	{3/{(2.42, 4.650), (2.520, 4.750), (2.62, 4.500)}},%
	{3/{(2.15, 2.4), (2.35, 2.6), (2.55, 2.2)}},%
	{3/{(2.45, 2.8), (2.65, 2.4), (2.85, 2.6)}},%
	{3/{(4.050, 2.050), (4.200, 1.450), (4.350, 1.700)}},%
	{3/{(4.60, 2.150), (4.700, 1.950), (4.750, 2.050)}},%
	{3/{(4.300, 2.500), (4.400, 2.300), (4.500, 2.400)}},%
	{3/{(1.950, 4.50), (2.500, 4.050), (4.500, 1.950)}},%
	{3/{(1.950, 4.575), (2.600, 4.100), (4.150, 2.400)}},%
	{3/{(1.950, 4.300), (4.050, 2.230), (4.500, 1.700)}},%
	{3/{(1.950, 4.400), (4.050, 2.330), (4.30, 2.050)}},%
	{3/{(2.500, 4.350), (2.850, 4.100), (4.600, 1.950)}},%
	{3/{(2.750, 4.500), (2.900, 4.350), (4.200, 2.800)}},%
	{3/{(2.950, 4.650), (4.650, 2.850), (4.850, 1.950)}},%
	{3/{(2.900, 4.500), (4.250, 2.950), (4.350, 2.850)}},%
	{3/{(4.550, 4.400), (4.750, 2.850), (4.950, 1.950)}},%
	{3/{(4.650, 4.300), (4.850, 2.950), (4.90, 2.800)}},%
	{4/{(4.400, 2.50), (4.50, 2.600), (4.600, 2.400), (4.70, 2.300)}},%
	{4/{(4.20, 4.800), (4.300, 4.90), (4.40, 4.70), (4.50, 4.60)}},%
	{4/{(2.20, 4.80), (2.30, 4.700), (2.40, 4.500), (2.50, 4.600)}},%
	{4/{(4.500, 4.90), (4.60, 4.800), (4.700, 4.600), (4.80, 4.70)}},%
	{2/{(1.250, 4.350), (2.100, 2.900)}},%
	{2/{(1.100, 4.300), (1.250, 4.100)}},%
	{2/{(2.050, 4.050), (2.400, 2.900)}},%
	{3/{(1.950, 4.950), (2.850, 4.750), (4.400, 4.550)}},%
	{3/{(2.700, 4.950), (2.950, 4.800), (4.300, 4.650)}}%
}{}%
};

			\draw (0,-0.7) to[treeedge] (-2.25,-1.65);
			\draw (0,-0.7) to[treeedge] (3.5,-1.65);
		\end{tikzpicture}
	\end{center}

	 The next step is another simultaneous row and column separation on $\TT_{13}$, following the same mechanism as the earlier separation of $\TT_4$. The crossing obstruction $(21, ((3,3), (5,3)))$ forces an ordering between two active cells sharing row $3$, and its symmetric counterpart $(21, ((3,5), (3,3)))$ does the same in column $3$. Separating both simultaneously yields the $6 \times 6$ tiling $\TT_{14}$.

	\begin{center}
		\begin{tikzpicture}
			\node at (-2.4, 1.7) {$\TT_{13}$};
			\node at (0,0) {\def\TilingScale{0.8}
\tiling{\TilingScale}{5}{5}{0/3, 3/0}%
{%
	{2/{(2.900, 2.100), (4.350, 1.250)}},%
	{2/{(2.900, 2.400), (4.050, 2.125)}},%
	{2/{(4.100, 1.250), (4.300, 1.100)}},%
	{3/{(1.725, 4.80), (1.875, 4.90), (2.100, 4.750)}},%
	{3/{(1.850, 4.600), (2.050, 4.70), (2.300, 4.450)}},%
	{3/{(2.42, 4.650), (2.520, 4.750), (2.62, 4.500)}},%
	{3/{(2.15, 2.4), (2.35, 2.6), (2.55, 2.2)}},%
	{3/{(2.45, 2.8), (2.65, 2.4), (2.85, 2.6)}},%
	{3/{(4.050, 2.050), (4.200, 1.450), (4.350, 1.700)}},%
	{3/{(4.60, 2.150), (4.700, 1.950), (4.750, 2.050)}},%
	{3/{(4.300, 2.500), (4.400, 2.300), (4.500, 2.400)}},%
	{3/{(1.950, 4.50), (2.500, 4.050), (4.500, 1.950)}},%
	{3/{(1.950, 4.575), (2.600, 4.100), (4.150, 2.400)}},%
	{3/{(1.950, 4.300), (4.050, 2.230), (4.500, 1.700)}},%
	{3/{(1.950, 4.400), (4.050, 2.330), (4.30, 2.050)}},%
	{3/{(2.500, 4.350), (2.850, 4.100), (4.600, 1.950)}},%
	{3/{(2.750, 4.500), (2.900, 4.350), (4.200, 2.800)}},%
	{3/{(2.950, 4.650), (4.650, 2.850), (4.850, 1.950)}},%
	{3/{(2.900, 4.500), (4.250, 2.950), (4.350, 2.850)}},%
	{3/{(4.550, 4.400), (4.750, 2.850), (4.950, 1.950)}},%
	{3/{(4.650, 4.300), (4.850, 2.950), (4.90, 2.800)}},%
	{4/{(4.400, 2.50), (4.50, 2.600), (4.600, 2.400), (4.70, 2.300)}},%
	{4/{(4.20, 4.800), (4.300, 4.90), (4.40, 4.70), (4.50, 4.60)}},%
	{4/{(2.20, 4.80), (2.30, 4.700), (2.40, 4.500), (2.50, 4.600)}},%
	{4/{(4.500, 4.90), (4.60, 4.800), (4.700, 4.600), (4.80, 4.70)}},%
	{2/{(1.250, 4.350), (2.100, 2.900)}},%
	{2/{(1.100, 4.300), (1.250, 4.100)}},%
	{2/{(2.050, 4.050), (2.400, 2.900)}},%
	{3/{(1.950, 4.950), (2.850, 4.750), (4.400, 4.550)}},%
	{3/{(2.700, 4.950), (2.950, 4.800), (4.300, 4.650)}}%
}{}%
};

			\node at (2.6,0) {$\cong$};

			\node at (2.9, 2.1) {$\TT_{14}$};
			\node at (5.7,0) {\def\TilingScale{0.8}
\tiling{\TilingScale}{6}{6}{0/4, 4/0}%
{%
	{2/{(2.900, 2.100), (5.350, 1.250)}},%
	{2/{(5.100, 1.250), (5.300, 1.100)}},%
	{3/{(1.725, 5.80), (1.875, 5.90), (3.100, 5.750)}},%
	{3/{(1.850, 5.600), (3.050, 5.70), (3.300, 5.450)}},%
	{3/{(2.15, 2.4), (2.35, 2.6), (2.55, 2.2)}},%
	{3/{(2.45, 2.8), (2.65, 2.4), (2.85, 2.6)}},%
	{3/{(3.42, 5.650), (3.520, 5.750), (3.62, 5.500)}},%
	{3/{(5.050, 3.050), (5.200, 1.450), (5.350, 1.700)}},%
	{3/{(5.60, 3.150), (5.700, 1.950), (5.750, 3.050)}},%
	{3/{(5.300, 3.500), (5.400, 3.300), (5.500, 3.400)}},%
	{3/{(1.950, 5.50), (3.500, 5.050), (5.500, 1.950)}},%
	{3/{(1.950, 5.575), (3.600, 5.100), (5.150, 3.400)}},%
	{3/{(1.950, 5.300), (5.050, 3.230), (5.500, 1.700)}},%
	{3/{(1.950, 5.400), (5.050, 3.330), (5.30, 3.050)}},%
	{3/{(3.500, 5.350), (3.850, 5.100), (5.600, 1.950)}},%
	{3/{(3.750, 5.500), (3.900, 5.350), (5.200, 3.800)}},%
	{3/{(3.950, 5.650), (5.650, 3.850), (5.850, 1.950)}},%
	{3/{(3.900, 5.500), (5.250, 3.950), (5.350, 3.850)}},%
	{3/{(5.550, 5.400), (5.750, 3.850), (5.950, 1.950)}},%
	{3/{(5.650, 5.300), (5.850, 3.950), (5.90, 3.800)}},%
	{4/{(5.400, 3.50), (5.50, 3.600), (5.600, 3.400), (5.70, 3.300)}},%
	{4/{(5.20, 5.800), (5.300, 5.90), (5.40, 5.70), (5.50, 5.60)}},%
	{4/{(3.20, 5.80), (3.30, 5.700), (3.40, 5.500), (3.50, 5.600)}},%
	{4/{(5.500, 5.90), (5.60, 5.800), (5.700, 5.600), (5.80, 5.70)}},%
	{2/{(1.250, 5.350), (2.100, 2.900)}},%
	{2/{(1.100, 5.300), (1.250, 5.100)}},%
	{3/{(1.950, 5.950), (3.850, 5.750), (5.400, 5.550)}},%
	{3/{(3.700, 5.950), (3.950, 5.800), (5.300, 5.650)}}%
}{}%
};
		\end{tikzpicture}
	\end{center}

	The point cells $(1,5)$ and $(5,1)$ of $\TT_{14}$ can be factored out, and again form a copy of $\TT_6$. The remaining six active cells form the $4 \times 4$ tiling $\TT_{15}$. This factor rule is shown below.

	\begin{center}
		\begin{tikzpicture}
			\node at (-2.75, 2.1) {$\TT_{14}$};
			\node at (0,0) {\def\TilingScale{0.8}};

			\node at (-3.9, -3.75) {$\TT_6$};
			\node at (-2.8, -4.2) {\def\TilingScale{0.8}
\tiling{\TilingScale}{2}{2}{0/1, 1/0}{}{}%
};

			\node at (0.05, -3.75) {$\TT_{15}$};
			\node at (2.0, -5.0) {\def\TilingScale{0.8}
\tiling{\TilingScale}{4}{4}{}%
{%
	{2/{(1.900, 1.100), (3.350, 0.250)}},%
	{2/{(3.100, 0.250), (3.300, 0.100)}},%
	{3/{(0.725, 3.80), (0.875, 3.90), (2.100, 3.750)}},%
	{3/{(0.850, 3.600), (2.050, 3.70), (2.300, 3.450)}},%
	{3/{(2.42, 3.650), (2.520, 3.750), (2.62, 3.500)}},%
	{3/{(1.15, 1.4), (1.35, 1.6), (1.55, 1.2)}},%
	{3/{(1.45, 1.8), (1.65, 1.4), (1.85, 1.6)}},%
	{3/{(3.050, 2.050), (3.200, 0.450), (3.350, 0.700)}},%
	{3/{(3.60, 2.150), (3.700, 0.950), (3.750, 2.050)}},%
	{3/{(3.300, 2.500), (3.400, 2.300), (3.500, 2.400)}},%
	{3/{(0.950, 3.50), (2.500, 3.050), (3.500, 0.950)}},%
	{3/{(0.950, 3.575), (2.600, 3.100), (3.150, 2.400)}},%
	{3/{(0.950, 3.300), (3.050, 2.230), (3.500, 0.700)}},%
	{3/{(0.950, 3.400), (3.050, 2.330), (3.30, 2.050)}},%
	{3/{(2.500, 3.350), (2.850, 3.100), (3.600, 0.950)}},%
	{3/{(2.750, 3.500), (2.900, 3.350), (3.200, 2.800)}},%
	{3/{(2.950, 3.650), (3.650, 2.850), (3.850, 0.950)}},%
	{3/{(2.900, 3.500), (3.250, 2.950), (3.350, 2.850)}},%
	{3/{(3.550, 3.400), (3.750, 2.850), (3.950, 0.950)}},%
	{3/{(3.650, 3.300), (3.850, 2.950), (3.90, 2.800)}},%
	{4/{(3.400, 2.50), (3.50, 2.600), (3.600, 2.400), (3.70, 2.300)}},%
	{4/{(3.20, 3.800), (3.300, 3.90), (3.40, 3.70), (3.50, 3.60)}},%
	{4/{(2.20, 3.80), (2.30, 3.700), (2.40, 3.500), (2.50, 3.600)}},%
	{4/{(3.500, 3.90), (3.60, 3.800), (3.700, 3.600), (3.80, 3.70)}},%
	{2/{(0.250, 3.350), (1.100, 1.900)}},%
	{2/{(0.100, 3.300), (0.250, 3.100)}},%
	{3/{(0.950, 3.950), (2.850, 3.750), (3.400, 3.550)}},%
	{3/{(2.700, 3.950), (2.950, 3.800), (3.300, 3.650)}}%
}{}%
};

			\draw (0,-2.4) to[treeedge] (-2.8,-3.4);
			\draw (0,-2.4) to[treeedge] (2.0,-3.4);
		\end{tikzpicture}
	\end{center}

	Row placement on $\TT_{15}$ targets the bottommost row, which again contains a single off-diagonal cell, $(4,1)$. If the row is empty, then column $1$ is empty by symmetry, and the remaining cells form the $3 \times 3$ tiling $\TT_{16}$. If the row is nonempty, placing the bottommost entry in cell $(4,1)$ and its mirror image in cell $(1,4)$ produces the tiling $\TT_{17}$. Note that we have omitted many crossing $321$ obstructions from the depiction of $\TT_{17}$ to make it more readable. We now point out how a few of the obstructions on $\TT_{15}$ evolve to obstructions on $\TT_{17}$ once the points are placed. Cell $(2,2)$ in $\TT_{15}$ becomes empty in $\TT_{17}$ (where it is now cell $(3,3)$) because of the crossing obstructions $(21,((1,4),(2,2)))$ and $(21,((2,2),(4,1)))$ in $\TT_{15}$. Cell $(4,4)$ in $\TT_{15}$ splits into the four cells $(4,4)$, $(4,6)$, $(6,4)$, and $(6,6)$ in $\TT_{17}$. Cell $(3,4)$ in $\TT_{15}$, which previously carried local obstructions $231$ and $4312$, becomes the monotone cell $(3,6)$ in $\TT_{17}$ because the crossing obstruction $(321,((3,4),(3,4),(4,1)))$ in $\TT_{15}$ reduces to a local $21$ once the point is placed. By symmetry, cell $(6,3)$ in $\TT_{17}$ is also monotone.

	\begin{center}
		\begin{tikzpicture}
			\node at (-1.95, 1.3) {$\TT_{15}$};
			\node at (0,0) {\def\TilingScale{0.8}};

			\node at (-4.55, -2.95) {$\TT_{16}$};
			\node at (-3.0, -3.8) {\def\TilingScale{0.8}
\tiling{\TilingScale}{3}{3}{}%
{%
	{3/{(0.15, 0.4), (0.35, 0.6), (0.55, 0.2)}},%
	{3/{(0.45, 0.8), (0.65, 0.4), (0.85, 0.6)}},%
	{3/{(1.3, 2.7), (1.450, 2.850), (1.6, 2.550)}},%
	{4/{(1.15, 2.60), (1.30, 2.450), (1.45, 2.150), (1.60, 2.300)}},%
	{3/{(2.150, 1.450), (2.300, 1.150), (2.450, 1.300)}},%
	{4/{(2.400, 1.450), (2.550, 1.600), (2.700, 1.300), (2.85, 1.150)}},%
	{4/{(2.15, 2.70), (2.250, 2.85), (2.35, 2.55), (2.45, 2.40)}},%
	{4/{(2.50, 2.85), (2.6, 2.700), (2.70, 2.400), (2.8, 2.66)}},%
	{3/{(1.750, 2.400), (1.900, 2.250), (2.200, 1.6500)}},%
	{3/{(1.850, 2.500), (2.30, 1.850), (2.45, 1.70)}},%
	{3/{(2.550, 2.150), (2.70, 1.850), (2.850, 1.700)}},%
	{3/{(1.700, 2.850), (1.850, 2.700), (2.150, 2.550)}}%
}{}%
};

			\node at (0.75, -2.95) {$\TT_{17}$};
			\node at (3.5, -5.0) {\def\TilingScale{0.8}
\tiling{\TilingScale}{6}{6}{0/4, 4/0}%
{%
	{2/{(1.100, 5.300), (1.250, 5.100)}},%
	{2/{(1.900, 5.200), (2.150, 5.050)}},%
	{2/{(1.400, 5.100), (3.100, 3.200)}},%
	{2/{(2.250, 5.250), (2.400, 5.100)}},%
	{2/{(3.400, 5.100), (3.750, 3.900)}},%
	{3/{(1.700, 5.550), (1.900, 5.700), (3.300, 5.400)}},%
	{3/{(1.950, 5.300), (2.850, 5.450), (3.300, 5.150)}},%
	{3/{(1.700, 5.750), (3.150, 5.900), (3.350, 5.550)}},%
	{3/{(2.700, 5.200), (2.850, 5.350), (3.150, 5.050)}},%
	{3/{(2.850, 5.650), (3.100, 5.800), (3.250, 5.500)}},%
	{3/{(3.350, 3.450), (3.550, 3.600), (3.750, 3.300)}},%
	{3/{(3.600, 5.750), (3.750, 5.900), (3.900, 5.600)}},%
	{3/{(3.150, 3.500), (3.300, 3.200), (3.500, 3.350)}},%
	{4/{(5.100, 5.450), (5.300, 5.600), (5.450, 5.300), (5.700, 5.150)}},%
	{4/{(3.350, 5.750), (3.550, 5.600), (3.700, 5.300), (3.900, 5.450)}},%
	{4/{(5.300, 5.850), (5.500, 5.700), (5.650, 5.400), (5.850, 5.550)}},%
	{2/{(5.100, 1.250), (5.300, 1.100)}},%
	{2/{(3.200, 3.100), (5.100, 1.400)}},%
	{2/{(2.550, 5.100), (3.100, 3.850)}},%
	{3/{(5.500, 3.450), (5.750, 1.700), (5.900, 3.300)}},%
	{3/{(5.400, 3.300), (5.550, 1.700), (5.700, 1.900)}},%
	{2/{(5.050, 2.150), (5.200, 1.900)}},%
	{3/{(5.050, 3.100), (5.300, 1.950), (5.450, 2.150)}},%
	{2/{(5.050, 2.350), (5.200, 2.200)}},%
	{2/{(3.850, 3.100), (5.100, 2.550)}},%
	{3/{(5.100, 3.250), (5.250, 2.800), (5.400, 2.950)}},%
	{3/{(5.550, 3.250), (5.700, 2.850), (5.850, 3.100)}},%
	{2/{(3.900, 3.750), (5.100, 3.450)}},%
	{4/{(5.300, 3.700), (5.450, 3.900), (5.600, 3.550), (5.750, 3.350)}},%
	{3/{(5.600, 3.900), (5.750, 3.600), (5.900, 3.750)}}%
}{}%
};

			\draw (0,-1.6) to[treeedge] (-3.0,-2.6);
			\draw (0,-1.6) to[treeedge] (3.5,-2.6);
		\end{tikzpicture}
	\end{center}

	The four active cells of $\TT_{16}$ are $(1,1)$, $(2,3)$, $(3,2)$, and $(3,3)$. Cell $(1,1)$ shares no row or column with the other three, and no crossing obstruction spans between them, so the tiling factors into two tilings. The first is the $1 \times 1$ tiling containing cell $(1,1)$ of $\TT_{16}$, which we have already seen as $\TT_7$. The second is the $2 \times 2$ tiling containing the obstructions in cells $(2,3)$, $(3,2)$, and $(3,3)$ of $\TT_{16}$, which we have already seen as $\TT_{12}$. So, $\TT_{16}$ is decomposed into two tilings that we have already seen, as shown in the rule below, and no further expansion down this branch is required.
	
	\begin{center}
		\begin{tikzpicture}
			\node at (-1.55, 0.85) {$\TT_{16}$};
			\node at (0,0) {\def\TilingScale{0.8}
\tiling{\TilingScale}{3}{3}{}%
{%
	{3/{(0.15, 0.4), (0.35, 0.6), (0.55, 0.2)}},%
	{3/{(0.45, 0.8), (0.65, 0.4), (0.85, 0.6)}},%
	{3/{(1.3, 2.7), (1.450, 2.850), (1.6, 2.550)}},%
	{4/{(1.15, 2.60), (1.30, 2.450), (1.45, 2.150), (1.60, 2.300)}},%
	{3/{(2.150, 1.450), (2.300, 1.150), (2.450, 1.300)}},%
	{4/{(2.400, 1.450), (2.550, 1.600), (2.700, 1.300), (2.85, 1.150)}},%
	{4/{(2.15, 2.70), (2.250, 2.85), (2.35, 2.55), (2.45, 2.40)}},%
	{4/{(2.50, 2.85), (2.6, 2.700), (2.70, 2.400), (2.8, 2.66)}},%
	{3/{(1.750, 2.400), (1.900, 2.250), (2.200, 1.6500)}},%
	{3/{(1.850, 2.500), (2.30, 1.850), (2.45, 1.70)}},%
	{3/{(2.550, 2.150), (2.70, 1.850), (2.850, 1.700)}},%
	{3/{(1.700, 2.850), (1.850, 2.700), (2.150, 2.550)}}%
}{}%
};

			\node at (-2.5, -2.6) {$\TT_7\;$\def\TilingScale{0.8}
\tiling{\TilingScale}{1}{1}{}%
{%
	{3/{(0.15, 0.4), (0.35, 0.6), (0.55, 0.2)}},%
	{3/{(0.45, 0.8), (0.65, 0.4), (0.85, 0.6)}}%
}{}%
};

			\node at (0.85, -2.5) {$\TT_{12}$};
			\node at (2.0, -3.0) {\def\TilingScale{0.8}
\tiling{\TilingScale}{2}{2}{}%
{%
	{3/{(0.7, 1.2), (0.8, 1.3), (0.9, 1.1)}},%
	{4/{(0.15, 1.85), (0.3, 1.7), (0.45, 1.4), (0.60, 1.55)}},%
	{3/{(1.1, 0.9), (1.2, 0.7), (1.3, 0.8)}},%
	{4/{(1.40, 0.45), (1.55, 0.6), (1.7, 0.3), (1.85, 0.15)}},%
	{4/{(1.3, 1.6), (1.45, 1.75), (1.6, 1.45), (1.75, 1.30)}},%
	{4/{(1.45, 1.9), (1.6, 1.75), (1.75, 1.45), (1.90, 1.6)}},%
	{3/{(0.4, 1.3), (0.55, 1.15), (1.15, 0.55)}},%
	{3/{(0.8, 1.9), (0.9, 1.8), (1.1, 1.6)}},%
	{3/{(0.45, 1.1), (1.05, 0.5), (1.2, 0.35)}},%
	{3/{(1.6, 1.1), (1.8, 0.9), (1.9, 0.8)}}%
}{}%
};

			\draw (0,-1.2) to[treeedge] (-2.25,-2.2);
			\draw (0,-1.2) to[treeedge] (2.0,-2.2);
		\end{tikzpicture}
	\end{center}

	Returning to $\TT_{17}$, the crossing $21$ obstructions again enable a simultaneous row and column separation. The obstruction $(21, ((4,4),(6,4)))$ forces an ordering of the entries between the two active cells sharing row $4$, and its symmetric counterpart $(21, ((4,6),(4,4)))$ does the same in column $4$. Separating both simultaneously yields the $7 \times 7$ tiling $\TT_{18}$,  shown below.
	
	\begin{center}
		\begin{tikzpicture}
			\node at (-2.35, 2.05) {$\TT_{17}$};
			\node at (0.4,0) {\def\TilingScale{0.8}
\tiling{\TilingScale}{6}{6}{0/4, 4/0}%
{%
	{2/{(1.100, 5.300), (1.250, 5.100)}},%
	{2/{(1.900, 5.200), (2.150, 5.050)}},%
	{2/{(1.400, 5.100), (3.100, 3.200)}},%
	{2/{(2.250, 5.250), (2.400, 5.100)}},%
	{2/{(3.400, 5.100), (3.750, 3.900)}},%
	{3/{(1.700, 5.550), (1.900, 5.700), (3.300, 5.400)}},%
	{3/{(1.950, 5.300), (2.850, 5.450), (3.300, 5.150)}},%
	{3/{(1.700, 5.750), (3.150, 5.900), (3.350, 5.550)}},%
	{3/{(2.700, 5.200), (2.850, 5.350), (3.150, 5.050)}},%
	{3/{(2.850, 5.650), (3.100, 5.800), (3.250, 5.500)}},%
	{3/{(3.350, 3.450), (3.550, 3.600), (3.750, 3.300)}},%
	{3/{(3.600, 5.750), (3.750, 5.900), (3.900, 5.600)}},%
	{3/{(3.150, 3.500), (3.300, 3.200), (3.500, 3.350)}},%
	{4/{(5.100, 5.450), (5.300, 5.600), (5.450, 5.300), (5.700, 5.150)}},%
	{4/{(3.350, 5.750), (3.550, 5.600), (3.700, 5.300), (3.900, 5.450)}},%
	{4/{(5.300, 5.850), (5.500, 5.700), (5.650, 5.400), (5.850, 5.550)}},%
	{2/{(5.100, 1.250), (5.300, 1.100)}},%
	{2/{(3.200, 3.100), (5.100, 1.400)}},%
	{2/{(2.550, 5.100), (3.100, 3.850)}},%
	{3/{(5.500, 3.450), (5.750, 1.700), (5.900, 3.300)}},%
	{3/{(5.400, 3.300), (5.550, 1.700), (5.700, 1.900)}},%
	{2/{(5.050, 2.150), (5.200, 1.900)}},%
	{3/{(5.050, 3.100), (5.300, 1.950), (5.450, 2.150)}},%
	{2/{(5.050, 2.350), (5.200, 2.200)}},%
	{2/{(3.850, 3.100), (5.100, 2.550)}},%
	{3/{(5.100, 3.250), (5.250, 2.800), (5.400, 2.950)}},%
	{3/{(5.550, 3.250), (5.700, 2.850), (5.850, 3.100)}},%
	{2/{(3.900, 3.750), (5.100, 3.450)}},%
	{4/{(5.300, 3.700), (5.450, 3.900), (5.600, 3.550), (5.750, 3.350)}},%
	{3/{(5.600, 3.900), (5.750, 3.600), (5.900, 3.750)}}%
}{}%
};

			\node at (3.3,0) {$\cong$};

			\node at (3.45, 2.45) {$\TT_{18}$};
			\node at (6.6,0) {\def\TilingScale{0.8}
\tiling{\TilingScale}{7}{7}{0/5, 5/0}%
{%
	{2/{(1.100, 6.300), (1.250, 6.100)}},%
	{2/{(1.900, 6.200), (2.150, 6.050)}},%
	{2/{(1.400, 6.100), (3.100, 3.200)}},%
	{2/{(2.250, 6.250), (2.400, 6.100)}},%
	{3/{(1.700, 6.550), (1.900, 6.700), (4.300, 6.400)}},%
	{3/{(1.950, 6.300), (2.850, 6.450), (4.300, 6.150)}},%
	{3/{(1.700, 6.750), (4.150, 6.900), (4.350, 6.550)}},%
	{3/{(2.700, 6.200), (2.850, 6.350), (4.150, 6.050)}},%
	{3/{(2.850, 6.650), (4.100, 6.800), (4.250, 6.500)}},%
	{3/{(3.350, 3.450), (3.550, 3.600), (3.750, 3.300)}},%
	{3/{(4.600, 6.750), (4.750, 6.900), (4.900, 6.600)}},%
	{3/{(3.150, 3.500), (3.300, 3.200), (3.500, 3.350)}},%
	{4/{(6.100, 6.450), (6.300, 6.600), (6.450, 6.300), (6.700, 6.150)}},%
	{4/{(4.350, 6.750), (4.550, 6.600), (4.700, 6.300), (4.900, 6.450)}},%
	{4/{(6.300, 6.850), (6.500, 6.700), (6.650, 6.400), (6.850, 6.550)}},%
	{2/{(6.100, 1.250), (6.300, 1.100)}},%
	{2/{(3.200, 3.100), (6.100, 1.400)}},%
	{2/{(2.550, 6.100), (3.100, 3.850)}},%
	{3/{(6.500, 4.450), (6.750, 1.700), (6.900, 4.300)}},%
	{3/{(6.400, 4.300), (6.550, 1.700), (6.700, 1.900)}},%
	{2/{(6.050, 2.150), (6.200, 1.900)}},%
	{3/{(6.050, 4.100), (6.300, 1.950), (6.450, 2.150)}},%
	{2/{(6.050, 2.350), (6.200, 2.200)}},%
	{2/{(3.850, 3.100), (6.100, 2.550)}},%
	{3/{(6.100, 4.250), (6.250, 2.800), (6.400, 2.950)}},%
	{3/{(6.550, 4.250), (6.700, 2.850), (6.850, 4.100)}},%
	{4/{(6.300, 4.700), (6.450, 4.900), (6.600, 4.550), (6.750, 4.350)}},%
	{3/{(6.600, 4.900), (6.750, 4.600), (6.900, 4.750)}}%
}{}%
};
		\end{tikzpicture}
	\end{center}

	We may again factor. The point cells $(1,6)$ and $(6,1)$ of $\TT_{18}$ again form a copy of $\TT_6$, while the remaining eight active cells together form the new $5 \times 5$ tiling $\TT_{19}$.

	\begin{center}
		\begin{tikzpicture}
			\node at (-3.15, 2.45) {$\TT_{18}$};
			\node at (0,0) {\def\TilingScale{0.8}};

			\node at (-3.9, -4.5) {$\TT_6$};
			\node at (-2.8, -5.0) {\def\TilingScale{0.8}
\tiling{\TilingScale}{2}{2}{0/1, 1/0}{}{}%
};

			\node at (0, -4.8) {$\TT_{19}$};
			\node at (2.5, -6.2) {\def\TilingScale{0.8}
\tiling{\TilingScale}{5}{5}{}%
{%
	{2/{(0.100, 4.300), (0.250, 4.100)}},%
	{2/{(0.900, 4.200), (1.150, 4.050)}},%
	{2/{(0.400, 4.100), (2.100, 2.200)}},%
	{2/{(1.250, 4.250), (1.400, 4.100)}},%
	{3/{(0.700, 4.550), (0.900, 4.700), (3.300, 4.400)}},%
	{3/{(0.950, 4.300), (1.850, 4.450), (3.300, 4.150)}},%
	{3/{(0.700, 4.750), (3.150, 4.900), (3.350, 4.550)}},%
	{3/{(1.700, 4.200), (1.850, 4.350), (3.150, 4.050)}},%
	{3/{(1.850, 4.650), (3.100, 4.800), (3.250, 4.500)}},%
	{3/{(2.350, 2.450), (2.550, 2.600), (2.750, 2.300)}},%
	{3/{(3.600, 4.750), (3.750, 4.900), (3.900, 4.600)}},%
	{3/{(2.150, 2.500), (2.300, 2.200), (2.500, 2.350)}},%
	{4/{(4.100, 4.450), (4.300, 4.600), (4.450, 4.300), (4.700, 4.150)}},%
	{4/{(3.350, 4.750), (3.550, 4.600), (3.700, 4.300), (3.900, 4.450)}},%
	{4/{(4.300, 4.850), (4.500, 4.700), (4.650, 4.400), (4.850, 4.550)}},%
	{2/{(4.100, 0.250), (4.300, 0.100)}},%
	{2/{(2.200, 2.100), (4.100, 0.400)}},%
	{2/{(1.550, 4.100), (2.100, 2.850)}},%
	{3/{(4.500, 3.450), (4.750, 0.700), (4.900, 3.300)}},%
	{3/{(4.400, 3.300), (4.550, 0.700), (4.700, 0.900)}},%
	{2/{(4.050, 1.150), (4.200, 0.900)}},%
	{3/{(4.050, 3.100), (4.300, 0.950), (4.450, 1.150)}},%
	{2/{(4.050, 1.350), (4.200, 1.200)}},%
	{2/{(2.850, 2.100), (4.100, 1.550)}},%
	{3/{(4.100, 3.250), (4.250, 1.800), (4.400, 1.950)}},%
	{3/{(4.550, 3.250), (4.700, 1.850), (4.850, 3.100)}},%
	{4/{(4.300, 3.700), (4.450, 3.900), (4.600, 3.550), (4.750, 3.350)}},%
	{3/{(4.600, 3.900), (4.750, 3.600), (4.900, 3.750)}}%
}{}%
};

			\draw (0,-2.8) to[treeedge] (-2.8,-4.2);
			\draw (0,-2.8) to[treeedge] (2.52,-4.2);
		\end{tikzpicture}
	\end{center}

	Tiling $\TT_{19}$ is the first and only tiling we encounter in this subsection to which the fusion strategy can be applied. As described in Section~\ref{subsection:strategies}, two adjacent rows can be fused when they are indistinguishable: every obstruction involving cells in one row has a partner involving the corresponding cells in the other. In $\TT_{19}$, rows $1$ and $2$ satisfy this condition. For example, the obstruction $(21, ((3,3),(5,1)))$ is partnered with $(21, ((3,3),(5,2)))$. When an obstruction involves entries from both of these rows, every distribution of those entries across the two rows appears as a separate obstruction. For instance, a $21$ pattern among entries in rows $1$ and $2$ of column $5$ is forbidden regardless of which row the entries lie in, giving rise to the three obstructions $(21, ((5,1),(5,1)))$, $(21, ((5,2),(5,1)))$, and $(21, ((5,2),(5,2)))$. After fusion merges these two rows, the corresponding cell inherits a local $21$ obstruction. Similarly, the three obstructions $(312, ((5,4),(5,1),(5,1)))$, $(312, ((5,4),(5,1),(5,2)))$, and $(312, ((5,4),(5,2),(5,2)))$ all appear on $\TT_{19}$. Fusing rows $1$ and $2$ merges them into a single row, and by involution symmetry columns $1$ and $2$ are fused simultaneously. The resulting $4 \times 4$ tiling is identical to $\TT_{15}$, although we warn the reader that the crossing $321$ obstructions that were hidden on $\TT_{19}$ are now shown again on $\TT_{15}$.

	\begin{center}
		\begin{tikzpicture}
			\node at (-2.37, 1.65) {$\TT_{19}$};
			\node at (0,0) {\def\TilingScale{0.8}};

			\node at (-1.95, -3.03) {$\TT_{15}$};
			\node at (0,-4.3) {\def\TilingScale{0.8}};

			\draw (0,-2) to[treeedge] (0,-2.7);
		\end{tikzpicture}
	\end{center}

	Since $\TT_{15}$ appears earlier in the proof tree and has already been expanded, this fusion closes the last open branch of exploration.

\begin{figure}[p]
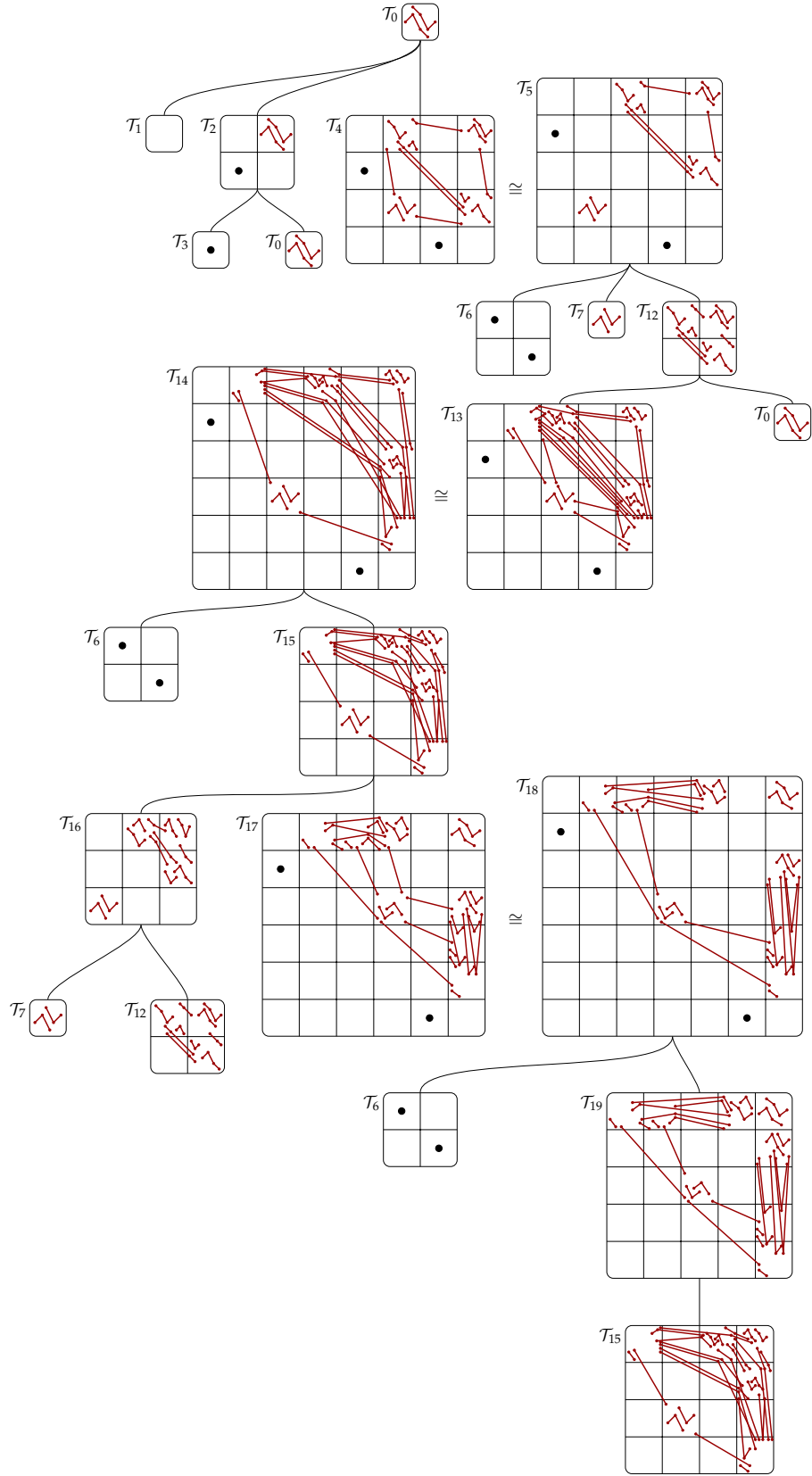

\centering
\scalebox{0.69}{
\begin{tikzpicture}[
	every node/.style={inner sep=0pt, outer sep=0pt},
	tilinglabel/.style={anchor=north east, font=\small, xshift=-2pt, yshift=-2pt, overlay},
	treeedge/.style={out=-90, in=90, looseness=0.4},
	ann/.style={font=\scriptsize, text=gray!70!black},
	scale=1
]
\node (T0) at (0, 0) {\def\TilingScale{0.8}\tilingroundedfalse\input{tikz-files/small-class-T0}\tilingroundedtrue};
\node[tilinglabel] at (T0.north west) {$\TT_0$};

\node[anchor=north] (T1) at (-5.5, -2) {\def\TilingScale{0.8}\tilingroundedfalse\input{tikz-files/small-class-T1}\tilingroundedtrue};
\node[tilinglabel] at (T1.north west) {$\TT_1$};

\node[anchor=north] (T2) at (-3.5, -2) {\def\TilingScale{0.8}\tilingroundedfalse\input{tikz-files/small-class-T2}\tilingroundedtrue};
\node[tilinglabel] at (T2.north west) {$\TT_2$};

\node[anchor=north] (T4) at (0, -2) {\def\TilingScale{0.8}\tilingroundedfalse\input{tikz-files/small-class-T4}\tilingroundedtrue};
\node[tilinglabel] at (T4.north west) {$\TT_4$};
\node at (2.05, -3.6) {$\cong$};
\node[anchor=north] (T5) at (4.5, -1.2) {\def\TilingScale{0.8}\tilingroundedfalse\input{tikz-files/small-class-T5}\tilingroundedtrue};
\node[tilinglabel] at (T5.north west) {$\TT_5$};

\draw (T0.south) to[treeedge] (T1.north);
\draw (T0.south) to[treeedge] (T2.north);
\draw (T0.south) to[treeedge] (T4.north);

\node[anchor=north] (T3) at (-4.5, -4.5) {\def\TilingScale{0.8}\tilingroundedfalse\input{tikz-files/small-class-T3}\tilingroundedtrue};
\node[tilinglabel] at (T3.north west) {$\TT_3$};
\node[anchor=north] (T0ref1) at (-2.5, -4.5) {\def\TilingScale{0.8}\tilingroundedfalse\input{tikz-files/small-class-T0}\tilingroundedtrue};
\node[tilinglabel] at (T0ref1.north west) {$\TT_0$};

\draw (T2.south) to[treeedge] (T3.north);
\draw (T2.south) to[treeedge] (T0ref1.north);

\node[anchor=north] (T6) at (2, -6) {\def\TilingScale{0.8}\tilingroundedfalse\input{tikz-files/small-class-T6}\tilingroundedtrue};
\node[tilinglabel] at (T6.north west) {$\TT_6$};
\node[anchor=north] (T7) at (4, -6) {\def\TilingScale{0.8}\tilingroundedfalse\input{tikz-files/small-class-T7}\tilingroundedtrue};
\node[tilinglabel] at (T7.north west) {$\TT_7$};
\node[anchor=north] (T12) at (6, -6) {\def\TilingScale{0.8}\tilingroundedfalse\input{tikz-files/small-class-T12}\tilingroundedtrue};
\node[tilinglabel] at (T12.north west) {$\TT_{12}$};

\draw (T5.south) to[treeedge] (T6.north);
\draw (T5.south) to[treeedge] (T7.north);
\draw (T5.south) to[treeedge] (T12.north);

\node[anchor=north] (T0ref2) at (8, -8.2) {\def\TilingScale{0.8}\tilingroundedfalse\input{tikz-files/small-class-T0}\tilingroundedtrue};
\node[tilinglabel] at (T0ref2.north west) {$\TT_0$};
\node[anchor=north] (T13) at (3, -8.2) {\def\TilingScale{0.8}\tilingroundedfalse\input{tikz-files/small-class-T13}\tilingroundedtrue};
\node[tilinglabel] at (T13.north west) {$\TT_{13}$};
\node at (0.45, -10.2) {$\cong$};
\node[anchor=north] (T14) at (-2.5, -7.4) {\def\TilingScale{0.8}\tilingroundedfalse\input{tikz-files/small-class-T14}\tilingroundedtrue};
\node[tilinglabel] at (T14.north west) {$\TT_{14}$};

\draw (T12.south) to[treeedge] (T0ref2.north);
\draw (T12.south) to[treeedge] (T13.north);

\node[anchor=north] (T6ref1) at (-6, -13) {\def\TilingScale{0.8}\tilingroundedfalse\input{tikz-files/small-class-T6}\tilingroundedtrue};
\node[tilinglabel] at (T6ref1.north west) {$\TT_6$};
\node[anchor=north] (T15) at (-1, -13) {\def\TilingScale{0.8}\tilingroundedfalse\input{tikz-files/small-class-T15}\tilingroundedtrue};
\node[tilinglabel] at (T15.north west) {$\TT_{15}$};

\draw (T14.south) to[treeedge] (T6ref1.north);
\draw (T14.south) to[treeedge] (T15.north);

\node[anchor=north] (T16) at (-6, -17) {\def\TilingScale{0.8}\tilingroundedfalse\input{tikz-files/small-class-T16}\tilingroundedtrue};
\node[tilinglabel] at (T16.north west) {$\TT_{16}$};
\node[anchor=north] (T17) at (-1, -17) {\def\TilingScale{0.8}\tilingroundedfalse\input{tikz-files/small-class-T17}\tilingroundedtrue};
\node[tilinglabel] at (T17.north west) {$\TT_{17}$};
\node at (2.05, -19.3) {$\cong$};
\node[anchor=north] (T18) at (5.42, -16.2) {\def\TilingScale{0.8}\tilingroundedfalse\input{tikz-files/small-class-T18}\tilingroundedtrue};
\node[tilinglabel] at (T18.north west) {$\TT_{18}$};

\draw (T15.south) to[treeedge] (T16.north);
\draw (T15.south) to[treeedge] (T17.north);

\node[anchor=north] (T7ref) at (-8, -21) {\def\TilingScale{0.8}\tilingroundedfalse\input{tikz-files/small-class-T7}\tilingroundedtrue};
\node[tilinglabel] at (T7ref.north west) {$\TT_7$};
\node[anchor=north] (T12ref) at (-5, -21) {\def\TilingScale{0.8}\tilingroundedfalse\input{tikz-files/small-class-T12}\tilingroundedtrue};
\node[tilinglabel] at (T12ref.north west) {$\TT_{12}$};

\draw (T16.south) to[treeedge] (T7ref.north);
\draw (T16.south) to[treeedge] (T12ref.north);

\node[anchor=north] (T6ref2) at (0, -23) {\def\TilingScale{0.8}\tilingroundedfalse\input{tikz-files/small-class-T6}\tilingroundedtrue};
\node[tilinglabel] at (T6ref2.north west) {$\TT_6$};
\node[anchor=north] (T19) at (6, -23) {\def\TilingScale{0.8}\tilingroundedfalse\input{tikz-files/small-class-T19}\tilingroundedtrue};
\node[tilinglabel] at (T19.north west) {$\TT_{19}$};

\draw (T18.south) to[treeedge] (T6ref2.north);
\draw (T18.south) to[treeedge] (T19.north);

\node[anchor=north] (T15ref) at (6, -28) {\def\TilingScale{0.8}\tilingroundedfalse\input{tikz-files/small-class-T15}\tilingroundedtrue};
\node[tilinglabel] at (T15ref.north west) {$\TT_{15}$};

\draw (T19.south) -- (T15ref.north);
\end{tikzpicture}
}
\caption{The complete proof tree for $\I(3421)$.}
\label{figure:proof-tree}
\end{figure}

	The complete proof tree is shown in Figure~\ref{figure:proof-tree}. Four tilings in this proof tree are not expanded by any rule. In the language of~\cite{combinatorial-exploration}, these are \emph{verified} tilings, meaning their enumeration is independently known. The simplest is $\TT_1$, the empty tiling, with generating function $F_1(x) = 1$. Next is $\TT_3$, a single point cell, with $F_3(x) = x$. The tiling $\TT_6$ is a $2 \times 2$ tiling whose only active cells are the two point cells $(1,2)$ and $(2,1)$, with generating function $F_6(x) = x^2$. Finally, $\TT_7$ is a $1 \times 1$ tiling that avoids $231$ and $312$. Simion and Schmidt~\cite{simion:restricted-permutations} showed that $|\I_n(231)| = 2^{n-1}$ for $n \geq 1$, giving $F_7(x) = (1-x)/(1-2x)$.

\subsubsection{Tracking}

	Each rule in the proof tree translates into an equation for the generating function of the parent in terms of the generating functions of the children. Rules arising from row placement and factoring correspond to sums and products, respectively. The generating functions of tilings related by row and column separation are equal, as separation is a size-preserving bijection. The fusion rule, however, requires more care.

	Recall from Section~\ref{subsection:strategies} that fusion merges two adjacent rows (or columns) into one. The reverse operation, unfusion, splits a row back into two. If a gridded involution on the fused tiling has $k$ entries in the fused row, it can be unfused in $k + 1$ ways, one for each position of the boundary between the two rows. In our proof tree, the fusion $\TT_{19} \to \TT_{15}$ merges rows $1$ and $2$ of $\TT_{19}$ (and simultaneously columns $1$ and $2$) into row $1$ (and column $1$) of $\TT_{15}$. To write the generating function equation for this rule, it is not enough to know the total number of gridded involutions of each size on $\TT_{15}$. Instead, we must know how many gridded involutions there are of size $n$ with $k$ entries in cell $(4,1)$ for all values of $n$ and $k$, as any such gridded involution corresponds to $k+1$ gridded involutions of size $n$ on $\TT_{19}$.

	We call this \emph{tracking} cell $(4,1)$ of $\TT_{15}$. By involution symmetry, the mirror cell $(1,4)$ always has the same number of entries, so we do not need to separately track cell $(1,4)$. Since $\TT_{15}$ is itself expanded by a rule in the proof tree, the generating function of $\TT_{15}$ is determined by the generating functions of its children. To track cell $(4,1)$ of $\TT_{15}$, we must therefore track the corresponding cells on each child tiling. This propagation turns out to be straightforward.

	The rule expanding $\TT_{15}$ is a row placement, which corresponds to a sum of generating functions. The two children are $\TT_{16}$ (the empty-row case) and $\TT_{17}$ (the nonempty case). On $\TT_{16}$, the bottommost row is empty. Cell $(4,1)$ lies in this row, so it contains no entries and there is nothing to track.

	On $\TT_{17}$, the bottommost entry is placed in the tracked cell $(4,1)$ of $\TT_{15}$. After the grid expansion, the tracked cell splits into the point cell $(5,1)$ of $\TT_{17}$, containing the newly placed entry, and cell $(6,2)$, which contains all other entries. After separation and factoring, the tracked point cell becomes one of the two points of $\TT_6$, and the remaining tracked entries become cell $(5,1)$ of $\TT_{19}$.

	The tracking closes neatly at the fusion step. Cell $(5,1)$ of $\TT_{19}$ is one of the two cells that merge under fusion to produce cell $(4,1)$ of $\TT_{15}$, so tracking this cell is exactly what is needed. Unfusing a gridded involution on $\TT_{15}$ with $k$ entries in cell $(4,1)$ distributes those entries between cells $(5,1)$ and $(5,2)$ of $\TT_{19}$, and the number of entries that land in the tracked cell $(5,1)$ is recorded by the tracking variable.

	\subsubsection{Generating functions}

	To extract generating functions from the proof tree, we associate to each tiling $\TT_i$ a generating function $F_i$. For tilings without tracking, $F_i(x) = \sum_{n \geq 0} a_n x^n$, where $a_n$ is the number of gridded involutions of size $n$ on $\TT_i$. For tilings that carry tracking information, we use a bivariate generating function $F_i(x, y) = \sum_{n,k \geq 0} a_{n,k}\, x^n y^k$, where $a_{n,k}$ counts the gridded involutions of size $n$ on $\TT_i$ having exactly $k$ entries in total in the tracked cells. As discussed above, each tracked cell is one representative of a symmetric pair, and by involution symmetry its partner always has the same number of entries, so a single variable $y$ suffices.

	Translating each rule mechanically, the proof tree yields the system of $16$ equations in Figure~\ref{figure:3421-system}, one for each tiling (recall that $\TT_8$ through $\TT_{11}$ are the expansion of the verified tiling $\TT_7$ and are not shown).
	\begin{figure}
	\begin{align*}
		F_0(x) &= F_1(x) + F_2(x) + F_4(x)  &  F_{12}(x) &= F_0(x) + F_{13}(x) \\
		F_1(x) &= 1                         &  F_{13}(x) &= F_{14}(x) \\
		F_2(x) &= F_0(x) \cdot F_3(x)       &  F_{14}(x) &= F_6(x) \cdot F_{15}(x, 1) \\
		F_3(x) &= x                         &  F_{15}(x, y) &= F_{16}(x) + F_{17}(x, y) \\
		F_4(x) &= F_5(x)                    &  F_{16}(x) &= F_7(x) \cdot F_{12}(x) \\
		F_5(x) &= F_6(x) \cdot F_7(x) \cdot F_{12}(x)  &  F_{17}(x, y) &= F_{18}(x, y) \\
		F_6(x) &= x^2                       &  F_{18}(x, y) &= x^2 y \cdot F_{19}(x, y) \\
		F_7(x) &= \frac{1-x}{1-2x}          &  F_{19}(x, y) &= \frac{F_{15}(x, 1) - y\, F_{15}(x, y)}{1 - y}
	\end{align*}
	\caption{The system of $16$ equations arising from the proof tree for $\I(3421)$.}
	\label{figure:3421-system}
	\end{figure}
	The equation for $F_{14}(x)$ uses $F_{15}(x, 1)$ because we do not need the tracking information passed to $\TT_{14}$. The equation for $F_{18}(x,y)$ uses $x^2 y$ rather than $F_6(x)$ because the copy of $\TT_6$ appearing in this factoring carries tracking on one of its point cells. A gridded involution on $\TT_{15}$ of size $n$ with $k$ entries in the tracked cell unfuses to $k+1$ gridded involutions on $\TT_{19}$, one for each number $j \in \{0, \ldots, k\}$ of entries that land in the new tracked cell $(5,1)$. The contribution to $F_{19}(x,y)$ is therefore
	\[
		a_{n,k}\, x^n(1 + y + \cdots + y^k) = a_{n,k}\, x^n \cdot\frac{1-y^{k+1}}{1-y},
	\]
	and summing over all $n$ and $k$ gives the stated expression for $F_{19}(x,y)$.

	In general, the theory of productive strategies developed in~\cite{combinatorial-exploration} guarantees that any combinatorial specification arising from Combinatorial Exploration determines unique counting sequences. The extension of this theory to bivariate systems with tracked cells is the subject of forthcoming work~\cite{bean:fusion}. For the present system, uniqueness can be verified by inspection.
	
	We solve this system by using algebraic elimination methods (e.g., Gr\"obner bases) to extract equations to which the generalized kernel methods of Bousquet-M\'elou and Jehanne~\cite{bousquet-melou:poly-eqs} can be applied. Each such application of the kernel method produces a new equation between the variables and can be added into the system. This process is repeated until the resulting ideal is zero-dimensional and a polynomial in terms of $x$ and $F_0(x)$ can be computed.

\smallclasstheorem*

The exponential growth rate of $\I(3421)$ is the reciprocal of the smallest positive real root of $1 - 2x - x^2 - 2x^3$, which is
\[
\frac{6(71 + 6\sqrt{177})^{1/3}}{(71 + 6\sqrt{177})^{2/3} - (71 + 6\sqrt{177})^{1/3} - 11} \approx 2.6590.
\]

\subsection{Involutions avoiding the pattern 2431}
\label{subsec:I-2431}

The class $\I(2431) = \I(2431, 4132)$ is enumerated using the same Combinatorial Exploration framework and the same four strategies applied in Section~\ref{subsec:I-3421}. The specification we obtained has 38 tilings, too many to display as a single figure. Rather than trace the entire proof tree, we highlight one instance each of row placement, factoring, and fusion, and refer the reader to the interactive HTML rendering\footnote{\url{https://jaypantone.github.io/pattern-avoiding-involutions/av2431_4132.html}} of the full specification, where every tiling, obstruction, and rule can be inspected directly. As with $\I(3421)$, the specification is also available as an arXiv ancillary file and in the Zenodo dataset~\cite{bean:invs-data-set}.

\begin{figure}
	\centering
	\includegraphics[width=\linewidth]{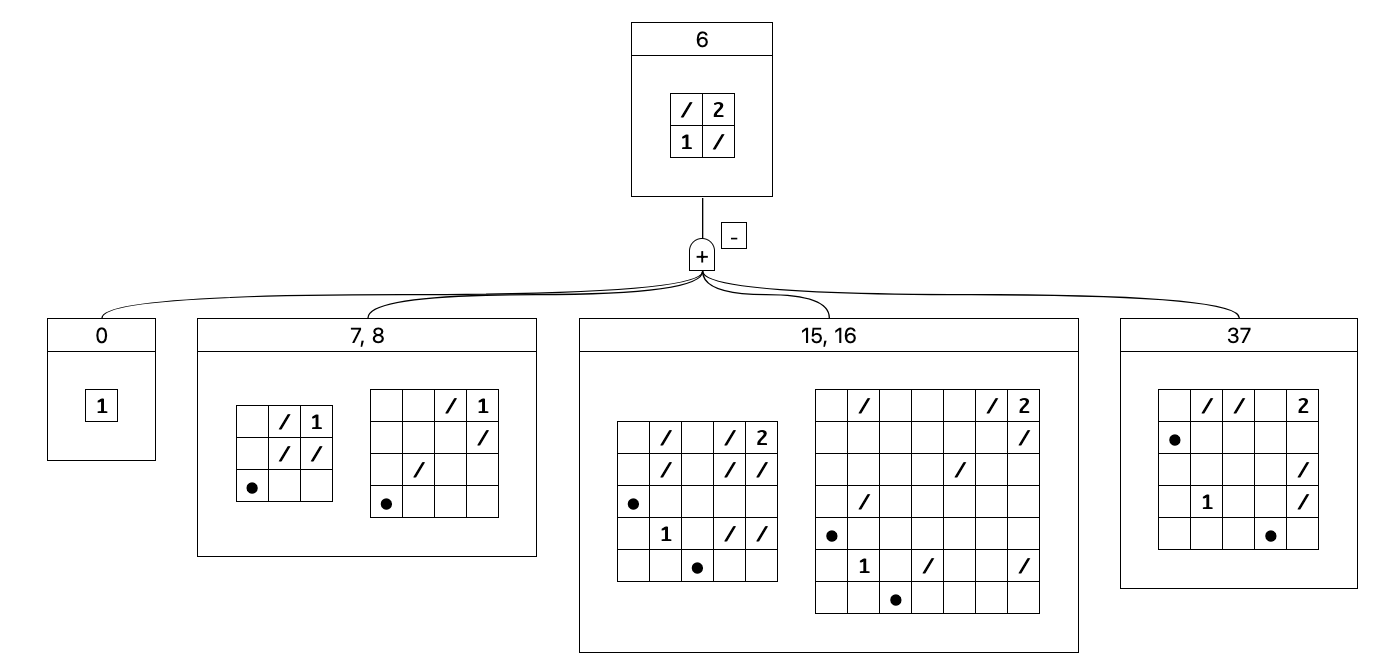}
	\caption{Row placement applied to $\TT_6$.}
	\label{figure:2431-row-placement}
\end{figure}

\textbf{Row placement.} Figure~\ref{figure:2431-row-placement} shows the row placement rule that expands the tiling $\TT_6$, a $2 \times 2$ tiling with an $\Av(132)$ cell at $(1,1)$, an $\Av(2431, 4132)$ cell at $(2,2)$, and monotone cells at $(2,1)$ and $(1,2)$. Row placement targets the bottom row and produces four children.

The figures in this section are screenshots of the HTML rendering of the specification, and we now make a few remarks on how to read them. Each nonempty cell of a tiling is annotated with a symbol: a solid dot marks a point cell, a diagonal line marks a monotone cell (a forward slash for $\Av(21)$ and a backslash for $\Av(12)$), and positive integers index the remaining local classes. The correspondence between integer labels and local obstructions is given by a legend specific to each tiling that can be viewed by clicking the tiling in the HTML version. The same integer may denote different classes across different tilings. Crossing obstructions are not displayed. They can also be inspected by clicking on any tiling in the HTML version. When two tilings are related by row and column separation, we display them together in a single panel, as with $\TT_7 \cong \TT_8$ and $\TT_{15} \cong \TT_{16}$ in Figure~\ref{figure:2431-row-placement}. Each rule appears as a small labeled node between a parent tiling and its children: $+$ marks a disjoint union, as produced by row placement, $\times$ marks a Cartesian product, as produced by factoring, and $?$ marks a fusion. A tiling filled in green is verified, meaning that its enumeration is already known. The boxed minus sign is used in the HTML to collapse or expand a rule.

The four children of $\TT_6$ correspond to four structural cases for row $1$.

The first child $\TT_0$ corresponds to the case that row $1$ is empty: column $1$ is empty by symmetry, and the remaining $\Av(2431, 4132)$ cell is a copy of the root.

The second child $\TT_7 \cong \TT_8$ corresponds to the case that the bottommost entry is a fixed point at the diagonal cell $(1,1)$. In the legend for $\TT_7$, the integer $1$ denotes $\Av(2431, 4132)$, so this tiling consists of a point cell at $(1,1)$, the class $\Av(2431, 4132)$ at $(3,3)$, and three monotone cells.

The third child $\TT_{15} \cong \TT_{16}$ corresponds to the case that the bottommost entry lies in the diagonal cell $(1,1)$ but is not a fixed point. By involution symmetry, its mirror image is simultaneously placed as the leftmost entry of column $1$. The mirror pair appears as the two points at cells $(1,3)$ and $(3,1)$ of $\TT_{15}$. In the legend for $\TT_{15}$, the integer $1$ denotes $\Av(132)$, appearing at cell $(2,2)$, and the integer $2$ denotes $\Av(2431, 4132)$, appearing at cell $(5,5)$. The separation $\TT_{15} \cong \TT_{16}$ separates row $4$ into three rows due to the three crossing $21$ obstructions between cells $(2,4)$, $(4,4)$, and $(5,4)$. By involution symmetry, column $4$ splits into three columns.

The fourth child $\TT_{37}$ corresponds to the case that the bottommost entry lies in the off-diagonal cell $(2,1)$. Being off the diagonal, this entry cannot be a fixed point. Its mirror image is simultaneously placed in cell $(1,2)$, producing the mirror pair at cells $(4,1)$ and $(1,4)$. The legend for $\TT_{37}$ matches that of $\TT_{15}$: the integer $1$ is $\Av(132)$ at cell $(2,2)$ and the integer $2$ is $\Av(2431, 4132)$ at cell $(5,5)$.

A few crossing obstructions give a flavor of what the figure suppresses. The tiling $\TT_7$ carries the crossing obstruction $(21, ((2,2),(3,2)))$, which forces every entry of cell $(2,2)$ to lie below every entry of cell $(3,2)$, as well as the crossing obstruction $(2431, ((2,3),(3,3),(3,3),(3,3)))$, a manifestation of the basis pattern $2431$ spread across cells. In $\TT_{37}$, the crossing obstruction $(132, ((2,2),(2,2),(5,2)))$ forbids a $132$ pattern formed from two entries of the $\Av(132)$ cell at $(2,2)$ and one entry of the monotone cell at $(5,2)$, a constraint stronger than what the local classes alone would impose.

\begin{figure}
	\centering
	\includegraphics[width=0.5\linewidth]{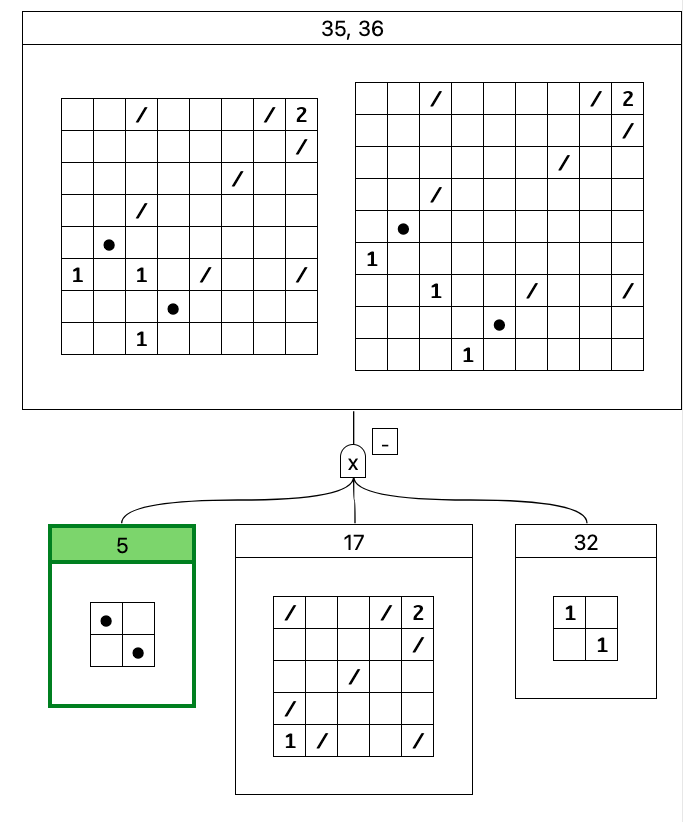}
	\caption{Factoring applied to $\TT_{36}$, which is derived from $\TT_{35}$ by row and column separation.}
	\label{figure:2431-factor}
\end{figure}

\textbf{Factoring.} Figure~\ref{figure:2431-factor} shows a factoring rule applied to $\TT_{36}$, which itself is obtained from $\TT_{35}$ by row and column separation. Its active cells partition into three involution-symmetric subsets with no crossing obstructions between them, producing the factors $\TT_5$, $\TT_{17}$, and $\TT_{32}$.

The factor $\TT_5$ is a $2 \times 2$ tiling containing a mirror pair of points, and $\TT_{32}$ is a $2 \times 2$ tiling containing a mirror pair of $\Av(132)$ cells. Both $\TT_{17}$ and $\TT_{32}$ are expanded further elsewhere in the proof tree.

\begin{figure}
	\centering
	\includegraphics[width=0.15\linewidth]{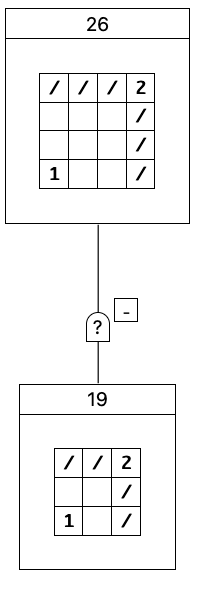}
	\caption{Fusion applied to $\TT_{26}$.}
	\label{figure:2431-fusion}
\end{figure}

\textbf{Fusion.} Figure~\ref{figure:2431-fusion} shows a fusion rule that fuses $\TT_{26}$, which is a $4 \times 4$ tiling with an $\Av(132)$ cell at $(1,1)$, an $\Av(2431, 4132)$ cell at $(4,4)$, and monotone cells filling the rest of row $4$ and column $4$. Rows $2$ and $3$ are indistinguishable. They fuse into a single row, and by involution symmetry columns $2$ and $3$ fuse as well, producing $\TT_{19}$.

\subsubsection{Generating functions}

Each of the 38 tilings in the specification contributes one equation, in the manner of Section~\ref{subsec:I-3421}, and fusion rules again require tracking. The full system is given in Figure~\ref{figure:2431-system}.

\begin{figure}
\begin{align*}
		F_0(x) &= F_1(x) + F_2(x) + F_4(x)  &  F_{19}(x,y) &= F_{12}(x,y) + F_{20}(x,y) + F_{22}(x,y) + F_{25}(x,y)\\
		F_1(x) &= 1  &  F_{20}(x,y) &= F_{21}(x,y)\\
		F_2(x) &= F_3(x) \cdot F_0(x)  &  F_{21}(x,y) &= F_3(x) \cdot F_9(x) \cdot F_{14}(x,y)\\
		F_3(x) &= x  &  F_{22}(x,y) &= F_{23}(x,y)\\
		F_4(x) &= F_5(x) \cdot F_6(x)  &  F_{23}(x,y) &= F_5(x) \cdot F_{24}(x,y)\\
		F_5(x) &= x^2  &  F_{24}(x,y) &= \frac{F_{17}(x,1)-y\,F_{17}(x,y)}{1-y}\\
		F_6(x) &= F_0(x) + F_7(x) + F_{15}(x) + F_{37}(x)  &  F_{25}(x,y) &= F_5(x) \cdot F_{26}(x,y)\\
		F_7(x) &= F_8(x)  &  F_{26}(x,y) &= \frac{F_{19}(x,1)-y\,F_{19}(x,y)}{1-y}\\
		F_8(x) &= F_3(x) \cdot F_9(x) \cdot F_{12}(x,1)  &  F_{27}(x,y) &= F_3(x) \cdot F_9(x) \cdot F_{12}(x,y) \cdot F_{28}(x)\\
		F_9(x) &= F_1(x) + F_{10}(x) + F_{11}(x)  &  F_{28}(x) &= F_1(x) + F_{29}(x) + F_{30}(x)\\
		F_{10}(x) &= F_3(x) \cdot F_9(x)  &  F_{29}(x) &= F_3(x) \cdot F_{28}(x)\\
		F_{11}(x) &= 0  &  F_{30}(x) &= F_{31}(x)\\
		F_{12}(x,y) &= F_0(x) + F_{13}(x,y)  &  F_{31}(x) &= F_5(x) \cdot F_{32}(x) \cdot F_{28}(x)\\
		F_{13}(x,y) &= x^2y \cdot F_{14}(x,y)  &  F_{32}(x) &= F_1(x) + F_{33}(x)\\
		F_{14}(x,y) &= \frac{F_{12}(x,1)-y\,F_{12}(x,y)}{1-y}  &  F_{33}(x) &= F_{34}(x)\\
		F_{15}(x) &= F_{16}(x)  &  F_{34}(x) &= F_5(x) \cdot F_{32}(x)^2\\
		F_{16}(x) &= F_5(x) \cdot F_{17}(x,1)  &  F_{35}(x,y) &= F_{36}(x,y)\\
		F_{17}(x,y) &= F_{18}(x,y) + F_{35}(x,y)  &  F_{36}(x,y) &= F_5(x) \cdot F_{17}(x,y) \cdot F_{32}(x)\\
		F_{18}(x,y) &= F_{19}(x,y) + F_{27}(x,y)  &  F_{37}(x) &= F_5(x) \cdot F_{19}(x,1)
	\end{align*}
\caption{The system of $38$ equations arising from the specification for $\I(2431)$.}
\label{figure:2431-system}
\end{figure}

As before, we solve this system using the kernel method combined with algebraic elimination. The result is the following.

\bigclasstheorem*

The exponential growth rate of $\I(2431)$ is the reciprocal of the smallest positive real root of $1-2x-6x^2+8x^3+8x^4+4x^5$, approximately $2.53041$.

\section{Enumerating 1324-avoiding involutions and 4231-avoiding involutions} 

\label{section:mosaic}

\subsection{The Mosaic method for involutions}

Section~\ref{section:two-exact-enumerations} gave exact enumerations for $\I(3421)$ and $\I(2431)$, leaving the Wilf classes $\I(1324)$ and $\I(4231)$. Exact enumeration for these two remains out of reach, so we instead adapt the \emph{Mosaic method}, a forthcoming counting algorithm of Bean and Pantone~\cite{bean:mosaic-method}, to the involution setting by using tilings. The Mosaic method produces many terms of the two counting sequences without producing a generating function.

The starting point is the insertion encoding of Albert, Linton, and Ru\v{s}kuc~\cite{albert:insertion-encoding}, as developed algorithmically by Vatter~\cite{vatter:regular-insertion-encoding}. In that encoding, a partial permutation is represented by a configuration with \emph{slots}, each of which marks an interval into which future entries must be inserted. One builds the permutation by inserting its entries one at a time in increasing order of value, so that each new entry is larger than all those already placed, with each insertion either filling a slot or replacing it with new slots. In the involution setting, an insertion that is not a fixed point commits the new entry together with its mirror image, so one step may commit two entries. If a class admits only finitely many valid slot configurations up to isomorphism, these transitions form a finite automaton, the language of insertion histories is regular, and the generating function is rational. The classes we target here do not have this finiteness property, so we truncate: fix $N \in \N$ and count only the involutions of length at most $N$. With the length bound in place, the reachable state space is finite.

The truncation gives a simple pruning rule. Let $k$ be the number of entries committed so far, and let $s$ be the number of diagonal slots plus twice the number of off-diagonal mirror pairs of slots, which is a lower bound on the number of additional entries needed to fill them. If $k + s > N$, no completion of the state has length at most $N$, so the state is discarded.

The tiling representation removes two bookkeeping costs that are expensive in a direct implementation of the insertion encoding. First, the dynamic program must decide when two states are isomorphic, so that their counts can be merged. Second, it must decide whether a state $\TT$ is completable, meaning that $\GridI(\TT)$ is nonempty. For tilings, canonicalization handles the first task and the obstruction and requirement machinery of Combinatorial Exploration handles the second. Conway, Guttmann, and Zinn-Justin~\cite{conway:4231-50-terms} accomplished much the same for $\Av(1324)$ using a state description tailored to that class. Tilings provide this bookkeeping without requiring a new state model for each class.

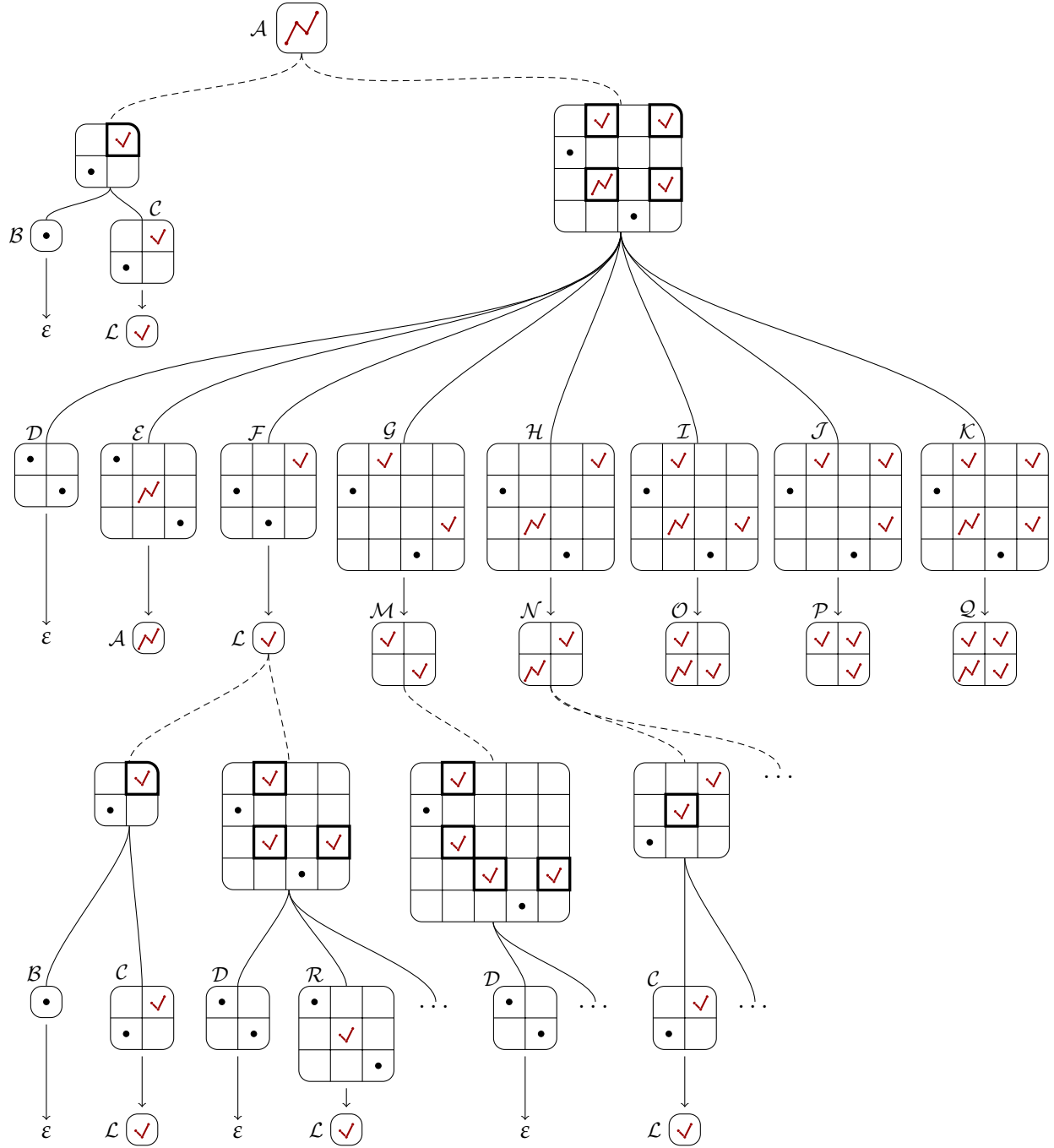
\begin{figure}
	\centering
	\begin{tikzpicture}[
		every node/.style={inner sep=0pt, outer sep=0pt},
		tlabel/.style={font=\small, anchor=south east, xshift=-2pt, yshift=-2pt, overlay},
		leaflab/.style={font=\scriptsize, anchor=north, yshift=-3pt, align=center},
		contlab/.style={font=\scriptsize, anchor=north, yshift=-3pt, align=center},
		prunedlab/.style={font=\scriptsize, anchor=north, yshift=-3pt, align=center, text=red!60!black},
		treeedge/.style={thin, out=-90, in=90, looseness=0.5},
	]
	
		\node (A) at (-2, 0) {%
			\tiling{0.8}{1}{1}{}%
			{{4/{(0.2, 0.2), (0.4, 0.6), (0.6, 0.4), (0.8, 0.8)}}}%
			{}%
		};
		
		\begin{scope}[shift={(-5,-2)}]
			\node (X) at (0, 0){
				\tiling{0.5}{2}{2}{0/0}{3/{(1.3, 1.5), (1.5, 1.3), (1.7, 1.7)}}{}
			};
			\draw[line width=1.25pt] (-0.05,0) -- (0.45,0) [rounded corners=5pt] --  (0.45, 0.5)  [rounded corners=0pt] -- (-0.05, 0.5) -- cycle;
		\end{scope}
		
		\begin{scope}[shift={(3,-2.2)}]
			\node (Y) at (0, 0){
				\tiling{0.5}{4}{4}{0/2, 2/0}
				{
					4/{(1.2, 1.2), (1.4, 1.6), (1.6, 1.4), (1.8, 1.8)},
					3/{(1.3, 3.5), (1.5, 3.3), (1.7, 3.7)},
					3/{(3.3, 3.5), (3.5, 3.3), (3.7, 3.7)},
					3/{(3.3, 1.5), (3.5, 1.3), (3.7, 1.7)}%
				}
				{}
			};
			\draw[line width=1.25pt] (0.45,0.5) -- (0.95,0.5) [rounded corners=5pt] --  (0.95, 1)  [rounded corners=0pt] -- (0.45, 1) -- cycle;
			\draw[line width=1.25pt] (-.55,-0.5) rectangle (-0.05, 0);
			\draw[line width=1.25pt] (.45,-0.5) rectangle (0.95, 0);
			\draw[line width=1.25pt] (-.55,0.5) rectangle (-0.05, 1);
		\end{scope}

		\node[anchor=north] (B) at (-6, -3) {\tiling{0.5}{1}{1}{0/0}{}{}};
		\node[anchor=north] (Bnext) at (-6, -4.65) {$\varepsilon$};

		\node[anchor=north] (C) at (-4.5, -3) {\tiling{0.5}{2}{2}{0/0}%
			{{3/{(1.3, 1.5), (1.5, 1.3), (1.7, 1.7)}}}%
			{}%
		};
		\node[anchor=north] (Cnext) at (-4.5, -4.5) {\tiling{0.5}{1}{1}{}%
			{{3/{(0.3, 0.5), (0.5, 0.3), (0.7, 0.7)}}}%
			{}%
			};

		\node[anchor=north] (D) at (-6, -6.5) {\tiling{0.5}{2}{2}{0/1, 1/0}{}{}};
		\node[anchor=north] (Dnext) at (-6, -9.45) {$\varepsilon$};

		\node[anchor=north] (E) at (-4.4, -6.5) {\tiling{0.5}{3}{3}{0/2, 2/0}%
			{{4/{(1.2, 1.2), (1.4, 1.6), (1.6, 1.4), (1.8, 1.8)}}}%
			{}%
		};
		\node[anchor=north] (Enext) at (-4.4, -9.3) {\tiling{0.5}{1}{1}{}%
			{{4/{(0.2, 0.2), (0.4, 0.6), (0.6, 0.4), (0.8, 0.8)}}}%
			{}%
			};

		\node[anchor=north] (F) at (-2.52, -6.5) {\tiling{0.5}{3}{3}{0/1, 1/0}%
			{{3/{(2.3, 2.5), (2.5, 2.3), (2.7, 2.7)}}}%
			{}%
		};
		\node[anchor=north] (Fnext) at (-2.52, -9.3) {\tiling{0.5}{1}{1}{}%
			{{3/{(0.3, 0.5), (0.5, 0.3), (0.7, 0.7)}}}%
			{}%
			};
		\begin{scope}[shift={(-4.7,-12)}]
			\node (Fnn) at (0, 0){
				\tiling{0.5}{2}{2}{0/0}{3/{(1.3, 1.5), (1.5, 1.3), (1.7, 1.7)}}{}
			};
			\draw[line width=1.25pt] (-0.05,0) -- (0.45,0) [rounded corners=5pt] --  (0.45, 0.5)  [rounded corners=0pt] -- (-0.05, 0.5) -- cycle;
		\end{scope}
		\begin{scope}[shift={(-2.2,-12.5)}]
			\node (Fnnn) at (0, 0){
				\tiling{0.5}{4}{4}{0/2, 2/0}{
					3/{(1.3, 1.5), (1.5, 1.3), (1.7, 1.7)},
					3/{(1.3, 3.5), (1.5, 3.3), (1.7, 3.7)},
					3/{(3.3, 1.5), (3.5, 1.3), (3.7, 1.7)}%
					}{}
			};
			\draw[line width=1.25pt] (-.55,-0.5) rectangle (-0.05, 0);
			\draw[line width=1.25pt] (-.55,0.5) rectangle (-0.05, 1);
			\draw[line width=1.25pt] (0.45,-0.5) rectangle (0.95, 0);
		\end{scope}
		
		\begin{scope}[shift={(0,-15)}]
			\node[anchor=north] (FB) at (-6, 0) {\tiling{0.5}{1}{1}{0/0}{}{}};
			\node[anchor=north] (FBnext) at (-6, -2.2) {$\varepsilon$};
			
			\node[anchor=north] (FC) at (-4.5, 0) {\tiling{0.5}{2}{2}{0/0}%
				{{3/{(1.3, 1.5), (1.5, 1.3), (1.7, 1.7)}}}%
				{}%
			};
			\node[anchor=north] (FCnext) at (-4.5, -2) {\tiling{0.5}{1}{1}{}%
				{{3/{(0.3, 0.5), (0.5, 0.3), (0.7, 0.7)}}}%
				{}%
				};
		\end{scope}
		
		\node[anchor=north] (Fnnn1) at (-3, -15) {\tiling{0.5}{2}{2}{0/1, 1/0}{}{}};
		\node[anchor=north] (Fnnn1next) at (-3, -17.2) {$\varepsilon$};
		\node[anchor=north] (Fnnn2) at (-1.3, -15) {\tiling{0.5}{3}{3}{0/2, 2/0}%
			{{3/{(1.3, 1.5), (1.5, 1.3), (1.7, 1.7)}}}%
			{}%
		};
		\node[anchor=north] (Fnnn2next) at (-1.3, -17) {\tiling{0.5}{1}{1}{}%
			{{3/{(0.3, 0.5), (0.5, 0.3), (0.7, 0.7)}}}%
			{}%
			};
		\node[anchor=north] (Fnnn3) at (0.1, -15.25) {$\cdots$};

		\node[anchor=north] (G) at (-0.4, -6.5) {\tiling{0.5}{4}{4}{0/2, 2/0}%
			{%
				{3/{(1.3, 3.5), (1.5, 3.3), (1.7, 3.7)}},%
				{3/{(3.3, 1.5), (3.5, 1.3), (3.7, 1.7)}}%
			}%
			{}
		};
		\node[anchor=north] (Gnext) at (-0.4, -9.3) {\tiling{0.5}{2}{2}{}%
			{{3/{(0.3, 1.5), (0.5, 1.3), (0.7, 1.7)}},
			{3/{(1.3, 0.5), (1.5, 0.3), (1.7, 0.7)}}}%
			{}%
			};
		\begin{scope}[shift={(1, -11.5)}]
			\node[anchor=north] (Gnn) at (0,0) {\tiling{0.5}{5}{5}{0/3, 3/0}%
				{%
					{3/{(1.3, 2.5), (1.5, 2.3), (1.7, 2.7)}},
					{3/{(1.3, 4.5), (1.5, 4.3), (1.7, 4.7)}},
					{3/{(2.3, 1.5), (2.5, 1.3), (2.7, 1.7)}},
					{3/{(4.3, 1.5), (4.5, 1.3), (4.7, 1.7)}}%
				}%
				{}
			};
			\draw[line width=1.25pt] (-.8,-0.5) rectangle (-0.3, 0);
			\draw[line width=1.25pt] (-.8,-1.5) rectangle (-0.3, -1);
			\draw[line width=1.25pt] (-.3,-2) rectangle (0.2, -1.5);
			\draw[line width=1.25pt] (0.7,-2) rectangle (1.2, -1.5);
		\end{scope}
		
		\node[anchor=north] (Gnnn1) at (1.5, -15) {\tiling{0.5}{2}{2}{0/1, 1/0}{}{}};
		\node[anchor=north] (Gnnn1next) at (1.5, -17.2) {$\varepsilon$};
		\node[anchor=north] (Gnnn2) at (2.6, -15.25) {$\cdots$};

		\node[anchor=north] (H) at (1.9, -6.5) {\tiling{0.5}{4}{4}{0/2, 2/0}%
			{%
				{4/{(1.2, 1.2), (1.4, 1.6), (1.6, 1.4), (1.8, 1.8)}},%
				{3/{(3.3, 3.5), (3.5, 3.3), (3.7, 3.7)}}%
			}%
			{}%
		};
		\node[anchor=north] (Hnext) at (1.9, -9.3) {\tiling{0.5}{2}{2}{}%
			{{4/{(0.2, 0.2), (0.4, 0.6), (0.6, 0.4),(0.8,0.8)}},
			{3/{(1.3, 1.5), (1.5, 1.3), (1.7, 1.7)}}}%
			{}%
			};
		\begin{scope}[shift={(4, -11.5)}]
			\node[anchor=north] (Hnn) at (0,0) {\tiling{0.5}{3}{3}{0/0}%
				{%
					{3/{(1.3, 1.5), (1.5, 1.3), (1.7, 1.7)}},
					{3/{(2.3, 2.5), (2.5, 2.3), (2.7, 2.7)}}%
				}%
				{}
			};
			\draw[line width=1.25pt] (-.3,-0.5) rectangle (0.2, -1);
		\end{scope}
		
		\node[anchor=north] (Hnnn) at (5.5, -11.65) {$\cdots$};
		
		\node[anchor=north] (Hnn1) at (4, -15) {\tiling{0.5}{2}{2}{0/0}%
				{{3/{(1.3, 1.5), (1.5, 1.3), (1.7, 1.7)}}}%
				{}%
			};
		\node[anchor=north] (Hnn1next) at (4, -17) {\tiling{0.5}{1}{1}{}%
			{{3/{(0.3, 0.5), (0.5, 0.3), (0.7, 0.7)}}}%
			{}%
		};
		\node[anchor=north] (Hnn2) at (5.1, -15.25) {$\cdots$};

		\node[anchor=north] (I) at (4.2, -6.5) {\tiling{0.5}{4}{4}{0/2, 2/0}%
			{%
				{4/{(1.2, 1.2), (1.4, 1.6), (1.6, 1.4), (1.8, 1.8)}},%
				{3/{(1.3, 3.5), (1.5, 3.3), (1.7, 3.7)}},%
				{3/{(3.3, 1.5), (3.5, 1.3), (3.7, 1.7)}}%
			}%
			{}
		};
		\node[anchor=north] (Inext) at (4.2, -9.3) {\tiling{0.5}{2}{2}{}%
			{{4/{(0.2, 0.2), (0.4, 0.6), (0.6, 0.4),(0.8,0.8)}},
			{3/{(0.3, 1.5), (0.5, 1.3), (0.7, 1.7)}},
			{3/{(1.3, 0.5), (1.5, 0.3), (1.7, 0.7)}}}%
			{}%
			};

		\node[anchor=north] (J) at (6.4, -6.5) {\tiling{0.5}{4}{4}{0/2, 2/0}%
			{%
				{3/{(3.3, 3.5), (3.5, 3.3), (3.7, 3.7)}},%
				{3/{(1.3, 3.5), (1.5, 3.3), (1.7, 3.7)}},%
				{3/{(3.3, 1.5), (3.5, 1.3), (3.7, 1.7)}}%
			}%
			{}%
		};
		\node[anchor=north] (Jnext) at (6.4, -9.3) {\tiling{0.5}{2}{2}{}%
			{{3/{(1.3, 1.5), (1.5, 1.3), (1.7, 1.7)}},
			{3/{(0.3, 1.5), (0.5, 1.3), (0.7, 1.7)}},
			{3/{(1.3, 0.5), (1.5, 0.3), (1.7, 0.7)}}}%
			{}%
			};

		\node[anchor=north] (K) at (8.7, -6.5) {\tiling{0.5}{4}{4}{0/2, 2/0}%
			{%
				{4/{(1.2, 1.2), (1.4, 1.6), (1.6, 1.4), (1.8, 1.8)}},
				{3/{(1.3, 3.5), (1.5, 3.3), (1.7, 3.7)}},%
				{3/{(3.3, 1.5), (3.5, 1.3), (3.7, 1.7)}},%
				{3/{(3.3, 3.5), (3.5, 3.3), (3.7, 3.7)}}%
			}%
			{}%
		};
		\node[anchor=north] (Knext) at (8.7, -9.3) {\tiling{0.5}{2}{2}{}%
			{{4/{(0.2, 0.2), (0.4, 0.6), (0.6, 0.4),(0.8,0.8)}},
			{3/{(0.3, 1.5), (0.5, 1.3), (0.7, 1.7)}},
			{3/{(1.3, 0.5), (1.5, 0.3), (1.7, 0.7)}},
			{3/{(1.3, 1.5), (1.5, 1.3), (1.7, 1.7)}}}%
			{}%
			};

		\draw[treeedge, densely dashed] (A.south) to (X.north);
		\draw[treeedge, densely dashed] (A.south) to (Y.north);
		\draw[treeedge] (X.south) to (B.north);
		\draw[treeedge] (X.south) to (C.north);
		\draw[treeedge] (Y.south) to (D.north);
		\draw[treeedge] (Y.south) to (E.north);
		\draw[treeedge] (Y.south) to (F.north);
		\draw[treeedge] (Y.south) to (G.north);
		\draw[treeedge] (Y.south) to (H.north);
		\draw[treeedge] (Y.south) to (I.north);
		\draw[treeedge] (Y.south) to (J.north);
		\draw[treeedge] (Y.south) to (K.north);
		\draw[treeedge, densely dashed] (Gnext.south) to (Gnn.north);
		\draw[treeedge, densely dashed] (Hnext.south) to (Hnn.north);
		\draw[treeedge, densely dashed] (Hnext.south) to (Hnnn.north);
		\draw[treeedge, densely dashed] (Fnext.south) to (Fnn.north);
		\draw[treeedge, densely dashed] (Fnext.south) to (Fnnn.north);
		\draw[treeedge] (Fnn.south) to (FB.north);
		\draw[treeedge] (Fnn.south) to (FC.north);
		\draw[treeedge] (Fnnn.south) to (Fnnn1.north);
		\draw[treeedge] (Fnnn.south) to (Fnnn2.north);
		\draw[treeedge] (Fnnn.south) to (Fnnn3.north);
		\draw[treeedge] (Gnn.south) to (Gnnn1.north);
		\draw[treeedge] (Gnn.south) to (Gnnn2.north);
		\draw[treeedge] (Hnn.south) to (Hnn1.north);
		\draw[treeedge] (Hnn.south) to (Hnn2.north);
		
		\draw[treeedge, ->] ($ (B.south) + (0,-0.1) $) to ($ (Bnext.north) + (0,0.1) $);
		\draw[treeedge, ->] ($ (C.south) + (0,-0.1) $) to ($ (Cnext.north) + (0,0.1) $);
		\draw[treeedge, ->] ($ (D.south) + (0,-0.1) $) to ($ (Dnext.north) + (0,0.1) $);
		\draw[treeedge, ->] ($ (E.south) + (0,-0.1) $) to ($ (Enext.north) + (0,0.1) $);
		\draw[treeedge, ->] ($ (F.south) + (0,-0.1) $) to ($ (Fnext.north) + (0,0.1) $);
		\draw[treeedge, ->] ($ (G.south) + (0,-0.1) $) to ($ (Gnext.north) + (0,0.1) $);
		\draw[treeedge, ->] ($ (H.south) + (0,-0.1) $) to ($ (Hnext.north) + (0,0.1) $);
		\draw[treeedge, ->] ($ (I.south) + (0,-0.1) $) to ($ (Inext.north) + (0,0.1) $);
		\draw[treeedge, ->] ($ (J.south) + (0,-0.1) $) to ($ (Jnext.north) + (0,0.1) $);
		\draw[treeedge, ->] ($ (K.south) + (0,-0.1) $) to ($ (Knext.north) + (0,0.1) $);
		\draw[treeedge, ->] ($ (FB.south) + (0,-0.1) $) to ($ (FBnext.north) + (0,0.1) $);
		\draw[treeedge, ->] ($ (FC.south) + (0,-0.1) $) to ($ (FCnext.north) + (0,0.1) $);
		\draw[treeedge, ->] ($ (Fnnn1.south) + (0,-0.1) $) to ($ (Fnnn1next.north) + (0,0.1) $);
		\draw[treeedge, ->] ($ (Fnnn2.south) + (0,-0.1) $) to ($ (Fnnn2next.north) + (0,0.1) $);
		\draw[treeedge, ->] ($ (Gnnn1.south) + (0,-0.1) $) to ($ (Gnnn1next.north) + (0,0.1) $);
		\draw[treeedge, ->] ($ (Hnn1.south) + (0,-0.1) $) to ($ (Hnn1next.north) + (0,0.1) $);

		\node[left=4pt] at (A.west) {\footnotesize$\mathcal{A}$};
		\node[left=3pt] at (B.west) {\footnotesize$\mathcal{B}$};
		\node[above left=2pt and -9pt] at (C.north) {\footnotesize$\mathcal{C}$};
		\node[above left=2pt] at (D.north) {\footnotesize $\mathcal{D}$};
		\node[above left=2pt] at (E.north) {\footnotesize $\mathcal{E}$};
		\node[above left=2pt] at (F.north) {\footnotesize $\mathcal{F}$};
		\node[above left=2pt and 3pt] at (G.north) {\footnotesize $\mathcal{G}$};
		\node[above left=2pt and 3pt] at (H.north) {\footnotesize $\mathcal{H}$};
		\node[above left=2pt and 3pt] at (I.north) {\footnotesize $\mathcal{I}$};
		\node[above left=2pt and 4pt] at (J.north) {\footnotesize $\mathcal{J}$};
		\node[above left=2pt and 4pt] at (K.north) {\footnotesize $\mathcal{K}$};
		\node[left=3pt] at (Cnext.west) {\footnotesize $\mathcal{L}$};
		\node[left=3pt] at (Enext.west) {\footnotesize $\mathcal{A}$};
		\node[left=3pt] at (Fnext.west) {\footnotesize $\mathcal{L}$};
		\node[above left=2pt and 4pt] at (Gnext.north) {\footnotesize $\mathcal{M}$};
		\node[above left=2pt and 4pt] at (Hnext.north) {\footnotesize $\mathcal{N}$};
		\node[above left=2pt and 4pt] at (Inext.north) {\footnotesize $\mathcal{O}$};
		\node[above left=2pt and 4pt] at (Jnext.north) {\footnotesize $\mathcal{P}$};
		\node[above left=2pt and 4pt] at (Knext.north) {\footnotesize $\mathcal{Q}$};
		\node[above left=3pt] at (FB.north) {\footnotesize$\mathcal{B}$};
		\node[above left=3pt and 6pt] at (FC.north) {\footnotesize$\mathcal{C}$};
		\node[above left=2pt and 4pt] at (Fnnn1.north) {\footnotesize $\mathcal{D}$};
		\node[above left=2pt and 10pt] at (Fnnn2.north) {\footnotesize $\mathcal{R}$};
		\node[above left=0pt and 11pt] at (Gnnn1.north) {\footnotesize $\mathcal{D}$};
		\node[above left=0pt and 11pt] at (Hnn1.north) {\footnotesize $\mathcal{C}$};
		\node[left=3pt] at (FCnext.west) {\footnotesize $\mathcal{L}$};
		\node[left=3pt] at (Fnnn2next.west) {\footnotesize $\mathcal{L}$};
		\node[left=3pt] at (Hnn1next.west) {\footnotesize $\mathcal{L}$};

	\end{tikzpicture}%
	\caption{The tiling expansion step to count $1324$-avoiding involutions up to length $4$.}
	\label{figure:1324-tree}
\end{figure}

We explain the Mosaic method by tracing through Figure~\ref{figure:1324-tree}, which computes the number of involutions avoiding $1324$ up to length $4$. The root, labeled $\mc{A}$, represents all involutions avoiding $1324$. The first step is to place the bottommost point. The left branch is the case where the bottommost point is a fixed point, while the right branch is the case where it is not. The cells with bold outlines are newly created by this insertion and, in the style of the original insertion encoding, we decide at this moment whether each such cell is empty or nonempty. The tiling along the left branch has a single cell to decide. The case where it is empty yields tiling $\mc{B}$, while the case where it is nonempty yields tiling $\mc{C}$.\footnote{We do not draw requirements in this figure, but every nonempty cell whose outline is not bold can be assumed to contain a requirement marking it as containing at least one point.} The last processing step for each tiling is to factor away the placed points. Their presence is not required to determine which future insertions are allowed, because as always the obstructions on each tiling are updated to record this information the moment the points are initially placed. The tiling $\mc{B}$ contains only a point, so it becomes the empty tiling $\varepsilon$. The tiling $\mc{C}$ becomes the $1 \times 1$ tiling avoiding $213$, which we label $\mc{L}$.

The right branch has four new cells whose emptiness must be decided, at coordinates $(2,2)$, $(2,4)$, $(4,2)$, and $(4,4)$. We point out here that this tiling has many crossing obstructions that we have hidden, as we have for all tilings in this figure. The cells on the diagonal can be made empty or not independently, but for the mirror image pair $(2,4)$ and $(4,2)$ either both cells must be empty or they both must be nonempty. This leads to eight cases, shown as tilings $\mc{D}$ through $\mc{K}$. Below each we have factored away the points. Note that $\mc{D}$ becomes the empty tiling, $\mc{E}$ becomes the root tiling $\mc{A}$, and $\mc{F}$ becomes another copy of $\mc{L}$, while the rest are new.

What we have described up to this point is one expansion step: take each tiling that has not yet been expanded (in this case, just the root tiling), form one branch for each cell that could contain the bottommost point (in this case, just one cell), allowing for the fact that a point in a diagonal cell may or may not be a fixed point, and then factor away the points. Depending on whether a fixed point was inserted or not, either one point is factored away or two are. For the purposes of counting, this is important information to retain. Therefore, we represent the expansion of a tiling under this process by a rule
\[
	\mc{T} \to (S_1, S_2)
\]
where $\mc{T}$ is the tiling being expanded, $S_1$ and $S_2$ are the multisets of child tilings obtained by factoring away one point and two points, respectively. For the expansion of the root described above, we would record
\[
	\mc{A} \to (\{\varepsilon, \mc{L}\}, \{\varepsilon, \mc{A}, \mc{L},  \mc{M}, \mc{N}, \mc{O}, \mc{P}, \mc{Q}\}).
\]
However, with the foreknowledge that we are only seeking to enumerate 1324-avoiding involutions up to length $4$, we can prune this rule. The expansion from $\mc{A}$ to $\mc{O}$, for example, comes from placing two points. The tiling $\mc{O}$ has three nonempty cells, therefore requiring at least three more points to be placed to fulfill its requirements. Thus, any complete involution derived from $\mc{O}$ will have length at least $5$. As a result, we do not need to expand $\mc{O}$, $\mc{P}$, or $\mc{Q}$.

Tilings $\mc{L}$, $\mc{M}$, and $\mc{N}$ still need to be expanded, so we choose one and repeat the expansion process. Inserting a bottommost fixed point into $\mc{L}$ produces the tilings $\varepsilon$ and $\mc{L}$. Inserting a bottommost non-fixed point into $\mc{L}$ and factoring away the points produces the tilings $\varepsilon$, $\mc{L}$, and two others that we do not draw because they contain two or three nonempty cells. These two tilings can be pruned because we first encountered $\mc{L}$ after placing one point (in the rule from $\mc{A}$), and fulfilling the requirements of these tilings from $\mc{L}$ requires placing two points into $\mc{L}$ and then filling the two or three nonempty cells, leading to involutions of length $5$ or $6$. So, we write the rule
\[
	\mc{L} \to (\{\varepsilon, \mc{L}\}, \{\varepsilon, \mc{L}\}).
\]
Expanding tilings $\mc{M}$ and $\mc{N}$ and pruning appropriately give the rules
\[
	\mc{M} \to (\{\}, \{\varepsilon\})
\]
and
\[
	\mc{N} \to (\{\mc{L}\}, \{\}).
\]
At this point, we have expanded all of the tilings needed in order to count the 1324-avoiding involutions up to length $4$. To produce the actual counts, we use a dynamic programming approach, iterating the four rules
\begin{align*}
	\mc{A} &\to (\{\varepsilon, \mc{L}\}, \{\varepsilon, \mc{A}, \mc{L},  \mc{M}, \mc{N}\})\\
	\mc{L} &\to (\{\varepsilon, \mc{L}\}, \{\varepsilon, \mc{L}\})\\
	\mc{M} &\to (\{\}, \{\varepsilon\})\\
	\mc{N} &\to (\{\mc{L}\}, \{\})
\end{align*}
with additional pruning where appropriate.

\begin{figure}
	\centering
	\begin{tikzpicture}
		\node (A0) at (0,2) {$\mc{A}:1$};

		\draw[dashed] (1,-0.5) -- (1,4.5);

		\node (e1) at (2,4) {$\varepsilon:1$};
		\node (L1) at (2,3) {$\mc{L}:1$};
		
		\draw (A0.east) -- (e1.west);
		\draw (A0.east) -- (L1.west);
		
		\draw[dashed] (3,-0.5) -- (3,4.5);
		
		\node (e2) at (4,4) {$\varepsilon:2$};
		\node (L2) at (4,3) {$\mc{L}:2$};
		\node (A2) at (4,2) {$\mc{A}:1$};
		\node (M2) at (4,0) {$\mc{M}:1$};
		\node (N2) at (4,1) {$\mc{N}:1$};
		
		\draw (A0.east) to[out=60,in=180] ([yshift=1pt]e2.west);
		\draw (A0.east) -- ([yshift=-2pt]L2.west);
		\draw (A0.east) -- (A2.west);
		\draw (A0.east) -- (M2.west);
		\draw (A0.east) -- (N2.west);
		\draw (L1.east) -- ([yshift=-1pt]e2.west);
		\draw (L1.east) -- (L2.west);
		
		\draw[dashed] (5,-0.5) -- (5,4.5);
		
		\node (e3) at (6,4) {$\varepsilon:4$};
		\node (L3) at (6,3) {$\mc{L}:5$};
		
		\draw (L1.east) -- ([yshift=2pt]e3.west);
		\draw (L1.east) to[out=15,in=165] ([yshift=2pt] L3.west);
		\draw (L2.east) -- ([yshift=0pt]e3.west);
		\draw (L2.east) -- (L3.west);
		\draw (A2.east) -- ([yshift=-3pt]e3.west);
		\draw (A2.east) -- ([yshift=-3pt]L3.west);
		\draw (N2.east) -- ([yshift=-6pt]L3.west);

		\draw[dashed] (7,-0.5) -- (7,4.5);
		
		\node (e4) at (8,4) {$\varepsilon:9$};
		\node (L4) at (8,3) {$\color{lightgray}{\mc{L}:8}$};
		\node (A4) at (8,2) {$\color{lightgray}{\mc{A}:1}$};
		\node (M4) at (8,0) {$\color{lightgray}{\mc{M}:1}$};
		\node (N4) at (8,1) {$\color{lightgray}{\mc{N}:1}$};

		\draw (A2.east) to[out=60,in=180] ([yshift=2pt]e4.west);
		\draw[gray] (A2.east) -- ([yshift=-2pt]L4.west);
		\draw[gray] (A2.east) -- (A4.west);
		\draw[gray] (A2.east) -- (M4.west);
		\draw[gray] (A2.east) -- (N4.west);
		\draw (M2.east) -- ([yshift=-1pt]e4.west);
		\draw (L3.east) -- ([yshift=0pt]e4.west);
		\draw[gray] (L3.east) -- (L4.west);
		\draw (L2.east) -- ([yshift=1pt]e4.west);
		\draw[gray] (L2.east) to[out=15,in=165] ([yshift=2pt] L4.west);
		
		\node at (0, 4.7) {\footnotesize length $0$};
		\node at (2, 4.7) {\footnotesize length $1$};
		\node at (4, 4.7) {\footnotesize length $2$};
		\node at (6, 4.7) {\footnotesize length $3$};
		\node at (8, 4.7) {\footnotesize length $4$};

	\end{tikzpicture}
	\caption{A representation of the dynamic programming step to enumerate the $1324$-avoiding involutions up to length $4$.}	
	\label{figure:dyn-prog}
\end{figure}
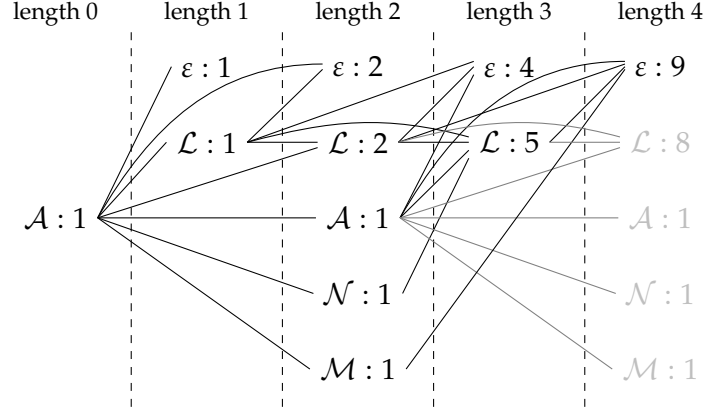

We initialize the counting by recording one occurrence of the length $0$ involution on the root tiling $\mc{A}$, represented by ``$\mc{A}:1$'' in the ``length 0'' column of Figure~\ref{figure:dyn-prog}. Lines in this figure connect each tiling with the children that result from its expansion, either one column to the right when one point is placed or two columns to the right when two points are placed. The counts associated with each label are determined by summing the counts of the tilings in the previous columns to which it is connected. The total number of $1324$-avoiding involutions of length $n$ is the count associated with the tiling $\varepsilon$ in column $n$.

We performed this computation for both $\I(1324)$ and $\I(4231)$ on Marquette University's Raj HPC cluster, producing the counting sequence for $\I(1324)$ through length $34$ and that of $\I(4231)$ through length $40$. The computations ran for several weeks using 60 CPU cores and 300 GB of RAM. In Section~\ref{subsection:asympt-est} we use these sequences to estimate the asymptotic behavior of both classes.

\subsection{Growth rates of involution classes}
\label{subsection:old-results}

In Section~\ref{subsection:asympt-est}, we use the initial terms of the counting sequences of $\I(1324)$ and $\I(4231)$ to empirically estimate their exponential growth rates. Before starting, we review some known results about growth rates of involution classes including a lower bound on the (upper) growth rate of $\I(1324)$. 

It will be helpful here to introduce a few pieces of notation. The \emph{sum} of two permutations $\alpha$ and $\beta$ is the permutation $\alpha \oplus \beta$ formed by putting the entries of $\beta$ above and to the right of the entries of $\alpha$. Formally, if $|\alpha| = k$ and $|\beta| = \ell$, then the $i$th entry of $\alpha \oplus \beta$ is
\[
	(\alpha \oplus \beta)(i) = \left\{\begin{array}{ll}
		\alpha(i), &\text{ if $1 \leq i \leq k$}\\
		\beta(i-k) + k, &\text{ if $k+1 \leq i \leq k+\ell$}
	\end{array}\right..
\]
Similarly, the \emph{skew sum} of $\alpha$ and $\beta$ is the permutation $\alpha \ominus \beta$ formed by putting the entries of $\beta$ below and to the right of the entries of $\alpha$, so that
\[
	(\alpha \ominus \beta)(i) = \left\{\begin{array}{ll}
		\alpha(i)+\ell, &\text{ if $1 \leq i \leq k$}\\
		\beta(i-k), &\text{ if $k+1 \leq i \leq k+\ell$}
	\end{array}\right..
\]

We call a permutation \emph{sum indecomposable} if it cannot be written as the sum of two nonempty permutations and \emph{skew indecomposable} if it cannot be written as the skew sum of two nonempty permutations.

For any permutation class $\CC$, we define the \emph{lower growth rate} $\lgr(\CC)$ and \emph{upper growth rate} $\ugr(\CC)$ by
\[
	\lgr(\CC) = \ds\liminf_{n \to \infty} \sqrt[n]{|\CC_n|} \qquad \text{ and } \qquad \ugr(\CC) = \limsup_{n \to\infty} \sqrt[n]{|\CC_n|}.
\]
When these quantities are equal, we call them the \emph{(proper) growth rate} of $\CC$ and denote it by $\gr(\CC)$. Arratia~\cite{arratia:stanley-wilf} proved that all \emph{principal classes} (those defined by avoiding a single pattern) have growth rates, which are finite by the Marcus--Tardos theorem~\cite{marcus:marcus-tardos}, but it is not known whether this fact extends to all permutation classes avoiding more than one pattern. We define lower, upper, and proper growth rates analogously for involution classes.

B\'ona, Homberger, Pantone, and Vatter~\cite{bona:pattern-avoiding-involutions} proved  the following.
\begin{proposition}[\cite{bona:pattern-avoiding-involutions}, Proposition 2.1]
	\label{proposition:sum-indec-gr}
	If every permutation in $B$ is sum indecomposable, then $\I(B)$ has a proper growth rate.
\end{proposition}
\begin{proposition}[\cite{bona:pattern-avoiding-involutions}, Proposition 2.2]
	\label{proposition:skew-indec-gr-bound}
	If $\beta$ is a skew indecomposable involution, then
	\[
		\ugr(\I(\beta)) \geq \sqrt{\gr(\Av(\beta))}.
	\]
\end{proposition}

Proposition~\ref{proposition:sum-indec-gr} tells us that $\I(4231)$ has a proper growth rate, but tells us nothing about its value, nor does it say anything about $\I(1324)$, as $1324$ is not sum indecomposable. Conversely, Proposition~\ref{proposition:skew-indec-gr-bound} provides a lower bound on $\ugr(\I(1324))$ but also does not tell us that $\I(1324)$ has a proper growth rate, nor does it say anything about $\I(4231)$.

When B\'ona, Homberger, Pantone, and Vatter~\cite{bona:pattern-avoiding-involutions} proved these propositions, the best known bounds on the growth rate of $\Av(1324)$ were a lower bound of $\approx 9.81$ due to Bevan~\cite{bevan:1324-9.81} and an upper bound of $\approx 13.74$ due to B\'ona~\cite{bona:new-record-1324}. Combining that lower bound with Proposition~\ref{proposition:skew-indec-gr-bound} led to the bound
\[
	\ugr(\I(1324)) > 3.13.
\]
In the time since, both the lower and upper bounds for the growth rate of $\Av(1324)$ have been improved. Using a staircase decomposition, Bevan, Brignall, Elvey Price, and Pantone~\cite{bevan:1324-staircase} found a lower bound of $\approx 10.271$ and an upper bound of $13.5$. Using this stronger lower bound, Proposition~\ref{proposition:skew-indec-gr-bound} now tells us that
\[
	\ugr(\I(1324)) > 3.20.
\]
We note finally that Conway, Guttmann, and Zinn-Justin~\cite{conway:4231-50-terms} conjecture that the growth rate of $\Av(1324)$ is $11.600 \pm 0.003$, which would imply $\ugr(\I(1324)) > 3.405$. They further observe that if this growth rate is algebraic, then $9 + 3\sqrt{3}/2 = 11.59807\ldots$ is an appealing candidate, which would imply $\ugr(\I(1324)) > 3.40559$. While Proposition~\ref{proposition:skew-indec-gr-bound} is an inequality, the numerical evidence of Section~\ref{subsection:asympt-est} indicates that for $1324$ the two quantities may actually be equal (which we conjecture in Conjecture~\ref{conjecture:1324}).

\subsection{Growth rate estimates}
\label{subsection:asympt-est}

For the six involution classes whose enumerations are known exactly, the counting sequences all take the simple power-law form $a_n \sim A\,\mu^n n^g$. Neither $\I(1324)$ nor $\I(4231)$ seems to behave this way. In both, the asymptotic forms seem to carry an additional factor $\mu_1^{n^\sigma}$ with $0 < \sigma < 1$, a \emph{stretched exponential}. Conway, Guttmann, and Zinn-Justin~\cite{conway:4231-50-terms} gave numerical evidence that $\sigma = 1/2$ for permutations avoiding $1324$, and we find the same value for both involution classes, leading to
\[
	a_n \sim B\,\mu^n\cdot \mu_1^{\sqrt{n}}\cdot n^g.
\]

Two caveats apply throughout. First, we are estimating the exponential growth rate $\mu = \lim_{n \to \infty} \sqrt[n]{a_n}$. Proposition~\ref{proposition:sum-indec-gr} guarantees that this limit exists for $\I(4231)$, but in the case of $\I(1324)$ we only conjecture that it exists.

Second, the counting sequences computed in Section~\ref{section:mosaic} are too short for several of these extrapolations to be useful, so we first extend each of them using the series extension method of Guttmann~\cite{guttmann:extension}. This produces $15$ further approximate terms for $\I(1324)$, namely $a_{35}$ through $a_{49}$, and $25$ further approximate terms for $\I(4231)$, namely $a_{41}$ through $a_{65}$. These are expected to be accurate to $10$ significant digits at the first extrapolated term, with the accuracy falling to $5$ or $6$ significant digits at the last. Every plot in this section uses the extended series, although we do not use all the extrapolated coefficients in all the plots. While 5 or 6 significant digit accuracy is sufficient for simple extrapolation methods such as simple ratios, some of our extrapolation methods intrinsically involve first- second- or even third-order differences, which require increasingly precise accuracy in the coefficient estimates. Accordingly, we use fewer approximate coefficients, (only the most precise), in such cases. We indicate in the text that follows when this is the case.

Both sequences show a period-$2$ oscillation in the ratios of consecutive coefficients, so throughout we work with
\[
	r_n = \sqrt{a_n/a_{n-2}}
\]
in place of $a_n/a_{n-1}$, which damps the oscillation. To avoid the noise of early terms, we  ignore the first 10 or 12 terms in our various extrapolations.

\subsubsection{Involutions avoiding 1324}

Let $a_n$ denote the number of $1324$-avoiding involutions, known exactly up to $a_{34}$ and, after the extension described above, approximately up to $a_{49}$. We now use these terms to empirically estimate the asymptotic behavior of the sequence.

\begin{figure}
    \centering
    \begin{minipage}[t]{0.48\linewidth}
        \centering
        \includegraphics[width=\linewidth]{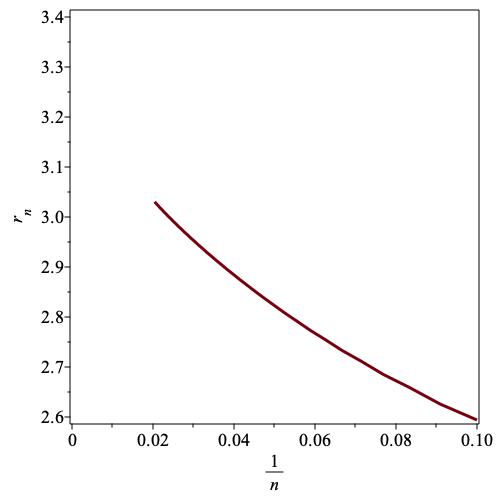}
        \caption{Ratio plot vs. $1/n$ for $\I(1324)$.}
        \label{figure:1324-ratios-n}
    \end{minipage}\hfill
    \begin{minipage}[t]{0.48\linewidth}
        \centering
        \includegraphics[width=\linewidth]{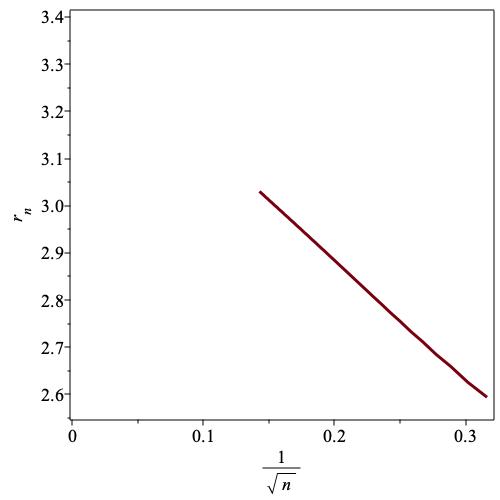}
        \caption{Ratio plot vs. $1/\sqrt{n}$ for $\I(1324)$.}
        \label{figure:1324-ratios-sqrtn}
    \end{minipage}
\end{figure}

A simple ratio plot of the coefficients, in which we plot the ratios $r_n$ against $1/n,$ is shown in Figure~\ref{figure:1324-ratios-n}, and displays significant curvature. This behavior of the ratios is reminiscent of that observed for $\Av(1324)$~\cite{conway:4231-50-terms} where it was argued that such behavior is the hallmark of a stretched-exponential singularity, $\mu_1^{n^\sigma}$.

If $\sigma=1/2,$ the ratios should approach linearity when plotted against $1/\sqrt{n}$. We show this plot in Figure~\ref{figure:1324-ratios-sqrtn}. Of course, we expect that there will still be terms ${\rm O}(1/n)$ in the expression for the ratios, and this is manifested by a residual mild degree of curvature in this plot.

From Figure~\ref{figure:1324-ratios-sqrtn} the ratios appear to be approaching a limit around 3.4. A quadratic fit through the data points gives the extrapolated limit as 3.4128, while a linear fit gives the limit as 3.3948. These are very close to the conjectured lower bound 3.4056.

To support our assumption that $\sigma=1/2,$ we will try to estimate the value of $\sigma$ directly from the data. From the conjectured expression for the coefficients, it follows that the ratios should behave as
\begin{equation}\label{equation:ratio-expansion}
	r_n \sim \mu \left ( 1 +\frac{\sigma \log {\mu_1}}{n^{1-\sigma}} +\frac{g}{n} + {\rm O} \left (\frac{1}{n^{2-2\sigma}}\right ) \right ).
\end{equation}

We estimate the value of $\sigma$ by first observing that
\[
	s_n = \log \left (1 -\frac{r_n}{\mu} \right ) \sim \text{const.} + (\sigma - 1)\log{n},
\]
so that $\sigma - 1 \sim \frac{s_n - s_{n-2}}{\log{n}-\log(n-2)}$.

\begin{figure}
    \centering
    \includegraphics[width=0.48\linewidth]{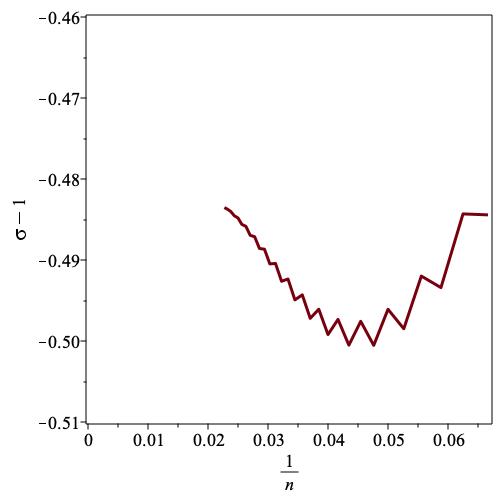}
    \caption{Estimate of the exponent $\sigma-1$ vs. $1/n$ for $\I(1324)$.}
    \label{figure:1324-sigma}
\end{figure}

In Figure~\ref{figure:1324-sigma} we show a plot of $\sigma-1$ estimates against $1/n$, where we have used the conjectured lower bound $\mu=3.4056$ in calculating $s_n.$ This plot can be seen to be approaching a value slightly above $\sigma = 1/2,$ but there is a hint of curvature that suggests it may be approaching a maximum. (We only used the first nine extrapolated coefficients in this plot).

The closeness of our estimate of $\mu$ to the conjectured lower bound suggests it is worth investigating whether they are indeed equal (or as close to equal as numerical tests can confirm).

\begin{figure}
    \centering
    \begin{minipage}[t]{0.48\linewidth}
        \centering
        \includegraphics[width=\linewidth]{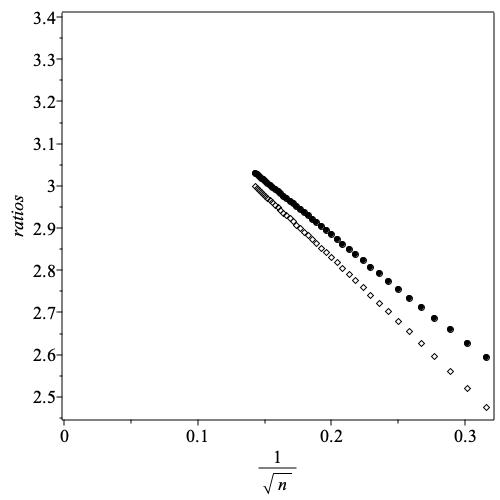}
        \caption{Upper curve: square root of ratios of involution coefficients. Lower curve: fourth root of ratios of the coefficients of $\Av(1324)$.}
        \label{figure:1324-overlay}
    \end{minipage}\hfill
    \begin{minipage}[t]{0.48\linewidth}
        \centering
        \includegraphics[width=\linewidth]{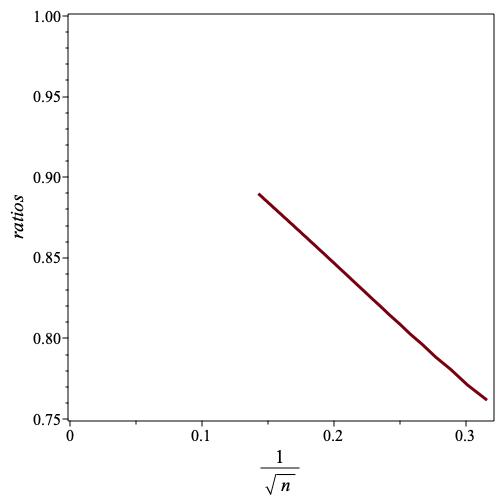}
        \caption{Estimate of the rescaled ratios vs. $1/\sqrt{n}$ for $\I(1324)$.}
        \label{figure:1324-normalized}
    \end{minipage}
\end{figure}

In Figure~\ref{figure:1324-overlay} we show two point plots. Let $b_n$ denote the number of $1324$-avoiding permutations. The upper curve (solid circles) consists of the points $r_n = \sqrt{a_n/a_{n-2}}$ plotted against $1/\sqrt{n},$ whereas the lower curve (open diamonds) consists of the points given by $\left(b_n/b_{n-2}\right)^{1/4},$ also plotted against $1/\sqrt{n}.$ If the two curves have the same limit as $n \to \infty,$ then the involution growth constant estimate is indeed equal to the lower bound. This seems plausible from the figure.

To investigate this further, we tentatively assume that the growth constant for $\I(1324)$ is $\mu=\sqrt{9+3\sqrt{3}/2}.$ Then we re-scale the ratio plot shown in Figure~\ref{figure:1324-ratios-sqrtn} by dividing the ratios by $\mu.$
A plot of the rescaled ratios against $1/\sqrt{n}$ is shown in Figure~\ref{figure:1324-normalized}, and it is clearly plausible that the limit is exactly 1. Indeed, a linear fit to the data points gives the limit as 0.9968, while a quadratic fit gives the limit as 1.0021.

\begin{figure}
    \centering
    \begin{minipage}[t]{0.48\linewidth}
        \centering
        \includegraphics[width=\linewidth]{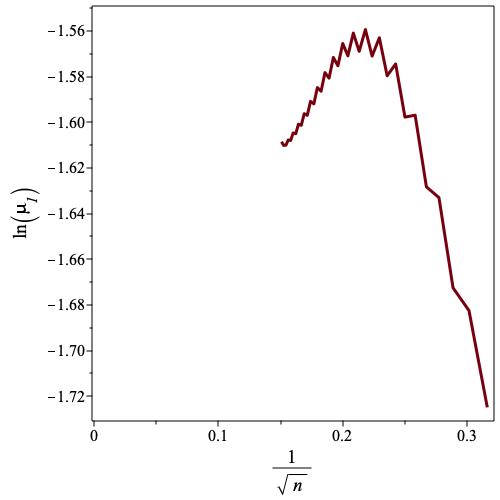}
        \caption{Estimate of $\log{\mu_1}$ vs. $1/\sqrt{n}$ for $\I(1324)$.}
        \label{figure:1324-mu1}
    \end{minipage}\hfill
    \begin{minipage}[t]{0.48\linewidth}
        \centering
        \includegraphics[width=\linewidth]{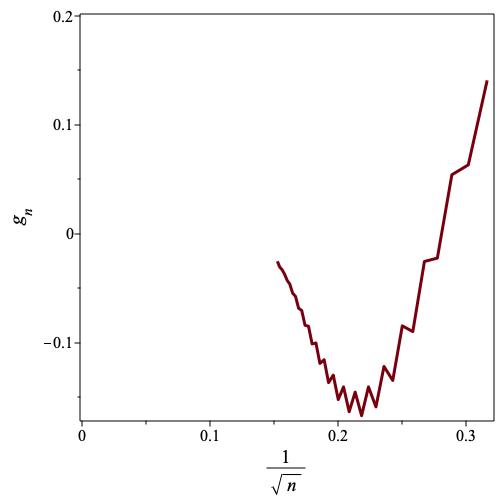}
        \caption{Estimate of exponent $g$ vs. $1/\sqrt{n}$ for $\I(1324)$.}
        \label{figure:1324-g}
    \end{minipage}
\end{figure}

Given that we expect that the coefficients behave as
\[
	a_n \sim B\,\mu^n\cdot \mu_1^{\sqrt{n}}\cdot n^g,
\]
we attempted to estimate the sub-exponential term $\mu_1,$ and exponent $g$ assuming the conjectured value for the growth constant. Note from the above equation that
\[
	\lambda_n \equiv \frac{\log{\tilde a}_n}{\sqrt{n}} \sim \log{\mu_1} + \frac{g\log{n}}{\sqrt{n}}+\frac{\log{B}}{\sqrt{n}}.
\]
We then solve the linear system given by the triple $\lambda_{k-2},\,\,\lambda_k,\,\,\lambda_{k+2}$ for the three unknowns $\log{\mu_1},$ $g,$ and $\log{B},$ increasing $k$ until we run out of known coefficients. We show a plot of successive estimates of $\log{\mu_1}$ against $1/\sqrt{n}$ in Figure~\ref{figure:1324-mu1} and of $g$ in Figure~\ref{figure:1324-g}. (In these plots we used only the first nine of the fifteen approximate coefficients). From these figures, any extrapolation is arguably brave or foolish or both, but we very tentatively suggest a value around $\log{\mu_1} \approx -1.67$ (hence $\mu_1 \approx 0.19$) and even more tentatively that $g$ takes a small positive value,  assuming the observed trends continue.

We have not even attempted an estimate of $\log{B}.$ We remark that for $\Av(1324)$, the corresponding estimates are $\mu_1 \approx 0.040$ and $g \approx -1.25.$\footnote{In~\cite{conway:4231-50-terms} the estimate $g = -1.1 \pm 0.1$ is given. Subsequent analysis leads us to believe that $g$ is closer to $-5/4.$} We have suggested that the growth constants are possibly related by $\gr(\I(1324))^2 = \gr(\Av(1324)).$ From the numerical evidence it is also possible that a similar relationship holds for the growth constant of the stretched exponential term: if the value of $\mu_1$ for $\I(1324)$ were the square root of its value for $\Av(1324)$, we would have $\mu_1 \approx 0.200$ for the involutions, quite close to our direct estimate $\mu_1 \approx 0.19.$

Let us investigate this possibility further. If, for $\Av(1324),$ the coefficients behave as
\[
    b_n \sim A\mu^n \mu_1^{\sqrt{n}} n^{g_1},
\]
and, for $\I(1324),$ the coefficients behave as
\[
    a_n \sim B\mu_I^n \mu_{1,I}^{\sqrt{n}} n^{g_2},
\]
and if $\mu_I=\sqrt{\mu}$ and $\mu_{1,I}=\sqrt{\mu_1},$ then
\[
    c_n \equiv \sqrt{b_n}/a_n \sim Cn^\theta,
\]
where $C=\sqrt{A}/B$ and $\theta=g_1/2-g_2.$ If the first ``if'' is violated, so that the growth constants are not related as proposed, then $c_n$ will diverge exponentially, rather than as a power-law. If the second ``if'' is violated, so that the stretched-exponential growth constants are not related as proposed, then $c_n$ will diverge sub-exponentially. The proposed power-law divergence will only hold if both conditions are satisfied. And in that case the ratios $r_n=c_n/c_{n-1}$ will extrapolate to 1 linearly when plotted against $1/n.$ If the first ``if'' is violated, the ratios will go to a value other than 1, and if the second ``if'' is violated, the ratio plot will not be linear, but will show curvature, as we have seen for example in Figure~\ref{figure:1324-ratios-n}.

\begin{figure}
    \centering
    \begin{minipage}[t]{0.48\linewidth}
        \centering
        \includegraphics[width=\linewidth]{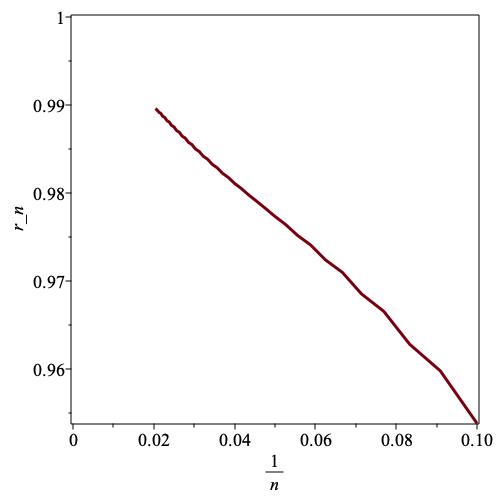}
        \caption{Ratios $\sqrt{c_n/c_{n-2}}$ vs. $1/n.$ }
        \label{figure:c-ratios}
    \end{minipage}\hfill
    \begin{minipage}[t]{0.48\linewidth}
        \centering
        \includegraphics[width=\linewidth]{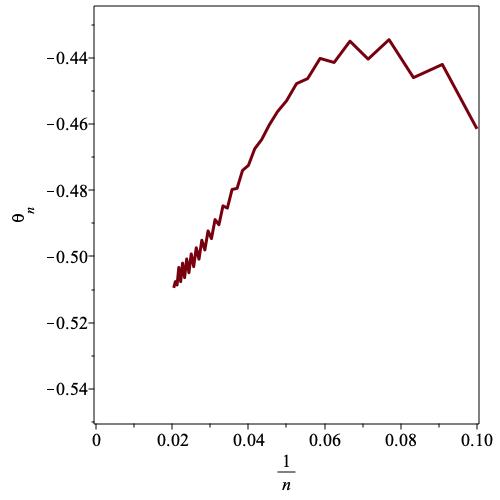}
        \caption{Estimate of exponent $\theta=g_1/2-g_2$ vs. $1/n.$}
        \label{figure:theta}
    \end{minipage}
\end{figure}

In Figure~\ref{figure:c-ratios} we show the ratios $r_n \equiv \sqrt{c_n/c_{n-2}}$ plotted against $1/n.$ This is seen to be convincingly approaching 1, and approaching linearly. This is strong numerical support for the two ``if''s. We next estimate the exponent $\theta,$ as it follows from the proposed power-law behavior that $\theta \sim (r_n-1)\cdot n.$ We define $\theta_n \equiv (r_n-1)\cdot n,$ and show in Figure~\ref{figure:theta} a plot of $\theta_n$ against $1/n.$ This appears to be going to a value around $-0.5$ or a little below. We estimate $\theta = -0.525 \pm 0.025.$ As mentioned above, our best estimate of $g_1$ is -1.25, so this would imply $g_2 \approx -0.1.$ (Note that this differs from our direct estimate above, where we tentatively suggested a small positive value for $g_2.$ This new estimate depends of course on the correctness of the estimate for the exponent $g_1.$)
 
In conclusion, we estimate that the coefficients behave as
\[
	a_n \sim B\,\mu^n\cdot \mu_1^{\sqrt{n}}\cdot n^g,
\]
where both $\mu_1$ and $\mu$ are conjectured to be given by the square root of the corresponding values for $\Av(1324)$, that is to say, around $\mu=3.4056$ and $\mu_1 = 0.200$. For the exponent $g$ we find $g \approx g_1/2+0.525,$ where $g_1$ is the power-law exponent associated with the coefficients of $\Av(1324).$ Our best estimate, $g_1 \approx -1.25$ implies $g \approx -0.1.$ We are unable to give a useful estimate of the amplitude $B$.

It was observed in \cite{bostan:stieltjes} that the coefficients of $\Av(1324)$ appeared to be expressible as a Stieltjes moment sequence, which then allows a quite tight lower bound on the growth constant to be obtained. The coefficients of $\I(1324)$ are not even log-convex (a necessary, but not sufficient condition for a sequence to be a Stieltjes moment sequence).

\subsubsection{Involutions avoiding 4231}

Let $a_n$ denote the number of $4231$-avoiding involutions, known exactly up to $a_{40}$ and, after the extension described above, approximately up to $a_{65}$. We now use these terms to empirically estimate the asymptotic behavior of the sequence.

\begin{figure}
    \centering
    \begin{minipage}[t]{0.48\linewidth}
        \centering
        \includegraphics[width=\linewidth]{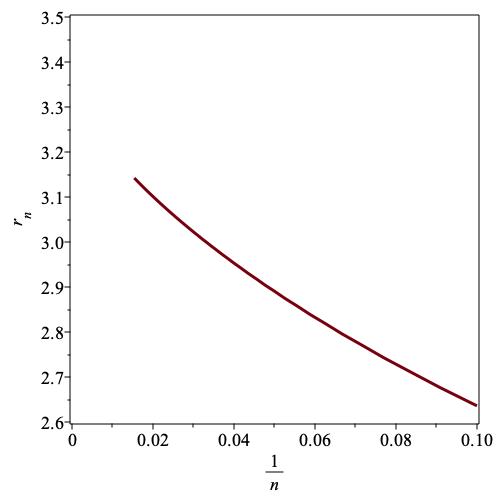}
        \caption{Ratio plot vs. $1/n$ for $\I(4231)$.}
        \label{figure:4231-ratios-n}
    \end{minipage}\hfill
    \begin{minipage}[t]{0.48\linewidth}
        \centering
        \includegraphics[width=\linewidth]{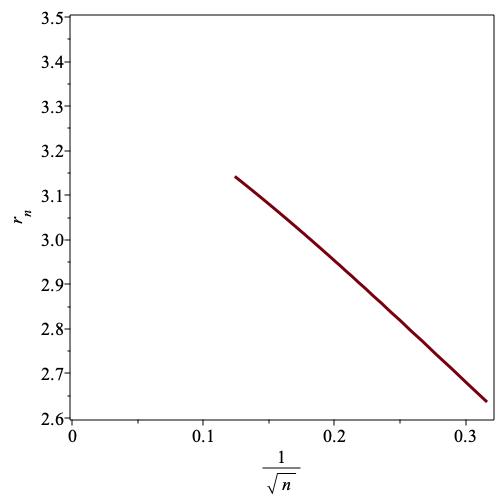}
        \caption{Ratio plot vs. $1/\sqrt{n}$ for $\I(4231)$.}
        \label{figure:4231-ratios-sqrtn}
    \end{minipage}
\end{figure}

As with the previous sequence, the ratios plotted against $1/n,$ shown in Figure~\ref{figure:4231-ratios-n}, display curvature, while the ratios are approximately linear when plotted against $1/\sqrt{n},$ as shown in Figure~\ref{figure:4231-ratios-sqrtn}, though note that there is more curvature than is evident in the corresponding graph for 1324-avoiding involutions. Presumably this is due to more significant higher-order terms. This figure looks very similar to that for $\I(1324)$, though the limits look slightly higher, perhaps 3.45 for the ratios, so our preliminary estimate is $\mu = 3.45 \pm 0.05.$ A linear fit through the displayed points yields $\mu \approx 3.471,$ and a quadratic fit yields $\mu \approx 3.422,$ so all these estimates are self-consistent.

\begin{figure}
    \centering
    \begin{minipage}[t]{0.48\linewidth}
        \centering
        \includegraphics[width=\linewidth]{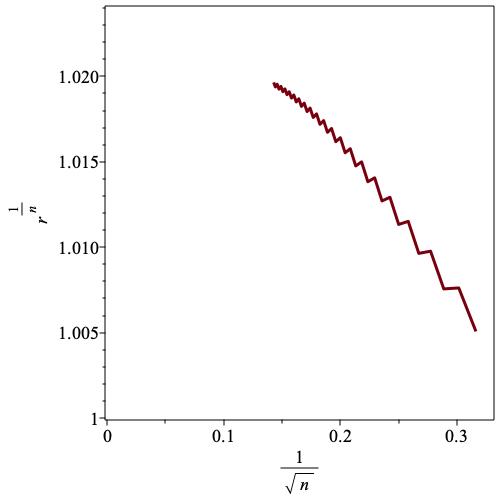}
        \caption{Coefficient ratio $\left(|\I_n(4231)|/|\I_n(1324)|\right)^{1/n}$ vs. $1/\sqrt{n}$.}
        \label{figure:mu-ratio}
    \end{minipage}\hfill
    \begin{minipage}[t]{0.48\linewidth}
        \centering
        \includegraphics[width=\linewidth]{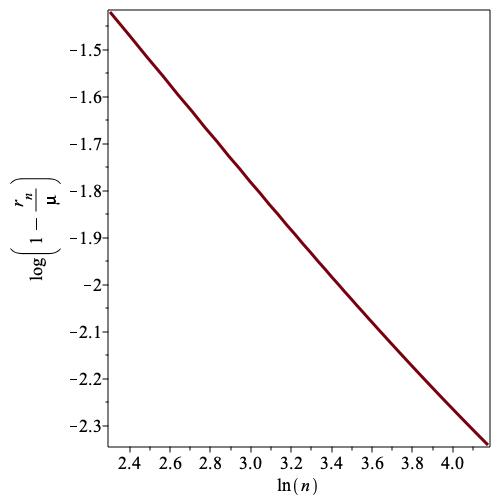}
        \caption{Plot of $\log{(1-r_n/\mu)}$ against $\log{n}$ for $\I(4231)$.}
        \label{figure:4231-sigma}
    \end{minipage}
\end{figure}

In order to more precisely estimate the growth constant, we utilize our estimate of the growth constant for the case of 1324-avoiding involutions just discussed, and form the ratio of the coefficients $\left(\frac{|\I_n(4231)|}{|\I_n(1324)|}\right)^{1/n}$, which should, in the large $n$ limit, approach the ratio of the growth constants $\gr(\I(4231))/\gr(\I(1324)).$ From this plot, shown in Figure~\ref{figure:mu-ratio}, the ratio appears to be approaching a limit of around 1.021, and so we obtain the estimate $\gr(\I(4231)) \approx 3.477,$ which we take as our best estimate.

Using the estimate $\mu \approx 3.477$, we show in Figure~\ref{figure:4231-sigma} a plot of $\log(1-r_n/\mu)$ against $\log{n}.$ From~\eqref{equation:ratio-expansion} this should have gradient $\sigma-1.$ A linear fit gives the gradient as $-0.486$, which is as close to $-1/2$ as we can reasonably expect, so this again is consistent with our assertion that the value of $\sigma$ is likely to be $1/2$ exactly.\footnote{The value of the gradient depends slightly on the range chosen. We have used $n \ge 12$ to avoid low-$n$ curvature.}

\begin{figure}
    \centering
    \begin{minipage}[t]{0.48\linewidth}
        \centering
        \includegraphics[width=\linewidth]{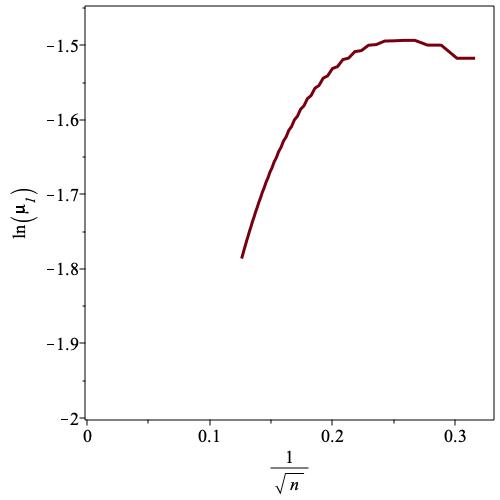}
        \caption{Three-parameter fit estimate of $\log{\mu_1}$ for $\I(4231)$ vs. $1/\sqrt{n}$.}
        \label{figure:4231-mu1}
    \end{minipage}\hfill
    \begin{minipage}[t]{0.48\linewidth}
        \centering
        \includegraphics[width=\linewidth]{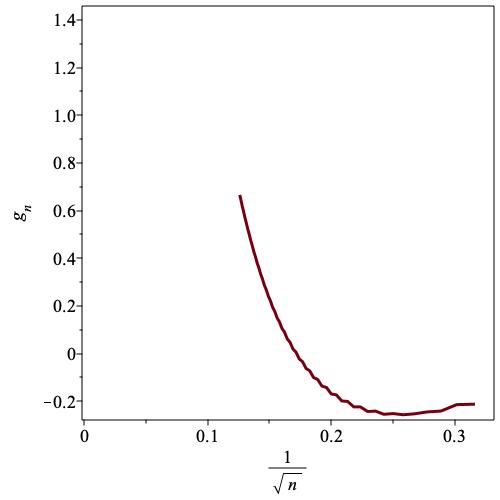}
        \caption{Three-parameter fit estimate of the exponent $g$ for $\I(4231)$ vs. $1/\sqrt{n}$.}
        \label{figure:4231-g}
    \end{minipage}
\end{figure}

Adopting this estimate of the growth constant, we repeated the analysis of the previous subsection to estimate the sub-exponential term $\mu_1$ and the exponent $g$. The relevant plots are shown in Figures~\ref{figure:4231-mu1} and~\ref{figure:4231-g} respectively. These are not reliably extrapolable, and so we are unable to give estimates for $\mu_1$ or $g.$

In conclusion, we estimate that the coefficients behave as
\[
	a_n \sim B\,\mu^n\cdot \mu_1^{\sqrt{n}}\cdot n^g,
\]
where $\mu \approx 3.477.$

Unlike $\I(1324),$ the coefficients of $\I(4231)$ are observed to be log-convex. However, construction of the relevant Hankel determinants shows that the sequence is not a Stieltjes moment sequence.

\section{Concluding remarks}
\label{section:concluding-remarks}

In this paper we adapted Combinatorial Exploration to the involution setting by introducing involution tilings, along with involution-preserving versions of the factor, row placement, row and column separation, and fusion strategies, and we used this machinery to derive the algebraic generating functions of $\I(3421)$ and $\I(2431)$. We also adapted the Mosaic method to count involutions, extending the known counting sequences of $\I(1324)$ and $\I(4231)$ to lengths $34$ and $40$, respectively. These new terms led to conjectured asymptotic forms for both sequences, each involving a stretched exponential factor $\mu_1^{\sqrt{n}}$, and, for $\I(1324)$, to the conjecture that its growth rate is exactly the square root of that of $\Av(1324)$, so that the bound of Proposition~\ref{proposition:skew-indec-gr-bound} is tight.

We return one last time to the ordering of the columns of Table~\ref{table:results}, which are sorted by the number of involutions of length $20$ in each class, replicating the table from~\cite{bona:pattern-avoiding-involutions} whose authors cautioned that this ordering was ``likely still incorrect.'' The growth rates collected in Table~\ref{table:results}, six of which are now known exactly and two of which are conjectured, would imply that for all sufficiently large $n$ the columns belong in the order
\[
	2431, \quad 2341, \quad 1342, \quad 3421, \quad 1234, \quad 2413, \quad 1324, \quad 4231,
\]
corresponding to the growth rates
\[
	2.53041 < 2.53999 < 2.61803 < 2.65897 < 3 < 3.14626 < 3.40559 < 3.477.
\]

\begin{figure}
	\centering
	\begin{footnotesize}
	\begin{tabular}{ccc}
	\begin{tikzpicture}[y=1.25cm, x=.15cm, font=\sffamily, scale=.8]
		\draw (0,0) -- (34,0);
		\draw (0,0) -- (0,4);
		\foreach \x in {0,1,...,34}
			\draw (\x,1pt) -- (\x,-3pt);
		\foreach \x in {0,5,...,30}
			\draw (\x,1pt) -- (\x,-3pt) node[anchor=north] {\x};
		\foreach \y in {0,1,...,4}
			\draw (1pt,\y) -- (-3pt,\y) node[anchor=east] {\y};
		\node[right=6pt] at (34,3.62) {$\displaystyle\frac{|\I_{n}(2413)|}{|\I_{n}(1324)|}$};
		\draw[gray] (0,1) -- (34,1);
		\draw plot[mark=*, mark size=.05cm] coordinates {
			(0,1.0000) (1,1.0000) (2,1.0000) (3,1.0000) (4,1.1111) (5,1.1429)
			(6,1.2549) (7,1.3175) (8,1.4206) (9,1.5049) (10,1.5991) (11,1.6979)
			(12,1.7840) (13,1.8924) (14,1.9713) (15,2.0858) (16,2.1580) (17,2.2764)
			(18,2.3424) (19,2.4629) (20,2.5232) (21,2.6445) (22,2.6993) (23,2.8204)
			(24,2.8699) (25,2.9900) (26,3.0344) (27,3.1529) (28,3.1922) (29,3.3084)
			(30,3.3427) (31,3.4561) (32,3.4854) (33,3.5957) (34,3.6200)
		};
	\end{tikzpicture}
	&&
	\begin{tikzpicture}[y=2.273cm, x=.15cm, font=\sffamily, scale=.8]
		\draw (0,0) -- (40,0);
		\draw (0,0) -- (0,2.2);
		\foreach \x in {0,1,...,40}
			\draw (\x,1pt) -- (\x,-3pt);
		\foreach \x in {0,5,...,40}
			\draw (\x,1pt) -- (\x,-3pt) node[anchor=north] {\x};
		\foreach \y in {0,1,2}
			\draw (1pt,\y) -- (-3pt,\y) node[anchor=east] {\y};
		\node[right=6pt] at (40,1.89) {$\displaystyle\frac{|\I_{n}(2413)|}{|\I_{n}(4231)|}$};
		\draw[gray] (0,1) -- (40,1);
		\draw plot[mark=*, mark size=.05cm] coordinates {
			(0,1.0000) (1,1.0000) (2,1.0000) (3,1.0000) (4,1.1111) (5,1.1429)
			(6,1.2549) (7,1.2969) (8,1.3945) (9,1.4382) (10,1.5202) (11,1.5620)
			(12,1.6299) (13,1.6679) (14,1.7236) (15,1.7568) (16,1.8020) (17,1.8301)
			(18,1.8662) (19,1.8890) (20,1.9171) (21,1.9347) (22,1.9556) (23,1.9683)
			(24,1.9828) (25,1.9907) (26,1.9995) (27,2.0028) (28,2.0065) (29,2.0056)
			(30,2.0047) (31,2.0000) (32,1.9948) (33,1.9866) (34,1.9778) (35,1.9664)
			(36,1.9543) (37,1.9400) (38,1.9250) (39,1.9082) (40,1.8906)
		};
	\end{tikzpicture}
	\end{tabular}
	\end{footnotesize}
	\caption{The ratio of the number of $2413$-avoiding involutions of length $n$ to the number of $1324$-avoiding involutions of length $n$ for $n \leq 34$ (left) and to the number of $4231$-avoiding involutions of length $n$ for $n \leq 40$ (right).}
	\label{figure:ratio-plots}
\end{figure}
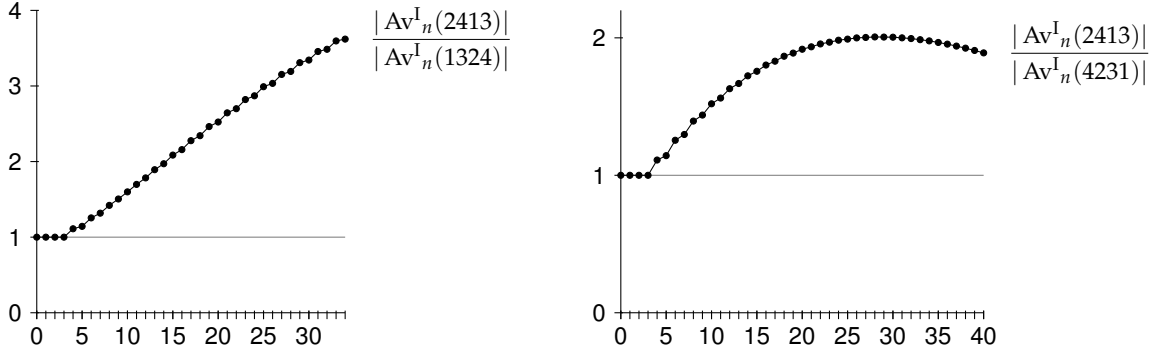

How large is sufficiently large for the columns to be in the correct order? The concluding section of~\cite{bona:pattern-avoiding-involutions} plotted the ratio $|\I_n(2413)|/|\I_n(1324)|$ for $n \leq 25$, observing that while this ratio should tend to $0$, the initial terms ``do not paint a very convincing picture of a sequence going to $0$.'' We can now replicate this plot for $n \leq 34$, shown on the left of Figure~\ref{figure:ratio-plots}, and although the curve still does not look like one that approaches $0$, it cannot keep climbing: the lower bound $\gr(\Av(1324)) > 10.27$ of Bevan, Brignall, Elvey Price, and Pantone~\cite{bevan:1324-staircase} combines with Proposition~\ref{proposition:skew-indec-gr-bound} to give
\[
	\ugr(\I(1324)) > 3.20 > \sqrt{2} + \sqrt{3} = \gr(\I(2413)).
\]
The ratio therefore falls below any positive threshold infinitely often, and tends to $0$ if $\I(1324)$ has a proper growth rate, as Conjecture~\ref{conjecture:1324} asserts.

The right plot of Figure~\ref{figure:ratio-plots} shows the corresponding ratio with $4231$ in place of $1324$, which we can plot all the way to $n = 40$. Here we cannot prove that the ratio tends to $0$, but the conjectured growth rate of $\I(4231)$ comfortably exceeds the known growth rate of $\I(2413)$, so it seems likely to be true. This ratio does turn around within the range we can compute, peaking just above $2$ at $n = 28$ and decreasing steadily thereafter. Extrapolating the ratios suggests it may not drop below $1$ until $n$ is between $70$ and $80$. The plot on the left may not drop below $1$ until $n$ is well into the hundreds.

Lastly, based on our numerical estimates, we pose the following.
\begin{conjecture}
\label{conjecture:gr-dominance}
$\gr(\I(4231)) > \ugr(\I(1324))$.
\end{conjecture}

A more tractable goal may be the following term-by-term comparison, which would follow from an injection from $\I_n(1324)$ into $\I_n(4231)$ and would imply a weak form of Conjecture~\ref{conjecture:gr-dominance}.

\begin{conjecture}
\label{conjecture:term-dominance}
For all $n \geq 0$, we have $|\I_n(4231)| \geq |\I_n(1324)|$.
\end{conjecture}

\section*{Acknowledgments}

JP's research was supported by grant 713579 from the Simons Foundation. The computations in Section~\ref{section:mosaic} were performed on the Raj cluster at Marquette University, funded in part by National Science Foundation award CNS-1828649, ``MRI: Acquisition of iMARC: High Performance Computing for STEM Research and Education in Southeast Wisconsin.''

\bibliographystyle{alpha}
\bibliography{refs}

\end{document}